\documentclass[12pt]{article}
\usepackage{amsfonts}
\usepackage{full page,
arcs, amssymb, amscd, amsmath, graphicx, %showkeys
wrapfig, stmaryrd, hyperref}
\allowdisplaybreaks

\newtheorem{thm}{Theorem}[section]
\newtheorem{defn}[thm]{Definition}
\newtheorem{prop}[thm]{Proposition}
\newtheorem{cor}[thm]{Corollary}
\newtheorem{lemma}[thm]{Lemma}
\newtheorem{rema}[thm]{Remark}

\newcommand{\halmos}{\rule{1ex}{1.4ex}}

\newcommand{\nn}{\nonumber \\}

 \newcommand{\res}{\mbox{\rm Res}}

\renewcommand{\hom}{\mbox{\rm Hom}}

 \newcommand{\pf}{{\it Proof.}\hspace{2ex}}
 
 \newcommand{\epfv}{\hspace*{\fill}\mbox{$\halmos$}\vspace{1em}}

\newcommand{\wt}{\mbox{\rm wt}\,}
\newcommand{\swt}{\mbox{\rm {\scriptsize wt}}\,}

\newcommand{\mbar}{\Big\vert}
\newcommand{\lbar}{\bigg\vert}

\newcommand{\Y}{\mathcal{Y}}
\newcommand{\C}{\mathbb{C}}
\newcommand{\Z}{\mathbb{Z}}
\newcommand{\R}{\mathbb{R}}
\newcommand{\Q}{\mathbb{Q}}
\newcommand{\N}{\mathbb{N}}

\newcommand{\one}{\mathbf{1}}

\renewcommand{\l}{\llfloor}
\renewcommand{\r}{\rrfloor}

\makeatletter
\DeclareFontFamily{U}{tipa}{}
\DeclareFontShape{U}{tipa}{m}{n}{<->tipa10}{}
\newcommand{\arc@char}{{\usefont{U}{tipa}{m}{n}\symbol{62}}}%

\newcommand{\arc}[1]{\mathpalette\arc@arc{#1}}

\newcommand{\arc@arc}[2]{%
  \sbox0{$\m@th#1#2$}%
  \vbox{
    \hbox{\resizebox{\wd0}{\height}{\arc@char}}
    \nointerlineskip
    \box0
  }%
}
\makeatother

\title{ {\bf  Associative algebras
and  intertwining operators} }
\date{}
\author{Yi-Zhi Huang}
\begin{document}

\bibliographystyle{alpha}
\maketitle
\begin{abstract}
Let $V$ be a vertex operator algebra and  $A^{\infty}(V)$ and $A^{N}(V)$ for $N\in \N$ the associative algebras
introduced by the author in \cite{H-aa-va}. For a lower-bounded generalized $V$-module $W$,
we give $W$ a structure of graded $A^{\infty}(V)$-module and 
we introduce an $A^{\infty}(V)$-bimodule $A^{\infty}(W)$ and an $A^{N}(V)$-bimodule $A^{N}(W)$. 
We prove that the space of (logarithmic)  intertwining operators  of type $\binom{W_{3}}{W_{1}W_{2}}$
for lower-bounded generalized $V$-modules $W_{1}$, $W_{2}$ and $W_{3}$
 is isomorphic to the space $\hom_{A^{\infty}(V)}(A^{\infty}(W_{1})\otimes_{A^{\infty}(V)}W_{2}, W_{3})$.
Assuming that $W_{2}$ and $W_{3}'$ are equivalent to certain universal lower-bounded 
generalized $V$-modules generated by their $A^{N}(V)$-submodules consisting of elements 
of levels less than or equal to $N\in \N$, we also prove that the space 
of (logarithmic) intertwining operators of type $\binom{W_{3}}{W_{1}W_{2}}$ 
is isomorphic to the space of $\hom_{A^{N}(V)}(A^{N}(W_{1})\otimes_{A^{N}(V)}\Omega_{N}^{0}(W_{2}), \Omega_{N}^{0}(W_{3}))$.
\end{abstract}

\renewcommand{\theequation}{\thesection.\arabic{equation}}
\setcounter{equation}{0}
\setcounter{thm}{0}
\section{Introduction}

In \cite{H-aa-va}, for a grading-restricted vertex algebra $V$, 
the author introduced associative algebras $A^{\infty}(V)$ and $A^{N}(V)$ 
for $N\in \N$ and proved basic properties of these algebras and their modules
in connection with $V$ and lower-bounded generalized  $V$-modules. 
In the present paper, we study the connections between (logarithmic) intertwining operators 
among lower-bounded generalized $V$-modules and module maps between 
suitable modules for these associative algebras in the case that $V$ is a vertex operator algebra
so that there is a conformal vector $\omega\in V$.
The results of the present paper will be used in the last step of a proof of the 
modular invariance of (logarithmic) intertwining operators in the case that 
$V$ satisfies the positive energy condition (or CFT type) and $C_{2}$-cofiniteness
condition in a paper in preparation.
For simplicity, we shall omit  ``(logarithmic)" in ``(logarithmic) intertwining operator" in this paper so that 
by an intertwining operator,
we always mean an intertwining operator which might contains logarithm of the variable.

In the case $N=0$, the algebra $A^{0}(V)$ was proved in \cite{H-aa-va}
 to be isomorphic to the Zhu algebra $A(V)$ (see \cite{Z}). In this case, 
there is a theorem of Frenkel and Zhu \cite{FZ}  stating that 
for irreducible $V$-modules $W_{1}$, $W_{2}$ and $W_{3}$, the space of intertwining operators
of type $\binom{W_{3}}{W_{1}W_{2}}$ is linearly isomorphic to 
$\hom_{A(V)}(A(W_{1})\otimes_{A(V)}\Omega(W_{2}), \Omega(W_{3}))$, where 
$\Omega(W_{2})$ and $\Omega(W_{3})$ are the lowest weight spaces of $W_{2}$ and $W_{3}$, respectively, 
and $A(W_{1})$ is an $A(V)$-bimodule introduced in the paper \cite{FZ}.
A proof of this theorem was given by Li in \cite{L} under the assumption 
that every lower-bounded generalized $V$-module (or $\N$-gradable weak $V$-module)
 is completely reducible. In the same paper \cite{L}, 
a counterexample was also given to show that without this assumption, 
the theorem is not true. 
Since this result needs this semisimplicity assumption,  
it cannot be used to study intertwining operators 
in the case that lower-bounded generalized $V$-modules might not be completely reducible. 
In \cite{HY1} and \cite{HY2}, under strong assumptions on $W_{2}$ and $W_{3}'$ but without 
this semisimplicity assumption, 
Yang and the author proved that the space of intertwining operators
of type $\binom{W_{3}}{W_{1}W_{2}}$ is linearly isomorphic to 
$\hom_{A_{N}(V)}(A_{N}(W_{1})\otimes_{A_{N}(V)}\Omega_{N}^{0}(W_{2}), \Omega_{N}^{0}(W_{3}))$,
where $\Omega_{N}^{0}(W_{2})$ and $\Omega_{N}^{0}(W_{3})$ are suitable subspaces of
$W_{2}$ and $W_{3}$, respectively,  and $A_{N}(W_{1})$ is a bimodule for the generalization $A_{N}(V)$ of 
$A(V)$ by Dong, Li and Mason \cite{DLM}.

In the present paper, for a lower-bounded generalized $V$-module $W$,
we introduce an $A^{\infty}(V)$-bimodule $A^{\infty}(W)$ and an $A^{N}(V)$-bimodule $A^{N}(W)$ for each $N\in \N$. 
We prove that the space of (logarithmic)  intertwining operators of type $\binom{W_{3}}{W_{1}W_{2}}$
for lower-bounded generalized $V$-modules $W_{1}$, $W_{2}$ and $W_{3}$
 is isomorphic to the space $\hom_{A^{\infty}(V)}(A^{\infty}(W)\otimes_{A^{\infty}(V)}W_{2}, W_{3})$. 
Assuming that $W_{2}$ and $W_{3}'$ are certain universal lower-bounded generated by their $A^{N}(V)$-submodules consisting of elements 
of levels less than or equal to $N\in \N$, we also prove that the space of 
(logarithmic)  intertwining operators  of type $\binom{W_{3}}{W_{1}W_{2}}$
is isomorphic to the space of $\hom_{A^{N}(V)}(A^{N}(W)\otimes_{A^{N}(V)}W_{2}, W_{3})$. 

Here we give more discussions on the main results. 
In \cite{H-aa-va}, the associative algebra $A^{\infty}(V)$ is defined using the associated graded 
spaces given by suitable ascending filtrations of lower-bounded generalized $V$-modules for a grading-restricted vertex algebra. 
These associated graded spaces are by definition nondegenerate graded $A^{\infty}(V)$-modules (see \cite{H-aa-va}). 
But it is not easy to work with the associated graded space of a lower-bounded module, since
in general it is very difficult to determine the associated graded space explicitly. 
In this paper we shall work with a canonical $\N$-grading 
$W=\coprod_{n\in \N}W_{\l n\r}$ of a lower-bounded generalized $V$-module $W$
and show that $W$ itself also has a structure of 
graded $A^{\infty}(V)$-module. We construct the structure of 
graded $A^{\infty}(V)$-module on $W$ by showing that $W$  is 
in fact a quotient of the associated graded space of another lower-bounded generalized $V$-module. 
But in general this structure of a graded $A^{\infty}(V)$-module on $W$ might not be nondegenerate
(see \cite{H-aa-va}). This result holds for a grading-restricted vertex algebra $V$ which might not have a conformal vector.

Given a lower-bounded generalized $V$-module $W$, let $U^{\infty}(W)$ be  
the space of all column-finite matrices with entries 
in $W$ and doubly index by $\N$.
We introduce  left and right actions of $A^{\infty}(V)$ on $U^{\infty}(W)$. 
Let $W_{2}$ and $W_{3}$ be lower-bounded generalized $V$-modules
and let $\Y$ be
an intertwining operator of type $\binom{W_{3}}{WW_{2}}$.
We introduce a linear map $\vartheta_{\Y}: U^{\infty}(W)\to \hom(W_{2}, W_{3})$.
Note that $A^{\infty}(V)$ is a quotient of a nonassociative algebra $U^{\infty}(V)$
of column-finite matrices with entries 
in $V$ and doubly index by $\N$.
Since $W_{2}$ and $W_{3}$ are left $A^{\infty}(V)$-modules, 
$\hom(W_{2}, W_{3})$ has a natural structure of $A^{\infty}(V)$-bimodule and in particular, 
a natural structure of $U^{\infty}(V)$-bimodule. We prove that 
$\vartheta_{\Y}$ commutes with the left and right actions of $A^{\infty}(V)$. 
Let $Q^{\infty}(W)$ be the intersection of $\ker \vartheta_{\Y}$ for 
all such lower-bounded generalized $V$-modules $W_{2}$ and $W_{3}$ and 
intertwining operators $\Y$ of type $\binom{W_{3}}{WW_{2}}$. Then 
$A^{\infty}(W)=U^{\infty}(W)/Q^{\infty}(W)$ is an $A^{\infty}(V)$-bimodule. 
The construction of the $A^{\infty}(V)$-bimodule $A^{\infty}(W)$ works for 
a grading-restricted vertex algebra $V$ which might not have a conformal vector. 

Given lower-bounded generalized $V$-modules $W_{1}$, $W_{2}$, $W_{3}$
 and an intertwining
operator $\Y$ of type $\binom{W_{3}}{W_{1}W_{2}}$, we prove that 
the map $\vartheta_{\Y}$ discussed above in the case $W=W_{1}$ gives an $A^{\infty}(V)$-module 
map $\rho(\Y): A^{\infty}(W_{1})\otimes_{A^{\infty}(V)} W_{2}\to W_{3}$.
In particular, we obtain a linear map $\rho: \mathcal{V}_{W_{1}W_{2}}^{W_{3}}\to 
\hom_{A^{\infty}(V)}(A^{\infty}(W_{1})\otimes_{A^{\infty}(V)} W_{2}, W_{3})$,
where $\mathcal{V}_{W_{1}W_{2}}^{W_{3}}$ is the space of 
intertwining operators of type $\binom{W_{3}}{W_{1}W_{2}}$.
Our first main theorem states that $\rho$ is an isomorphism. In the proof of 
the first main theorem, we need to assume that $V$ is a vertex operator algebra
which has a conformal vector. This first main theorem can be generalized to 
a grading-restricted vertex algebra. But in this more 
general case, since there is no conformal vector, one has to introduce additional structures 
so that left actions of $L(-1)$ on the $A^{\infty}(V)$-bimodules 
can be introduced. 

Let $N\in \N$ and $W$ a lower-bounded generalized $V$-module.  Let 
$\Omega_{N}^{0}(W)=\coprod_{n=0}^{N}W_{\l n\r}$. Then $\Omega_{N}^{0}(W)$ 
is in fact an $A^{N}(V)$-module. 
We consider the space $U^{N}(W)$ of  $(N+1)\times (N+1)$ matrices with entries in $W$. 
Then $U^{N}(W)$ can be viewed as a subspace of $U^{\infty}(W)$. 
In particular, the cosets in $A^{\infty}(W)$ containing elements of 
$U^{N}(W)$ 
form an $A^{N}(V)$-bimodule $A^{N}(W)$. The construction of the $A^{N}(V)$-bimodule $A^{N}(W)$ works for 
a grading-restricted vertex algebra $V$. 

Given lower-bounded generalized $V$-modules $W_{1}$, $W_{2}$, $W_{3}$ and
an intertwining operator $\Y$ of type $\binom{W_{3}}{W_{1}W_{2}}$, the 
$A^{\infty}(V)$-module  map $\rho(\Y): A^{\infty}(W_{1})\otimes_{A^{\infty}(V)} W_{2}\to W_{3}$
induces an $A^{N}(V)$-module  map 
$\rho^{N}(\Y): A^{N}(W_{1})\otimes_{A^{N}(V)}\Omega_{N}^{0}(W_{2})\to \Omega_{N}^{0}(W_{3})$. 
Then we obtain a linear map $\rho^{N}: \mathcal{V}_{W_{1}W_{2}}^{W_{3}}\to 
\hom_{A^{N}(V)}(A^{N}(W_{1})\otimes_{A^{N}(V)}\Omega_{N}^{0}(W_{2}), 
\Omega_{N}^{0}(W_{3})$. We prove that $\rho^{N}$ is injective. 
Our second main theorem states that $\rho^{N}$ is an isomorphism when $W_{2}$ and $W_{3}'$ are certain
universal lower-bounded generalized $V$-mdoules generated by the $A^{N}(V)$-modules 
$\Omega_{N}^{0}(W_{2})$ and $\Omega_{N}^{0}(W_{3}')$, respectively. 
This second main theorem is proved using the first main theorem above. In particular, 
$V$ is also assumed to be a vertex operator algebra. Again,
this theorem can be generalized to the case that $V$ is a grading-restricted vertex algebra. 

The main motivation of this paper is the modular invariance of intertwining operators
in the nonsemisimple (or logarithmic) case. 
In \cite{H-modular}, the author proved the conjecture of Moore and Seiberg on 
modular invariance of intertwining operators for rational conformal field theories
 (see \cite{MS}). Mathematically this modular invariance is for intertwining operators 
among modules for a vertex operator algebra
$V$ satisfying the conditions that  $V$ is of positive energy (or CFT-type), $V$ is $C_{2}$-cofinite and every 
lower-bounded generalized  $V$-module (or every $\N$-gradable weak $V$-module) 
is completely reducible. 
After the convergence and analytic extensions 
of shifted $q$-traces of products of intertwining operators and the genus-one
associativity is proved, the proof of the modular invariance is reduced to the 
proof that a genus-one one-point correlation 
function can always be written as the analytic extension of a shifted $q$-trace of an intertwining operator. 
It is in this last step of the proof that
the theorem of Frenkel and Zhu mentioned above is used.  

In \cite{F1} and \cite{F}, Fiordalisi proved the convergence and analytic extension property
of shifted pseudo-$q$-traces of products of  intertwining operators and the genus-one
associativity under the condition that $V$ is of positive energy (or CFT-type) and $V$ is $C_{2}$-cofinite, 
but without assuming the condition that 
every lower-bounded generalized  $V$-module (or every $\N$-gradable weak $V$-module) 
is completely reducible. The proof of the modular invariance in this nonsemisimple (or logarithmic)
case is then also reduced to the proof that a genus-one one-point correlation 
function can always be written as a shifted pseudo-$q$-trace of an intertwining operator. 
But in this case, the theorem of Frenkel and Zhu mentioned above cannot be used. 
In a paper  in preparation, we shall give this last step using the results obtained in the present paper.

This paper is organized as follows: We recall in Section 2
the associative algebras 
$A^{\infty}(V)$ and $A^{N}(V)$ and 
graded $A^{\infty}(V)$- and $A^{N}(V)$-modules introduced in \cite{H-aa-va}. 
In Section 3, we give an $A^{\infty}(V)$-module structure
to a lower-bounded generalized $V$-module $W$ . In Section 4, we construct the $A^{\infty}(V)$-bimodule 
$A^{\infty}(W)$ from $W$. Our first main theorem stating that $\rho$ is an isomorphism
is formulated and proved in Section 5. Our second main theorem stating that $\rho^{N}$
is an isomorphism when when $W_{2}$ and $W_{3}'$ are certain
universal lower-bounded generalized $V$-modules generated by the $A^{N}(V)$-modules 
$\Omega_{N}^{0}(W_{2})$ and $\Omega_{N}^{0}(W_{3}')$, respectively,
is formulated and proved in Section 6.

\paragraph{Acknowledgment} I am grateful to Robert McRae for helpful comments. 

\renewcommand{\theequation}{\thesection.\arabic{equation}}
\setcounter{equation}{0}
\setcounter{thm}{0}
\section{The associative algebras $A^{\infty}(V)$ and $A^{N}(V)$ 
and graded $A^{\infty}(V)$- and $A^{N}(V)$-modules}

In this section, we recall the associative algebras 
$A^{\infty}(V)$ and $A^{N}(V)$ for $N\in \N$ for a grading-restricted 
vertex algebra $V$ and  $A^{\infty}(V)$- and $A^{N}(V)$-graded modules
introduced in \cite{H-aa-va}. 

Let $V$ be a grading-restricted vertex algebra. 
Let $U^{\infty}(V)$ be
the space of  column-finite infinite
matrices with entries in $V$, but doubly indexed by $\N$ instead of $\Z_{+}$. 
Elements of  $U^{\infty}(V)$ are of the form 
$\mathfrak{v}=[v_{kl}]$
for $v_{kl}\in V$, $k, l\in \N$  such that for each fixed  $l\in \N$, 
there are only finitely many nonzero $v_{kl}$. 
For $k, l\in \N$ and $\in V$, let $[v]_{kl}$ be the element of $U^{\infty}(V)$
with the entry in the $k$-th row and $l$-th column equal to $v$ and 
all the other entries equal to $0$. 

We define a product $\diamond$ on  $U^{\infty}(V)$ by
$$\mathfrak{u}\diamond \mathfrak{v}=[(\mathfrak{u}\diamond\mathfrak{v})_{kl}]$$
for $\mathfrak{u}=[u_{kl}], \mathfrak{v}=[v_{kl}]\in U^{\infty}(V)$, where 
\begin{align}\label{defn-diamond}
(\mathfrak{u}\diamond\mathfrak{v})_{kl}&=\sum_{n=k}^{l}\res_{x}T_{k+l+1}((x+1)^{-k+n-l-1})
(1+x)^{l}Y_{V}((1+x)^{L_{V}(0)}u_{kn}, x)v_{nl}\nn
&=\sum_{n=k}^{n}\sum_{m=0}^{l}
\binom{-k+n-l-1}{m}\res_{x}
x^{-k+n-l-m-1}(1+x)^{l}Y_{V}((1+x)^{L_{V}(0)}u_{kn}, x)v_{nl}
\end{align}
for $k, l\in \N$, where 
$$T_{k+l+1}((x+1)^{-k+n-l-1})
=\sum_{m=0}^{n}\binom{-k+n-l-1}{m}x^{-k+n-l-m-1}$$
is the Taylor polynomial in $x^{-1}$ of order $k+l+1$ of $(x+1)^{-k+n-l-1}$. 
Then $U^{\infty}(V)$ equipped with $\diamond$ is an algebra 
but in general is not 
even associative. 
By definition, 
for $u, v\in V$ and $k, m, n, l\in \N$, by definition, 
$$[u]_{km}\diamond[v]_{nl}=0$$
when $m\ne n$ and 
\begin{align}\label{defn-diamond-1}
[u]_{kn}\diamond[v]_{nl}&=
\res_{x}T_{k+l+1}((x+1)^{-k+n-l-1})(1+x)^{l}\left[Y_{V}((1+x)^{L_{V}(0)}u, x)v\right]_{kl}\nn
&=\sum_{m=0}^{n}
\binom{-k+n-l-1}{m}\res_{x}
x^{-k+n-l-m-1}(1+x)^{l}\left[Y_{V}((1+x)^{L_{V}(0)}u, x)v\right]_{kl}.
\end{align}
Since $[u]_{km}\diamond[v]_{nl}=0$ when $m\ne n$, we need only consider 
$[u]_{kn}\diamond[v]_{nl}$ for $u, v\in V$ and $k, n, l\in \N$.

Let $\one^{\infty}$ be the element of $U^{\infty}(V)$ with diagonal entries being $\one \in V$ and
all the other entries being $0$. 

Let $W$ be a lower-bounded generalized $V$-module. 
For $n\in \N$, let 
$$\Omega_{n}(W)=\{w \in W \mid  (Y_{W})_{k}(v)w = 0\ {\rm for\ homogeneous}\ v\in V, 
\wt v - k - 1 < -n\}.$$
Then
$$\Omega_{n_{1}}(W)\subset \Omega_{n_{2}}(W)$$
for $n_{1}\le n_{2}$
and 
$$W=\bigcup_{n\in \N}\Omega_{n}(W).$$
So $\{\Omega_{n}(w)\}_{n\in \N}$ is an ascending filtration of $W$. 
Let 
$$Gr(W)=\sum_{n\in \N}Gr_{n}(W)$$
 be the 
associated graded space, where 
$$Gr_{n}(W)=\Omega_{n}(W)/\Omega_{n-1}(W).$$
Sometimes we shall use $[w]_{n}$ to denote the element $w+\Omega_{n-1}(W)$ of $Gr_{n}(W)$,
where $w\in \Omega_{n}(W)$. 

For $\mathfrak{v}=[v_{kl}]\in U^{\infty}(V)$, where $v_{kl}\in V$ and $k, l\in \N$, we have an operator
$\vartheta_{Gr(W)}(\mathfrak{v})$ on $Gr(W)$ defined
by
$$\vartheta_{Gr(W)}(\mathfrak{v})\mathfrak{w}=\sum_{k, l\in \N}
\res_{x}x^{l-k-1}Y_{W}(x^{L_{V}(0)}v_{kl}, x)\pi_{Gr_{l}(W)}\mathfrak{w},
$$
for $\mathfrak{w}\in Gr(W)$,
where $\pi_{Gr_{l}(W)}$ is the projection from 
$Gr(W)$ to $Gr_{l}(W)$. 
Then we have a linear map 
\begin{align*}
\vartheta_{Gr(W)}: U^{\infty}(V)&\to \mbox{\rm End}\; Gr(W)\nn
\mathfrak{v}&\mapsto \vartheta_{Gr(W)}(\mathfrak{v}). 
\end{align*}
Let $Q^{\infty}(V)$ be the intersection of $\ker \vartheta_{Gr(W)}$ for all lower-bounded generalized $V$-modules $W$ and
$A^{\infty}(V)=U^{\infty}(V)/Q^{\infty}(V)$. 

The following result gives the associative algebra $A^{\infty}(V)$: 

\begin{thm}[\cite{H-aa-va}]\label{assoc-alg}
The product $\diamond$ on $U^{\infty}(V)$ induces a product, denoted still by $\diamond$, on 
$A^{\infty}(V)=U^{\infty}(V)/Q^{\infty}(V)$ such that $A^{\infty}(V)$  equipped with $\diamond$ is an associative algebra
with $\one^{\infty}+Q^{\infty}(V)$ as identity. 
Moreover, the associated graded space $Gr(W)$ of the ascendant
filtration $\{\Omega_{n}(W)\}_{n\in \N}$ of a lower-bounded generalized $V$-module $W$  is an 
$A^{\infty}(V)$-module. 
\end{thm}

We also need the following notion introduced in \cite{H-aa-va}:

\begin{defn}\label{N-graded-A-inf-mod}
{\rm Let $G$ be an $A^{\infty}(V)$-module with the $A^{\infty}(V)$-module structure on $G$ given by a homomorphism
$\vartheta_{G}: A^{\infty}(V)\to {\rm End}\; G$ of 
associative algebras. We say that $G$ is a {\it graded $A^{\infty}(V)$-module}
if the following conditions are 
satisfied:

\begin{enumerate}

\item $G$ is graded by $\N$, that is, $G=\coprod_{n\in \N}G_{n}$, and for $v\in V$, $k, l\in \N$,
$\vartheta_{G}([v]_{kl}+Q^{\infty}(V))$ maps $G_{n}$ to 
$0$ when $n\ne l$ and to $G_{k}$ when $n=l$. 

\item $G$ is a direct sum of generalized eigenspaces of  an operator $L_{G}(0)$ on $G$, $G_{n}$ for $n\in \N$ 
are invariant under $L_{G}(0)$ and
the real parts of the eigenvalues of $L_{G}(0)$ have a lower bound. 

\item There is an operator $L_{G}(-1)$ on $G$ mapping $G_{n}$ to $G_{n+1}$ for $n\in \N$.

\item The commutator relations
\begin{align*}
{[L_{G}(0), L_{G}(-1)]}&=L_{G}(-1),\nn
{[L_{G}(0), \vartheta_{G}([v]_{kl}+Q^{\infty}(V))]}&=(k-l)\vartheta_{G}([v]_{kl}+Q^{\infty}(V)),\nn
{[L_{G}(-1), \vartheta_{G}([v]_{kl}+Q^{\infty}(V))]}&=\vartheta_{G}([L_{V}(-1)v]_{(k+1)l}+Q^{\infty}(V))
\end{align*}
hold for $v\in V$ and 
$k, l\in \N$

\end{enumerate}
A graded $A^{\infty}(V)$-algebra $G$ is said to be {\it nondegenerate} if it satisfies in addition the following condition:
For $\textsl{g}\in G_{l}$, if $\vartheta_{G}([v]_{0l}+Q^{\infty}(V))\textsl{g}=0$ for all $v\in V$, then $\textsl{g}=0$. 
Let $G_{1}$ and $G_{2}$ be graded $A^{\infty}(V)$-modules.
A {\it graded $A_{\infty}(V)$-module map} from 
$G_{1}$ to $G_{2}$ is an $A^{N}(V)$-module map $f: G_{1}\to G_{2}$ such that 
$f((G_{1})_{n})\subset (G_{2})_{n}$,  
$f\circ L_{G_{1}}(0)=L_{G_{2}}(0)\circ f$ and $f\circ L_{G_{1}}(-1)=L_{G_{2}}(-1)\circ f$.
A {\it graded $A^{\infty}(V)$-submodule} of a graded $A^{\infty}(V)$-module $G$ is an $A^{\infty}(V)$-submodule 
 of $G$ that is also an $\N$-graded subspace of $G$ and invariant under 
the operators $L_{G}(0)$ and $L_{G}(-1)$.
A graded $A^{\infty}(V)$-module $G$ is said to be {\it generated by a subset $S$} if 
$G$  is equal  to the smallest graded $A^{\infty}(V)$-submodule containing $S$, or equivalently, 
$G$  is spanned by homogeneous elements with respect to the $\N$-grading and the grading given by 
$L_{G}(0)$ obtained by applying elements of $A^{\infty}(V)$, $L_{G}(0)$ and $L_{G}(-1)$ to 
homogeneous summands of elements of $S$. 
 A graded $A^{\infty}(V)$-module is said to be 
{\it irreducible} if it has no nonzero proper graded $A^{\infty}(V)$-submodules. A graded $A^{\infty}(V)$-module 
is said to be {\it completely reducible} if it is a direct sum of irreducible graded $A^{\infty}(V)$-modules.}
\end{defn}

We now recall the subalgebras  $A^{N}(V)$ 
of $A^{\infty}(V)$ also introduced in
\cite{H-aa-va}. For $N\in \N$, 
let $U^{N}(V)$ be the space of all $(N+1)\times (N+1)$ 
matrices with entries in $V$. It is clear that $U^{N}(V)$ can be canonically
embedded into $U^{\infty}_{0}(V)$ as a subspace. We view $U^{N}(V)$ 
as a subspace of $U^{\infty}_{0}(V)$. 
As a subspace of   $U^{\infty}_{0}(V)$, $U^{N}(V)$ 
consists of infinite matrices in $U^{\infty}(V)$
whose $(k, l)$-th entries for $k> N$ or $l> N$ are all $0$ and is
 spanned by elements of the form 
$[v]_{kl}$ for $v\in V$, $k, l=0, \dots, N$.

By  (\ref{defn-diamond-1}), for $u, v\in V$ and $k, n, l=0, \dots, N$,
\begin{equation}\label{defn-diamond-1-N}
[u]_{kn}\diamond [v]_{nl}
=\res_{x}T_{k+l+1}((x+1)^{-k+n-l-1})(1+x)^{l}\left[Y_{V}((1+x)^{L(0)}u, x)v\right]_{kl}\in U^{N}(V).
\end{equation}
So $U^{N}(V)$ is closed under the product $\diamond$. 
Let
$$A^{N}(V)=\{\mathfrak{v}+Q^{\infty}(V)\;|\;\mathfrak{v}\in U^{N}(V)\}
=\pi_{A^{\infty}(V)}(U^{N}(V)),$$
where 
$\pi_{A^{\infty}(V)}$ is the projection from $U^{\infty}(V)$ to 
$A^{\infty}(V)$. Then $A^{N}(V)$ is
spanned by elements of the form 
$[v]_{kl}+Q^{\infty}(V)$ for $v\in V$ and $k, l=0, \dots, N$.
Let 
$\one^{N}=\sum_{k=0}^{N}[\one]_{kk},$
that is, $\one^{N}$ is the element of $U^{N}(V)$ 
with the only nonzero entries to be 
equal to $\one$ at the diagonal $(k, k)$-th entries for $k=0, \dots, N$.

\begin{prop}
The subspace $A^{N}(V)$ is closed under $\diamond$ and is thus a  subalgebra
of $A^{\infty}(V)$ with the identity $\one^{N}+Q^{\infty}(V)$. 
\end{prop}

We also have the following notion introduced in \cite{H-aa-va}:

\begin{defn}\label{gr-A-N-mod}
{\rm Let $M$ be an $A^{N}(V)$-module $M$ with the $A^{N}(V)$-module structure on $M$ given by 
$\vartheta_{M}: A^{N}(V)\to {\rm End}\; M$. We say that $M$ is a {\it graded $A^{N}(V)$-module}
if the following conditions are 
satisfied:
\begin{enumerate}

\item  $M=\coprod_{n=0}^{N}G_{n}(M)$ such that for $v\in V$ and $k, l=0, \dots, N$,
$\vartheta_{M}([v]_{kl}+Q^{\infty}(V))$ maps $G_{n}(M)$ for $0\le n\le N$ to 
$0$ when $n\ne l$ and to $G_{k}(M)$ when $n=l$.

\item $M$ is  a direct sum of generalized eigenspaces of of an operator $L_{M}(0)$
on $M$. $G_{n}(M)$ for $n\in \N$ are invariant under $L_{M}(0)$ and 
the real parts of the eigenvalues of $L_{M}(0)$ has a lower bound. 

\item There is a linear map $L_{M}(-1): \coprod_{n=0}^{N-1}G_{n-1}(M)\to 
\coprod_{n=1}^{N}G_{n}(M)$
mapping $G_{n}(M)$ to $G_{n+1}(M)$ for $n=0, \dots, N-1$.

\item The commutator relations
\begin{align*}
{[L_{M}(0), L_{M}(-1)]}&=L_{M}(-1),\nn
{[L_{M}(0), \vartheta_{M}([v]_{kl}+Q^{\infty}(V))]}&=(k-l)\vartheta_{M}([v]_{kl}+Q^{\infty}(V)),\nn
{[L_{M}(-1), \vartheta_{M}([v]_{pl}+Q^{\infty}(V))]}&=\vartheta_{M}([L_{V}(-1)v]_{(p+1)l}+Q^{\infty}(V))
\end{align*}
hold for $v\in V$, 
$k, l=0, \dots, N$ and $p=0, \dots, N-1$. 

\end{enumerate}
A graded $A^{N}(V)$-module $M$ is said to be {\it nondegenerate} if 
the following additional condition holds: For $w\in G_{l}(M)$, if $\vartheta_{M}([v]_{0l}+Q^{\infty}(V))w=0$ 
for all $v\in V$, then $w=0$. 
Let $M_{1}$ and $M_{2}$ be graded $A^{N}(V)$-modules.
An {\it graded $A_{N}(V)$-module map} from 
$M_{2}$ to $M_{2}$ is an $A^{N}(V)$-module map $f: M_{1}\to M_{2}$ 
such that $f(G_{n}(M_{1}))\subset G_{n}(M_{2})$ for $n=0 \dots, N$, 
$f\circ L_{M_{1}}(0)=L_{M_{2}}(0)\circ f$ and $f\circ L_{M_{1}}(-1)=L_{M_{2}}(-1)\circ f$.
A {\it graded $A^{N}(V)$-submodule} of a graded $A^{N}(V)$-module $M$ is an $A^{N}(V)$-submodule 
$M_{0}$ of $M$ such that  with the $A^{N}(V)$-module structure, the $\N$-grading induced from $M$
and the operators $L_{M}(0)\mbar_{M_{0}}$ and $L_{M}(-1)\mbar_{M_{0}}$, 
$M_{0}$ is a graded $A^{N}(V)$-module.
A graded $A^{\infty}(V)$-module $M$ is said to be {\it generated by a subset $S$}  if 
$M$  is equal  to the smallest graded $A^{N}(V)$-submodule containing $S$, or equivalently, 
$M$  is spanned by homogeneous elements 
obtained by applying elements of $A^{N}(V)$, $L_{M}(0)$ and $L_{M}(-1)$ to 
homogeneous summands of  elements of $S$.  A graded $A^{N}(V)$-module is said to be 
{\it irreducible} if it has no nonzero proper graded $A^{N}(V)$-modules. A graded $A^{N}(V)$-module 
is said to be {\it completely reducible} if it is a direct sum of irreducible graded $A^{N}(V)$-modules.}
\end{defn}

\renewcommand{\theequation}{\thesection.\arabic{equation}}
\setcounter{equation}{0}
\setcounter{thm}{0}
\section{Lower-bounded generalized $V$-modules and $A^{\infty}(V)$-modules}

In this section, $V$ is a grading-restricted vertex algebra. In particular, we do not assume that 
$V$ has a conformal vector. 
We give an $A^{\infty}(V)$-module structure to a lower-bounded generalized $V$-module
$W$ in this section.  For the associative algebra $A^{\infty}(V)$,
its subalgebras, its modules and related structures and 
results, see Section 2 and \cite{H-aa-va}. 

Let $W$ be a lower-bounded generalized $V$-module. 
Let $W^{\mu}$ for $\mu\in \C/\Z$ be the generalized $V$-submodule of $W$ spanned by 
homogeneous elements of weights in $\mu$. Let 
$$\Gamma(W)=\{\mu\in \C/\Z\mid W^{\mu}\ne 0\}.$$
We call $\Gamma(W)$ the {\it sets of congruence classes of weights of $W$}.

For $\mu \in \Gamma(W)$, there exists $h^{\mu}\in \C$ such that 
$$W^{\mu}=\coprod_{n\in \N}W_{[h^{\mu}+n]}$$
and $W_{[h^{\mu}]}\ne 0$. 
Then
$$W=\coprod_{\mu\in \Gamma(W)}W^{\mu}=\coprod_{\mu\in \Gamma(W)}\coprod_{n\in \N}W_{[h^{\mu}+n]}.$$
For $n\in \N$, let 
$$W_{\l n\r}=\coprod_{\mu\in \Gamma(W)}W_{[h^{\mu}+n]}.$$
Then $W$ has a canonical $\N$-grading
$$W
=\coprod_{n\in \N}W_{\l n\r}.$$
For $n\in \N$, we call the space $W_{\l n\r}$ the {\it $n$-the level of $W$} and for an element $w\in W_{\l n\r}$, 
we call $n$ the {\it level of $w$}. From the definition, we have $W_{\l n\r}
\subset \Omega_{n}(W)$.

We define a linear map 
$\vartheta_{W}: U^{\infty}(V) \to \text{End}\;W$ by
$$\vartheta_{W}([v]_{kl})w=\delta_{ln}\res_{x}x^{l-k-1}Y_{W}(x^{L_{V}(0)}v, x)w$$
for $k, l\in \N$, $v\in V$ and $w\in W_{\l n\r}$. 

We now want to show that $W$ is in fact a graded
$A^{\infty}(V)$-module. For simplicity, we discuss only the case that 
$W=W^{\mu}=\coprod_{n\in \N}W_{[h^{\mu}+n]}$ for some $\mu\in \Gamma(W)$.
The general case follows immediately  from the decomposition 
$W=\coprod_{\mu\in \Gamma(W)}W^{\mu}$. 

We need to use the construction in Section 5 of \cite{H-const-twisted-mod}. 
Take the generating fields for the grading-restricted vertex algebra $V$ to be 
$Y_{V}(v, x)$ for $v\in V$. By definition, 
$W$ is a direct sum of generalized eigenspaces of $L_{W}(0)$ and the real parts of the eigenvalues of 
$L_{W}(0)$ has a lower bound $\Re(h^{\mu})\in \R$. 
We take $M$ and $B$ in Section 5 of \cite{H-const-twisted-mod}
to be $W$ and $\Re(h^{\mu})$, respectively. 
Using the construction in  Section 5 of \cite{H-const-twisted-mod},
we obtain a universal lower-bounded generalized $V$-module 
$\widehat{W}^{[1_{V}]}_{\Re(h^{\mu})}$. For simplicity, we shall denote it simply by 
$\widehat{W}$. 

By Theorem 3.3 in \cite{H-exist-twisted-mod} and the construction in  Section 5 of \cite{H-const-twisted-mod}
and by identifying elements of the form 
$(\psi_{\widehat{W}}^{a})_{-1, 0}\one$ with basis elements $w^{a}\in W$ for $a\in A$
for a basis $\{w^{a}\}_{a\in A}$ of $W$, 
we see that $\widehat{W}$ is generated by $W$ (in the sense of Definition 3.1 in \cite{H-exist-twisted-mod}). 
Moreover,  after identifying $(\psi_{\widehat{W}}^{a})_{-1, 0}\one$ with basis elements 
$w^{a}\in G$ for $a\in A$,
Theorems 3.3 and 3.4 in \cite{H-exist-twisted-mod}
say that elements of the form 
$L_{\widehat{W}}(-1)^{p}w^{a}$ for $p\in \N$ and $a\in A$ are linearly independent and 
$\widehat{W}$ is spanned by elements obtained by applying the components of the vertex operators
to these elements.  In particular, $W$ can be embedded into $\widehat{W}$ as 
a subspace. So we shall view $W$ as a subspace of $\widehat{W}$. But note that $W$ is not a 
$V$-submodule of $\widehat{W}$ since $W$ is not invariant under the action of vertex operators on 
$\widehat{W}$. 

Since $W=\coprod_{n\in \N}W_{[h^{\mu}+n]}$, from the construction of $\widehat{W}$
in Section 5 of \cite{H-const-twisted-mod}, we also have 
$\widehat{W}=\coprod_{n\in \N}\widehat{W}_{[h^{\mu}+n]}$. For $n\in \N$, 
$W_{[h^{\mu}+n]}$ is a subspace of $\widehat{W}_{[h^{\mu}+n]}$.

\begin{lemma}\label{widehat-W}
For $n\in \N$, $W_{[h^{\mu}+n]}\cap\Omega_{n-1}(\widehat{W})=0$. In particular, 
we can view $W_{[h^{\mu}+n]}$ as a subspace 
of $Gr_{n}(\widehat{W})$. 
\end{lemma}
\pf
Let $w\in W_{[h^{\mu}+n]}\cap \Omega_{n-1}(\widehat{W})$. 
Then for homogeneous $v\in V$,
$$\res_{x}x^{n-1}Y_{\widehat{W}}(x^{L_{V}(0)}v, x)w=(Y_{\widehat{W}})_{\swt +n-1}(v)w=0.$$
When $w$ is a basis element $w^{a}$ of $W$ as in \cite{H-const-twisted-mod},
the element $(Y_{\widehat{W}})_{\swt +n-1}(v)w^{a}$ can be written as 
$(Y_{\widehat{W}})_{\swt v+n-1}(v)(\psi_{\widehat{W}}^{a})_{-1, 0}\one$. 
Since $w^{a}\ne 0$ as a basis element, 
\begin{equation}\label{widehat-W-1}
(Y_{\widehat{W}})_{\swt +n-1}(v)(\psi_{\widehat{W}}^{a})_{-1, 0}\one=0
\end{equation}
for $v\in V$ and $a\in A$ is a relation in $\widehat{W}$ with the left-hand side being of weight $h^{\mu}$. 
But relations in $\widehat{W}$ of weight $h^{\mu}$ (or of level $0$) and involving 
the operator $(\psi_{\widehat{W}}^{a})_{-1, 0}$  must
have elements of $V$ not equal to $\one$ to the right of 
 $(\psi_{\widehat{W}}^{a})_{-1, 0}$ (see Section 5 of \cite{H-const-twisted-mod}).
 So (\ref{widehat-W-1}) is not a relation 
in $\widehat{W}$. Thus we must have $w^{a}=0$ for $a\in A$, that is,
$W_{[h^{\mu}+n]}\cap\Omega_{n-1}(\widehat{W})=0$.

By the construction of $\widehat{W}$, $W_{[h^{\mu}+n]}\subset \Omega_{n}(\widehat{W})$. 
We define a linear map from $W_{[h^{\mu}+n]}$ to $Gr_{n}(\widehat{W})$ by 
$w\mapsto w+\Omega_{n-1}(\widehat{W})$. 
Since $W_{[h^{\mu}+n]}\cap\Omega_{n-1}(\widehat{W})=0$, this map is injective. 
So we can view $W_{[h^{\mu}+n]}$ as a subspace 
of $Gr_{n}(\widehat{W})$. 
\epfv

Let $J_{W}$  be  the generalized $V$-submodule 
of $\widehat{W}$ generated by elements of the forms
\begin{equation}\label{J-W}
\res_{x}x^{l-k-1}Y_{\widehat{W}}(x^{L_{V}(0)}v, x)w-\res_{x}x^{l-k-1}Y_{W}(x^{L_{V}(0)}v, x)w
\end{equation}
for $v\in V$, $k, l\in \N$, $w\in W_{[h^{\mu}+l]}$ and 
\begin{equation}\label{J-W-2}
L_{\widehat{W}}(-1)w-L_{W}(-1)w
\end{equation}
for $w\in W$. Let 
$$G_{n}(J_{W})=\{w+\Omega_{n-1}(\widehat{W})\mid w\in \Omega_{n}(J_{W})\}
\subset Gr_{n}(\widehat{W})$$
for $n\in \N$ and let
$$G(J_{W})=\coprod_{n\in \N}G_{n}(J_{W}).$$
Since $J_{W}$ is a generalized $V$-submodule of $\widehat{W}$ 
and $Gr_{n}(\widehat{W})$ is an $A^{\infty}(V)$-module, $G(J_{W})$ is an 
$A^{\infty}(V)$-submodule of $Gr(\widehat{W})$.

\begin{prop}\label{W-A-infty-mod}
The $\N$-graded space $W$ with the action of $U^{\infty}(V)$ given by $\vartheta_{W}$ induces a
graded $A^{\infty}(V)$-module structure on $W$ canonically equivalent to the quotient $A^{\infty}(V)$-module 
$Gr(\widehat{W})/G(J_{W})$. 
\end{prop}
\pf
As we have done above, we prove only the case that $W=W^{\mu}$ for some $\mu\in \Gamma(W)$. 
The general case follows immediately.

Since $\widehat{W}$ is generated by $W$, $\widehat{W}$ is 
spanned by elements of the form 
\begin{align*}
&\res_{x}x^{l-k-1}Y_{\widehat{W}}(x^{L_{V}(0)}v, x)w\nn
&\quad =(\res_{x}x^{l-k-1}Y_{\widehat{W}}(x^{L_{V}(0)}v, x)w-\res_{x}x^{l-k-1}Y_{W}(x^{L_{V}(0)}v, x)w)\nn
&\quad \quad +
\res_{x}x^{l-k-1}Y_{W}(x^{L_{V}(0)}v, x)w\nn
&\quad \in J_{W}+W
\end{align*}
for $v\in V$, $k, l\in \N$, $w\in W_{[h^{\mu}+l]}$
and 
$$L_{\widehat{W}}(-1)w=(L_{\widehat{W}}(-1)w-L_{W}(-1)w)+L_{W}(-1)w\in J_{W}+W$$
for $w\in W$. From the definition of $J_{W}$, we have $J_{W}\cap W=0$. 
Hence $\widehat{W}=J_{W}\oplus W$. 
From this decomposition 
and Lemma \ref{widehat-W}, we have $Gr_{n}(\widehat{W})=G_{n}(J_{W})\oplus W_{[h^{\mu}+n]}$
and $Gr(\widehat{W})=G(J_{W})\oplus W$. 
Thus the $\N$-graded space $W$ is canonically isomorphic to $Gr(\widehat{W})/G(J_{W})$. 
Since $Gr(\widehat{W})$ is an $A^{\infty}(V)$-module and $G(J_{W})$ is an $A^{\infty}(V)$-submodule of 
$Gr(\widehat{W})$, we see that $Gr(\widehat{W})/G(J_{W})$ as a quotient of $A^{\infty}(V)$ is also an
$A^{\infty}(V)$-module.  

Let $f$ be the canonical isomorphism from $W$ to $Gr(\widehat{W})/G(J_{W})$. Then
$$f(w)=(w+\Omega_{l-1}(\widehat{W}))+G_{l}(J_{W})$$ for $w\in W_{[h^{\mu}+l]}$.
We use $\vartheta_{Gr(\widehat{W})/G(J_{W})}([v]_{kn})\tilde{w}$ to denote the 
action of $[v]_{kn}$ on $\tilde{w}\in Gr(\widehat{W})/G(J_{W})$.
Then 
\begin{align*}
&\vartheta_{Gr(\widehat{W})/G(J_{W})}([v]_{kn})f(w)\nn
&\quad =\vartheta_{Gr(\widehat{W})/G(J_{W})}([v]_{kn})((w+\Omega_{l-1}(\widehat{W}))+G_{l}(J_{W}))\nn
&\quad =\delta_{nl}
(\res_{x}x^{l-k-1}Y_{\widehat{W}}(x^{L_{V}(0)}v, x)w+\Omega_{l-1}(\widehat{W}))+G_{l}(J_{W})\nn
&\quad =\delta_{nl}(\res_{x}x^{l-k-1}Y_{W}(x^{L_{V}(0)}v, x)w+\Omega_{l-1}(\widehat{W}))+G_{l}(J_{W})\nn
&\quad =(\vartheta_{W}([v]_{kn})w+\Omega_{l-1}(\widehat{W}))+G_{l}(J_{W})\nn
&\quad =f(\vartheta_{W}([v]_{kn})w)
\end{align*}
for $v\in V$, $k, l, n\in \N$ and $w\in W_{[h^{\mu}+l]}$. So the isomorphism 
$f$ from $W$ to $Gr(\widehat{W})/G(J_{W})$
commutes with the actions of $A^{\infty}(V)$ on $W$ and 
$Gr(\widehat{W})/G(J_{W})$. Thus $f$ is an $A^{\infty}(V)$-module map. Since 
$f$ is also a linear isomorphism, it is an equivalence of $A^{\infty}(V)$-modules. 

The four conditions for $W$ to be a graded $A^{\infty}(V)$-module 
in Definition \ref{N-graded-A-inf-mod}
are in fact the properties of the $\N$-grading of $W$, 
the operators $L_{W}(0)$ and $L_{W}(-1)$ and $\vartheta_{W}$. 
\epfv

\begin{rema}\label{graded-mod}
{\rm Note that in general, $W$ is not nondegenerate as a graded $A^{\infty}(V)$-module (see Definition \ref{N-graded-A-inf-mod}).}
\end{rema}

\begin{rema}\label{alt-defn-Q-infty-V}
{\rm It is easy to see that $\ker \vartheta_{W}\subset \ker \vartheta_{Gr(W)}$. In fact, 
let $\mathfrak{v}\in \ker \vartheta_{W}$. Then for $k, l\in \N$, $w\in W_{\l l \r}$,
$\pi_{W_{\l k\r}}\vartheta_{W}(\mathfrak{v})w=0$, where $\pi_{W_{\l k\r}}$ is the projection 
from $W$ to $W_{\l k\r}$. By the definition of $\vartheta_{W}$,
we know that $\pi_{W_{\l k\r}}\vartheta_{W}(\mathfrak{v})w$ must be of the form 
$\vartheta_{W}([v]_{kl})w$ for some $v\in V$. So we have $\vartheta_{W}([v]_{kl})w=0$, or explicitly, 
$\res_{x}x^{l-k-1}Y_{W}(x^{L_{V}(0)}v, x)w=0$. Then we also have 
$$\pi_{Gr_{k}(W)}\vartheta_{Gr(W)}(\mathfrak{v})[w]_{l}=
\vartheta_{Gr(W)}([v]_{kl})[w]_{l}=[\res_{x}x^{l-k-1}Y_{W}(x^{L_{V}(0)}v, x)w]_{k}
=0$$
in $Gr_{k}(W)$. Since $k, l$ and $w$ are arbitrary, we obtain $\vartheta_{Gr(W)}(\mathfrak{v})=0$,
that is, $\mathfrak{v}\in \ker \vartheta_{Gr(W)}$. So $\ker \vartheta_{W}\subset \ker \vartheta_{Gr(W)}$.
This fact together with Proposition \ref{W-A-infty-mod}
 means that $Q^{\infty}(V)$ is also equal to the intersection of 
$\ker \vartheta_{W}$ for all lower-bounded generalized $V$-modules $W$. 
This fact gives new definitions of $Q^{\infty}(V)$ and $A^{\infty}(V)$.
In this paper, we shall use these new definitions.}
\end{rema}

\renewcommand{\theequation}{\thesection.\arabic{equation}}
\setcounter{equation}{0}
\setcounter{thm}{0}
\section{Intertwining operators and $A^{\infty}(V)$-bimodules}

In this section, $V$ is still a grading-restricted vertex algebra which does not have to have a 
conformal vector. We construct an $A^{\infty}(V)$-bimodule $A^{\infty}(W)$ from a lower-bounded 
generalized $V$-module using all intertwining operators of type $\binom{W_{3}}{WW_{2}}$ 
for all lower-bounded generalized $V$-modules $W_{2}$ and $W_{3}$. 

Let $U^{\infty}(W)$ be the space of all column-finite infinite matrices with entries 
in $W$ and doubly index by $\N$. 
For $w\in W$ and $k, l\in \N$, we use $[w]_{kl}$ to denote the  infinite matrix
in $U^{\infty}(W)$
with the $(k, l)$-th entry equal to $w$ and all the other entries equal to $0$. 
Then  elements of $U^{\infty}(W)$ are infinite linear combinations of elements of the form 
$[w]_{kl}$
for $w\in W$ and $k, l\in \N$ with only finitely many elements of such a form for each fixed $k, l\in \N$. 
We can use  elements of the form 
$[w]_{kl}$ for 
$w\in W$ and $k, l\in \N$
to study $U^{\infty}(W)$. 

For $v\in V$, $w\in W$ and $k, m, n, l\in \N$, we define
$$[v]_{km}\diamond [w]_{nl}=0$$
when $m\ne n$ and 
\begin{align}\label{defn-diamond-1}
[v]_{kn}\diamond [w]_{nl}
&=
\res_{x}T_{k+l+1}((x+1)^{-k+n-l-1})(1+x)^{l}
\left[Y_{W}((1+x)^{L_{V}(0)}v, x)w\right]_{kl} \nn
&=\sum_{m=0}^{n}
\binom{-k+n-l-1}{m}\res_{x}
x^{-k+n-l-m-1}(1+x)^{l}
\left[Y_{W}((1+x)^{L_{V}(0)}v, x)w\right]_{kl}.
\end{align}
This is a left action of $U^{\infty}(V)$ on $U^{\infty}(W)$, that is, a linear map
from $U^{\infty}(V)\otimes U^{\infty}(W)$ to $U^{\infty}(W)$. 

For $v\in V$, $w\in W$ and $k, m, n, l\in \N$, we define 
$$[w]_{km} \diamond [v]_{nl}=0$$
when $m\ne n$ and 
\begin{align}\label{defn-diamond-2}
[w]_{kn} \diamond [v]_{nl}
& =\res_{x}T_{k+l+1}((x+1)^{-k+n-l-1})(1+x)^{k}
\left[Y_{W}((1+x)^{-L_{V}(0)}v, -x(1+x)^{-1})w\right]_{kl}\nn
& =\sum_{m=0}^{n}
\binom{-k+n-l-1}{m}\res_{x}
x^{-k+n-l-m-1}(1+x)^{k}\cdot\nn
&\quad\quad\quad\quad  \cdot 
\left[Y_{W}((1+x)^{-L_{V}(0)}v, -x(1+x)^{-1})w\right]_{kl}.
\end{align}
We then obtain a right action
of $U^{\infty}(V)$ on $U^{\infty}(W)$, 
that is, a linear map
from $U^{\infty}(W)\otimes U^{\infty}(V)$ to $U^{\infty}(W)$.  
With the left and right actions of $U^{\infty}(V)$,
$U^{\infty}(W)$ becomes a $U^{\infty}(V)$-bimodule. 
The definition (\ref{defn-diamond-2}) can also be rewritten using 
the $L(0)$-conjugation formula as
\begin{align}\label{defn-diamond-3}
[w]&_{kn} \diamond [v]_{nl}\nn
& =\res_{x}T_{k+l+1}((x+1)^{-k+n-l-1})(1+x)^{k}
\left[(1+x)^{-L_{W}(0)}Y_{W}(v, -x)(1+x)^{L_{W}(0)}w\right]_{kl}\nn
& =\sum_{m=0}^{n}
\binom{-k+n-l-1}{m}\res_{x}
x^{-k+n-l-m-1}(1+x)^{k}\cdot\nn
&\quad\quad\quad\quad  \cdot 
\left[(1+x)^{-L_{W}(0)}Y_{W}(v, -x)(1+x)^{L_{W}(0)}w\right]_{kl}.
\end{align}
The definition (\ref{defn-diamond-2}) can also be rewritten using the right vertex operator
map $Y_{WV}^{V}: W\otimes V\to W[[x, x^{-1}]]$ defined by 
\begin{equation}\label{right-vo}
Y_{WV}^{V}(w, x)v=e^{xL_{W}(-1)}Y_{W}(v, -x)w
\end{equation}
for $v\in V$ and $w\in W$ (see \cite{FHL}). 
To do this, we need a formula involving $L_{W}(-1)$ and $L_{W}(0)$.
It is straightforward to show 
$$\frac{d}{dx}e^{xyL_{W}(-1)}(1+x)^{L_{W}(0)}(1+x)^{-(yL_{W}(-1)+L_{W}(0))}=0$$
Then $e^{xyL_{W}(-1)}(1+x)^{L_{W}(0)}(1+x)^{-(yL_{W}(-1)+L_{W}(0))}$
must be independent of $x$. Setting $x=0$, we obtain 
$$e^{xyL_{W}(-1)}(1+x)^{L_{W}(0)}(1+x)^{-(yL_{W}(-1)+L_{W}(0))}=1_{W},$$
which is equivalent to
$$e^{xyL_{W}(-1)}(1+x)^{L_{W}(0)}=(1+x)^{yL_{W}(-1)+L_{W}(0)}.$$
Let $y=1$. We obtain
\begin{equation}\label{1+x-formula}
e^{xL_{W}(-1)}(1+x)^{L_{W}(0)}=(1+x)^{L_{W}(-1)+L_{W}(0)}.
\end{equation}
Using (\ref{1+x-formula}) and (\ref{right-vo}), we see that (\ref{defn-diamond-3})
can also be rewritten as 
\begin{align}\label{defn-diamond-4}
[w]_{kn} \diamond [v]_{nl}
& =\res_{x}T_{k+l+1}((x+1)^{-k+n-l-1})(1+x)^{k}\cdot\nn
&\quad\quad\quad\quad  \cdot 
\left[(1+x)^{-(L_{W}(-1)+L_{W}(0))}Y_{WV}^{W}((1+x)^{L_{W}(0)}w, x)v\right]_{kl}\nn
& =\sum_{m=0}^{n}
\binom{-k+n-l-1}{m}\res_{x}
x^{-k+n-l-m-1}(1+x)^{k}\cdot\nn
&\quad\quad\quad\quad  \cdot 
\left[(1+x)^{-(L_{W}(-1)+L_{W}(0))}Y_{WV}^{W}((1+x)^{L_{W}(0)}w, x)v\right]_{kl}.
\end{align}
The definition (\ref{defn-diamond-4}) is more conceptual since it says that the right action of 
$U^{\infty}(V)$ on $U^{\infty}(W)$ is defined using the right action (the right vertex operator map)
of $V$ on $W$. 

Now let $W_{2}$ and $W_{3}$ be lower-bounded generalized $V$-modules
and $\Y$ an intertwining operator of type $\binom{W_{3}}{WW_{2}}$. 
As we have mentioned in the 
introduction, for simplicity, by an intertwining operator, we always mean a logarithmic intertwining operator 
defined in Definition 3.10 in \cite{HLZ}, except that in the case that $V$ is a grading-restricted vertex algebra
so that there are
no $L_{V}(1)$ on $V$ and no $L_{W}(1)$ on a lower-bounded generalized $V$-module $W$, 
we do not require that the $L(1)$-commutator formula for intertwining operator hold.
Such an intertwining operator might 
contain the  logarithm of the variable and might even be an infinite power series 
in the logarithm of the variable. 
See \cite{HLZ} for more details. But in addition, we are interested in only 
 intertwining operators $\Y$  such that
the powers of the formal variable $x$ in $\Y(w, x)w_{2}$ for $w\in W$ and $w_{2}\in W_{2}$ 
belong to only finitely many congruence classes in $\C/Z$. 
For such an intertwining operator, its image is in a finite direct sum of $W_{3}^{\nu}$ for $\nu\in \C/\Z$.
We say that such an intertwining operator 
has {\it locally-finite sets of congruence classes of powers} or is an {\it intertwining operator with locally-finite sets of congruence classes of powers}.  
If $\Gamma(W_{3})$ is a finite set, then every intertwining operator of type $\binom{W_{3}}{WW_{2}}$ in the sense of 
 Definition 3.10 in \cite{HLZ} has locally-finite sets of congruence classes of powers.
In this paper,  by an intertwining operator, we always mean an intertwining operator with locally-finite sets of congruence classes of powers.
In particular, the space of intertwining operators 
also means the space of intertwining operators with locally-finite sets of congruence classes of powers.

As is discussed in the preceding section, there exist $h_{2}^{\mu}, h_{3}^{\nu}\in \C$ for $\mu
\in \Gamma(W_{2})$, $\nu\in \Gamma(W_{3})$ such that 
$$W_{2}=\coprod_{n\in \N}(W_{2})_{\l n\r}=\coprod_{n\in \N}
\coprod_{\mu\in \Gamma(W_{2})}(W_{2})_{[h_{2}^{\mu}+n]}, 
\; W_{3}=\coprod_{n\in \N}(W_{3})_{\l n\r}=\coprod_{n\in \N}
\coprod_{\nu\in \Gamma(W_{3})}(W_{3})_{[h_{3}^{\nu}+n]}.$$

For $w\in W$, we have 
$$\Y(w, x)=\sum_{k\in \N}\sum_{m\in \C}\Y_{m, k}(w)x^{-m-1}(\log x)^{k},$$
where for $m\in \C$, $k\in \N$ and homogeneous $w$, 
the map $\Y_{m, k}(w): W_{2}\to W_{3}$ is homogeneous of weight $\wt w-m-1$.
For $w\in W$, let 
$$\Y^{k}(w, x)=\sum_{m\in \C}\Y_{m, k}(w)x^{-m-1}.$$
Note that $\Y^{k}(w, x)$ for $k\in \N$ satisfy the same Jacobi identity as the one for intertwining 
operators but  do not satisfy the $L(-1)$-derivative property for intertwining operators. 

For a vector space $U$ and $k\in \N$, let $\text{Coeff}_{\log x}^{k}: U\{x\}[[\log x]]\to U\{x\}$ be the 
linear map given by taking the coefficients of $(\log x)^{k}$. Then we have 
$$\text{Coeff}_{\log x}^{k}\Y(w, x)=\Y^{k}(w, x)$$
for $w\in W$ and $k\in \N$. 
  
For $k, l\in \N$, $\mu
\in \Gamma(W_{2})$, $\nu\in \Gamma(W_{3})$,   
$w\in W$ and $w_{2}\in (W_{2})_{[h_{2}^{\mu}+l]}\subset \Omega_{l}^{0}(W_{2})$,
\begin{align*}
&\text{Coeff}_{\log x}^{0}\res_{x}x^{h_{2}^{\mu}-h_{3}^{\nu}+l-k-1}
\Y(x^{L_{W}(0)}w, x)w_{2}\nn
&\quad =\res_{x}x^{h_{2}^{\mu}-h_{3}^{\nu}+l-k-1}
\Y^{0}(x^{L_{W}(0)_{S}}w, x)w_{2}\nn
&\quad \in (W_{3})_{[h_{3}^{\nu}+k]} \subset W_{\l k\r}.
\end{align*}
For  $k, l, n\in \N$, 
$\mu \in \Gamma(W_{2})$, $w\in W$ and $w_{2}\in (W_{2})_{[h_{2}^{\mu}+l]}\subset (W_{2})_{\l l\r}$, 
we define 
\begin{align*}
\vartheta_{\Y}([w]_{kl})w_{2}
&=\sum_{\nu\in \Gamma(W_{3})}\text{Coeff}_{\log x}^{0}\res_{x}
x^{h_{2}^{\mu}-h_{3}^{\nu}+l-k-1}
\Y(x^{L_{W}(0)}w, x)w_{2}\nn
&=\sum_{\nu\in \Gamma(W_{3})}\res_{x}x^{h_{2}^{\mu}-h_{3}^{\nu}+l-k-1}
\Y^{0}(x^{L_{W}(0)_{S}}w, x)w_{2}\nn
&\in W_{\l k\r}.
\end{align*}
In this definition, the sum is finite since there are only finitely many congruence classes of powers of $x$ in 
$\Y^{0}(x^{L_{W}(0)_{S}}w, x)w_{2}$. 
We now have a  linear map
$$\vartheta_{\Y}([w]_{kl}): W_{2}\to W_{3},$$
or equivalently, we obtain an element 
$$\vartheta_{\Y}([w]_{kl})\in \hom(W_{2}, W_{3}).$$
The maps $\vartheta_{\Y}(
[w]_{kl})$ for 
$w\in W$ and $k, l\in \N$give a 
linear map
$$\vartheta_{\Y}: U^{\infty}(W)\to \hom(W_{2}, W_{3}).$$
Since $W_{2}$ and $W_{3}$ are both left $A^{\infty}(V)$-modules, 
$\hom(W_{2}, W_{3})$ is an  $A^{\infty}(V)$-bimodule. In particular, we have left and right actions 
of $U^{\infty}(V)$  on $\hom(W_{2}, W_{3})$ such that both the left and right actions of 
$Q^{\infty}(V)$ on $\hom(W_{2}, W_{3})$ are $0$.

\begin{prop}\label{l-r-action-comm}
The linear map $\vartheta_{\Y}$ commutes with the left and right actions 
of $U^{\infty}(V)$.
In particular, $U^{\infty}(W)/\ker \vartheta_{\Y}$ is an $A^{\infty}(V)$-bimodule. 
\end{prop}
\pf
We first prove that $\vartheta_{\Y}$ commutes with the left action of $U^{\infty}(V)$.
Let $k, m, n, l, p\in \N$,  $\mu\in \Gamma(W_{2})$, $v\in V$, $w\in W$, 
$w_{2}\in (W_{2})_{[h_{2}^{\mu}+p]}\subset (W_{2})_{\l p\r}$. 
In the case $m\ne n$ or  $p\ne l$
by definition, we have
$$\vartheta_{\Y}([v]_{km}\diamond [w]_{nl})w_{2}=0=
\vartheta_{W_{3}}([v]_{km})\vartheta_{\Y}( [w]_{nl})w_{2}.$$

In the case $m=n$ and $p=l$,
\begin{align}\label{l-r-action-comm-1}
\vartheta&_{\Y}([v]_{kn}\diamond  [w]_{nl})w_{2}\nn
&=\sum_{\nu\in \Gamma(W_{3})}\res_{x_{0}}T_{k+l+1}((x_{0}+1)^{-k+n-l-1})(1+x_{0})^{l}
\vartheta_{\Y}([Y_{W}((1+x_{0})^{L_{V}(0)}v, x_{0})w]_{kl})
w_{2}\nn
&=\sum_{\nu\in \Gamma(W_{3})}\res_{x_{0}}\res_{x_{2}}
T_{k+l+1}((x_{0}+1)^{-k+n-l-1})(1+x_{0})^{l}x_{2}^{h_{2}^{\mu}-h_{3}^{\nu}+l-k-1}\cdot\nn
&\quad\quad\quad\cdot
\Y^{0}(x_{2}^{L_{W}(0)_{S}}
Y_{W}((1+x_{0})^{L_{V}(0)}v, x_{0})w, x_{2})w_{2}\nn
&=\sum_{\nu\in \Gamma(W_{3})}\res_{x_{0}}\res_{x_{2}}
T_{k+l+1}((x_{0}+1)^{-k+n-l-1})(1+x_{0})^{l}x_{2}^{h_{2}^{\mu}-h_{3}^{\nu}+l-k-1}\cdot\nn
&\quad\quad\quad\cdot
\Y^{0}(Y_{W}((x_{2}+x_{0}x_{2})^{L_{V}(0)}v, x_{0}x_{2})x_{2}^{L_{W}(0)_{S}}w, x_{2})w_{2}\nn
&=\sum_{\nu\in \Gamma(W_{3})}\res_{x_{0}}\res_{x_{2}}\res_{x_{1}}x_{1}^{-1}\delta\left(
\frac{x_{2}+x_{0}x_{2}}{x_{1}}\right)
T_{k+l+1}((x_{0}+1)^{-k+n-l-1})\cdot\nn
&\quad\quad\quad\cdot x_{1}^{l}x_{2}^{h_{2}^{\mu}-h_{3}^{\nu}-k-1}
\Y^{0}(Y_{W}(x_{1}^{L_{V}(0)}v, x_{0}x_{2})x_{2}^{L_{W}(0)_{S}}w, x_{2})w_{2}\nn
&=\sum_{\nu\in \Gamma(W_{3})}\res_{x_{0}}\res_{x_{2}}\res_{x_{1}}x_{0}^{-1}x_{2}^{-1}
\delta\left(\frac{x_{1}-x_{2}}{x_{0}x_{2}}\right)
T_{k+l+1}((x_{0}+1)^{-k+n-l-1})\cdot\nn
&\quad\quad\quad\cdot x_{1}^{l}x_{2}^{h_{2}^{\mu}-h_{3}^{\nu}-k-1}
Y_{W_{3}}(x_{1}^{L_{V}(0)}v, x_{1})\Y^{0}(x_{2}^{L_{W}(0)_{S}}w, x_{2})w_{2}\nn
&\quad-\sum_{\nu\in \Gamma(W_{3})}\res_{x_{0}}\res_{x_{2}}\res_{x_{1}}x_{0}^{-1}
x_{2}^{-1}\delta\left(\frac{x_{2}-x_{1}}{x_{0}x_{2}}\right)
T_{k+l+1}((x_{0}+1)^{-k+n-l-1})\cdot\nn
&\quad\quad\quad\cdot x_{1}^{l}x_{2}^{h_{2}^{\mu}-h_{3}^{\nu}-k-1}
\Y^{0}(x_{2}^{L_{W}(0)_{S}}w, x_{2})Y_{W_{2}}(x_{1}^{L_{V}(0)}v, x_{1})w_{2}\nn
&=\sum_{\nu\in \Gamma(W_{3})}\res_{x_{2}}\res_{x_{1}}
T_{k+l+1}((x_{0}+1)^{-k+n-l-1})\lbar_{x_{0}=(x_{1}-x_{2})x_{2}^{-1}}\cdot\nn
&\quad\quad\quad\cdot x_{1}^{l}x_{2}^{h_{2}^{\mu}-h_{3}^{\nu}-k-2}
Y_{W_{3}}(x_{1}^{L_{V}(0)}v, x_{1})\Y^{0}(x_{2}^{L_{W}(0)_{S}}w, x_{2})w_{2}\nn
&\quad-\sum_{\nu\in \Gamma(W_{3})}\res_{x_{2}}\res_{x_{1}}
T_{k+l+1}((x_{0}+1)^{-k+n-l-1})\lbar_{x_{0}=(-x_{2}+x_{1})x_{2}^{-1}}\cdot\nn
&\quad\quad\quad\cdot x_{1}^{l}x_{2}^{h_{2}^{\mu}-h_{3}^{\nu}-k-2}
\Y^{0}(x_{2}^{L_{W}(0)_{S}}w, x_{2})Y_{W_{2}}(x_{1}^{L_{V}(0)}v, x_{1})w_{2}.
\end{align}

The second term in
the right hand side of (\ref{l-r-action-comm-1}) is $0$ since $w_{2}\in (W_{2})_{[h_{2}^{\mu}+l]}$.
Expanding $T_{k+l+1}((x_{0}+1)^{-k+n-l-1})$ explicitly, we see that  the first term in the right-hand side of 
(\ref{l-r-action-comm-1})  is equal to 
\begin{align}\label{l-r-action-comm-2}
&\sum_{\nu\in \Gamma(W_{3})}\sum_{m=0}^{n}\binom{-k+n-l-1}{m}
\res_{x_{2}}\res_{x_{1}}(x_{1}-x_{2})^{-k+n-l-m-1}x_{2}^{k-n+l+m+1}
\cdot\nn
&\quad\quad\quad\cdot x_{1}^{l}x_{2}^{h_{2}^{\mu}-h_{3}^{\nu}-k-2}
Y_{W_{3}}(x_{1}^{L_{V}(0)}v, x_{1})\Y^{0}(x_{2}^{L_{W}(0)_{S}}w, x_{2})w_{2}\nn
&\quad =\sum_{\nu\in \Gamma(W_{3})}
\sum_{m=0}^{n}\sum_{j\in \N}\binom{-k+n-l-1}{m}\binom{-k+n-l-m-1}{j}(-1)^{j}
\cdot\nn
&\quad\quad\quad\cdot \res_{x_{2}}\res_{x_{1}}x_{1}^{-k+n-m-1-j} 
x_{2}^{h_{2}^{\mu}-h_{3}^{\nu}-n+l+m-1+j}
Y_{W_{3}}(x_{1}^{L_{V}(0)}v, x_{1})\Y^{0}(x_{2}^{L_{W}(0)_{S}}w, x_{2})w_{2}.
\end{align}

Since $w_{2}\in (W_{2})_{[h_{2}^{\mu}+l]}$, we know that
$\res_{x}x^{h_{2}^{\mu}-h_{3}^{\nu}+q-1}\Y^{0}(x^{L_{W_{1}}(0)_{S}}w, x)w_{2}\in 
(W_{3})_{[h_{3}^{\nu}+l-q]}$ and thus is equal to $0$ if $q>l$.
In the case $j>n-m$, we have $-n+l+m+j>l$ and hence 
$$\res_{x_{2}}x_{2}^{h_{2}^{\mu}-h_{3}^{\nu}-n+l+m-1+j}\Y^{0}(x_{2}^{L_{W}(0)_{S}}w, x_{2})w_{2}=0.$$
In particular,  those terms in the right-hand side of (\ref{l-r-action-comm-2}) with $j>n-m$ is $0$.
So the right-hand side of (\ref{l-r-action-comm-2}) is equal to 
\begin{align}\label{l-r-action-comm-3}
&\sum_{\nu\in \Gamma(W_{3})}
\sum_{m=0}^{n}\sum_{j=0}^{n-m}\binom{-k+n-l-1}{m}\binom{-k+n-l-m-1}{j}(-1)^{j}
\cdot\nn
&\quad\quad\quad\cdot \res_{x_{2}}\res_{x_{1}}x_{1}^{-k+n-m-1-j} x_{2}^{h_{2}^{\mu}-h_{3}^{\nu}-n+l+m-1+j}
Y_{W_{3}}(x_{1}^{L_{V}(0)}v, x_{1})\Y^{0}(x_{2}^{L_{W}(0)_{S}}w, x_{2})w_{2}\nn
&\quad=\sum_{\nu\in \Gamma(W_{3})}
\sum_{m=0}^{n}\sum_{q=m}^{n}\binom{-k+n-l-1}{m}\binom{-k+n-l-m-1}{q-m}(-1)^{q-m}
\cdot\nn
&\quad\quad\quad\quad\quad\cdot \res_{x_{2}}\res_{x_{1}}x_{1}^{-k+n-1-q}
 x_{2}^{h_{2}^{\mu}-h_{3}^{\nu}-n+l-1+q}
Y_{W_{3}}(x_{1}^{L_{V}(0)}v, x_{1})\Y^{0}(x_{2}^{L_{W}(0)}v, x_{2})w\nn
&\quad=\sum_{\nu\in \Gamma(W_{3})}\sum_{q=0}^{n}\left(\sum_{m=0}^{q}\binom{-k+n-l-1}{m}
\binom{-k+n-l-m-1}{q-m}(-1)^{q-m}\right)
\cdot\nn
&\quad\quad\quad\quad\quad\cdot \res_{x_{2}}\res_{x_{1}}x_{1}^{-k+n-1-q}
x_{2}^{h_{2}^{\mu}-h_{3}^{\nu}-n+l-1+q}
Y_{W_{3}}(x_{1}^{L_{V}(0)}v, x_{1})\Y^{0}(x_{2}^{L_{W}(0)_{S}}v, x_{2})w.
\end{align}
Using (2.18) in \cite{H-aa-va}, that is,
\begin{equation}\label{binom-identity}
\sum_{m=0}^{q}\binom{-k+n-l-1}{m}\binom{-k+n-l-m-1}{q-m}(-1)^{q-m}
 =\binom{-k+n-l-1}{q}\delta_{q, 0}
\end{equation}
for $q=0, \dots, n$, 
we see that the right-hand side of (\ref{l-r-action-comm-3}) is equal to 
\begin{align}\label{l-r-action-comm-4}
&\sum_{\nu\in \Gamma(W_{3})}\res_{x_{2}}\res_{x_{1}}x_{1}^{-k+n-1} x_{2}^{h_{2}^{\mu}-h_{3}^{\nu}-n+l-1}
Y_{W_{3}}(x_{1}^{L_{V}(0)}v, x_{1})\Y^{0}(x_{2}^{L_{W}(0)_{S}}w, x_{2})w_{2}\nn
&\quad =\vartheta_{W_{3}}([v]_{kn})\sum_{\nu\in \Gamma(W_{3})}
\res_{x_{2}}x_{2}^{h_{2}^{\mu}-h_{3}^{\nu}+l-n-1}\Y^{0}(x_{2}^{L_{W}(0)_{S}}w, x_{2})w_{2}\nn
&\quad =\vartheta_{W_{3}}([v]_{kn})\vartheta_{\Y}([w]_{nl})w_{2}.
\end{align}

From (\ref{l-r-action-comm-1})--(\ref{l-r-action-comm-4}),
we obtain 
$$\vartheta_{\Y}([v]_{kn}\diamond [w]_{nl})w_{2}
=\vartheta_{W_{3}}([v]_{kn})\vartheta_{\Y}([w]_{nl})w_{2}.$$
We have now proved 
$$\vartheta_{\Y}([v]_{km}\diamond [w]_{nl})w_{2}=
\vartheta_{W_{3}}([v]_{km})\vartheta_{\Y}([w]_{nl})w_{2}$$
for $k, l, m, n\in \N$,  $v\in V$, $w\in W$ and $w_{2}\in W_{2}$.
This shows that $\vartheta_{\Y}$ commutes with the left actions of $U^{\infty}(V)$.

Next we prove that $\vartheta_{\Y}$ commutes with the right actions of $U^{\infty}(V)$.
Let $k, m, n, l, p\in \N$, $\mu\in \Gamma(W_{2})$,  $v\in V$, $w\in W$, 
$w_{2}\in (W_{2})_{[h_{2}^{\mu}+p]}\subset (W_{2})_{\l p\r}$.
In the case $m\ne n$ or $p\ne l$, by definition, 
we have
$$\vartheta_{\Y}([w]_{km}\diamond [v]_{nl})w_{2}=0=
\vartheta_{\Y}([w]_{km})\vartheta_{W_{2}}([v]_{nl})w_{2}.$$

In the case $m=n$ and $p= l$,  using (\ref{defn-diamond-3}), we have
\begin{align}\label{l-r-action-comm-5}
\vartheta&_{\Y}([w]_{kn}\diamond [v]_{nl})w_{2}\nn
&=\res_{x_{0}}T_{k+l+1}((x_{0}+1)^{-k+n-l-1})(1+x_{0})^{k}\cdot\nn
&\quad \quad\cdot\vartheta_{\Y}(
[(1+x_{0})^{-L_{W}(0)}Y_{W}(v, -x_{0})(1+x_{0})^{L_{W}(0)}w]_{kl})
w_{2}\nn
&=\sum_{\nu\in \Gamma(W_{3})}\res_{x_{0}}\text{Coeff}_{\log x_{2}}^{0}\res_{x_{2}}
T_{k+l+1}((x_{0}+1)^{-k+n-l-1})(1+x_{0})^{k}
x_{2}^{h_{2}^{\mu}-h_{3}^{\nu}+l-k-1}\cdot\nn
&\quad\quad\cdot
\Y(x_{2}^{L_{W}(0)}(1+x_{0})^{-L_{W}(0)}Y_{W}(v, -x_{0})(1+x_{0})^{L_{W}(0)}w, x_{2})w_{2}\nn
&=\sum_{\nu\in \Gamma(W_{3})}\res_{x_{0}}\text{Coeff}_{\log x_{2}}^{0}\res_{x_{2}}
T_{k+l+1}((x_{0}+1)^{-k+n-l-1})(1+x_{0})^{k}
x_{2}^{h_{2}^{\mu}-h_{3}^{\nu}+l-k-1}\cdot\nn
&\quad\quad\cdot (1+x_{0})^{-L_{W_{3}}(0)}
\Y(x_{2}^{L_{W}(0)}
Y_{W}(v, -x_{0})(1+x_{0})^{L_{W}(0)}w, x_{2}(1+x_{0}))(1+x_{0})^{L_{W_{2}}(0)}w_{2}\nn
&=\sum_{\nu\in \Gamma(W_{3})}\res_{x_{0}}\text{Coeff}_{\log x_{2}}^{0}\res_{x_{2}}
T_{k+l+1}((x_{0}+1)^{-k+n-l-1})\cdot\nn
&\quad\quad\cdot (1+x_{0})^{h_{2}^{\mu}-h_{3}^{\nu}+l}
x_{2}^{h_{2}^{\mu}-h_{3}^{\nu}+l-k-1}
(1+x_{0})^{-L_{W_{3}}(0)}\cdot\nn
&\quad\quad\cdot 
\Y(Y_{W}(x_{2}^{L_{W}(0)}v, -x_{0}x_{2})(x_{2}+x_{0}x_{2})^{L_{W}(0)}w, x_{2}+x_{0}x_{2})
(1+x_{0})^{L_{W_{2}}(0)}w_{2}\nn
&=\sum_{\nu\in \Gamma(W_{3})}\res_{x_{0}}\text{Coeff}_{\log x_{2}}^{0}\res_{x_{2}}\res_{x_{1}}x_{1}^{-1}\delta\left(
\frac{x_{2}+x_{0}x_{2}}{x_{1}}\right)
T_{k+l+1}((x_{0}+1)^{-k+n-l-1})\cdot\nn
&\quad\quad\cdot (x_{2}+x_{0}x_{2})^{h_{2}^{\mu}-h_{3}^{\nu}+l}x_{2}^{-k-1}
(1+x_{0})^{-L_{W_{3}}(0)_{N}}\cdot\nn
&\quad\quad\cdot\Y(Y_{W}(x_{2}^{L_{V}(0)}v, -x_{0}x_{2})
(x_{2}+x_{0}x_{2})^{L_{W}(0)}w, 
x_{2}+x_{0}x_{2})(1+x_{0})^{L_{W_{2}}(0)_{N}}w_{2}.
\end{align}

From the $L(0)$-conjugation property for intertwining operators, 
$$(1+x_{0})^{-L_{W_{3}}(0)}\Y((1+x_{0})^{L_{W}(0)}w, x_{2}+x_{0}x_{2})(1+x_{0})^{L_{W_{2}}(0)}
=\Y(w, x_{2}),$$
or equivalently,
\begin{align}\label{l-r-action-comm-5.3}
&\sum_{k=0}^{K}
(1+x_{0})^{-L_{W_{3}}(0)}\Y^{k}((1+x_{0})^{L_{W}(0)}w, x_{2}+x_{0}x_{2})(1+x_{0})^{L_{W_{2}}(0)}
(\log (x_{2}+x_{0}x_{2}))^{k}\nn
&\quad=\sum_{k=0}^{K}\Y^{k}(w, x_{2})(\log x_{2})^{k}.
\end{align}
Taking $\text{Coeff}_{\log x_{2}}^{0}$ on both sides of (\ref{l-r-action-comm-5.3}), we obtain 
\begin{align*}
&\sum_{k=0}^{K}
(1+x_{0})^{-L_{W_{3}}(0)}\Y^{k}((1+x_{0})^{L_{W}(0)}w, x_{2}+x_{0}x_{2})(1+x_{0})^{L_{W_{2}}(0)}
(\log (1+x_{0}))^{k}\nn
&\quad=\Y^{0}(w, x_{2}).
\end{align*}
Then we have
\begin{align*}
&\text{Coeff}_{\log x_{2}}^{0}(1+x_{0})^{-L_{W_{3}}(0)_{N}}\Y((1+x_{0})^{L_{W}(0)_{N}}w, x_{2}+x_{0}x_{2})
(1+x_{0})^{L_{W_{2}}(0)_{N}}\nn
&\quad=\sum_{k=0}^{K}
(1+x_{0})^{-L_{W_{3}}(0)_{N}}\Y^{k}((1+x_{0})^{L_{W}(0)_{N}}w, x_{2}+x_{0}x_{2})
(1+x_{0})^{L_{W_{2}}(0)_{N}}
(\log (1+x_{0}))^{k}\nn
&\quad=(1+x_{0})^{L_{W_{3}}(0)_{S}}\Y^{0}((1+x_{0})^{-L_{W}(0)_{S}}w, x_{2})
(1+x_{0})^{-L_{W_{2}}(0)_{S}}\nn
&\quad=\Y^{0}(w, x_{2}+x_{0}x_{2}).
\end{align*}
In particular, we have
\begin{align}\label{l-r-action-comm-5.7}
&\text{Coeff}_{\log x_{2}}^{0}(1+x_{0})^{-L_{W_{3}}(0)_{N}}\cdot\nn
&\quad\quad\quad \cdot
\Y(Y_{W}(x_{2}^{L_{V}(0)}v, -x_{0}x_{2})
(x_{2}+x_{0}x_{2})^{L_{W}(0)}w, 
x_{2}+x_{0}x_{2})(1+x_{0})^{L_{W_{2}}(0)_{N}}\nn
&\quad=\text{Coeff}_{\log x_{2}}^{0}(1+x_{0})^{-L_{W_{3}}(0)_{N}}\cdot\nn
&\quad\quad\quad \cdot\Y((1+x_{0})^{L_{W}(0)_{N}}Y_{W}(x_{2}^{L_{V}(0)}v, -x_{0}x_{2})
(x_{2}+x_{0}x_{2})^{L_{W}(0)_{S}}x_{2}^{L_{W}(0)_{N}}w, 
x_{2}+x_{0}x_{2})\cdot\nn
&\quad\quad\quad \cdot (1+x_{0})^{L_{W_{2}}(0)_{N}}\nn
&\quad =\Y^{0}(Y_{W}(x_{2}^{L_{V}(0)}v, -x_{0}x_{2})
(x_{2}+x_{0}x_{2})^{L_{W}(0)_{S}}w, x_{2}+x_{0}x_{2}).
\end{align}
Using (\ref{l-r-action-comm-5.7}), we see that the right-hand side of (\ref{l-r-action-comm-5})
is equal to 
\begin{align}\label{l-r-action-comm-5.8}
&\sum_{\nu\in \Gamma(W_{3})}\res_{x_{0}}\res_{x_{2}}\res_{x_{1}}x_{1}^{-1}\delta\left(
\frac{x_{2}+x_{0}x_{2}}{x_{1}}\right)
T_{k+l+1}((x_{0}+1)^{-k+n-l-1})\cdot\nn
&\quad\cdot (x_{2}+x_{0}x_{2})^{h_{2}^{\mu}-h_{3}^{\nu}+l}x_{2}^{-k-1}
\Y^{0}(Y_{W}(x_{2}^{L_{V}(0)}v, -x_{0}x_{2})
(x_{2}+x_{0}x_{2})^{L_{W}(0)_{S}}w, 
x_{2}+x_{0}x_{2})w_{2}.
\end{align}
On the other hand, 
$$(x_{2}+x_{0}x_{2})^{h_{2}^{\mu}-h_{3}^{\nu}+l}
\Y^{0}(Y_{W}(x_{2}^{L_{V}(0)}v, -x_{0}x_{2})(x_{2}+x_{0}x_{2})^{L_{W}(0)_{S}}w, 
x_{2}+x_{0}x_{2})w_{2}$$
contains only integral powers of $x_{2}+x_{0}x_{2}$ and does not contain $\log (x_{2}+x_{0}x_{2})$.
Then (\ref{l-r-action-comm-5.8})
is equal to
\begin{align}\label{l-r-action-comm-6}
&\sum_{\nu\in \Gamma(W_{3})}\res_{x_{0}}\res_{x_{2}}\res_{x_{1}}x_{2}^{-1}\delta\left(
\frac{x_{1}-x_{0}x_{2}}{x_{2}}\right)
T_{k+l+1}((x_{0}+1)^{-k+n-l-1})\cdot\nn
&\quad\quad\quad\;\cdot x_{1}^{h_{2}^{\mu}-h_{3}^{\nu}+l}
x_{2}^{-k-1}
\Y^{0}(Y_{W}(x_{2}^{L_{V}(0)}v, -x_{0}x_{2})x_{1}^{L_{W}(0)_{S}}w, 
x_{1})w_{2}.
\end{align}

The Jacobi identity for $\Y^{0}$ gives
\begin{align}\label{l-r-action-comm-6.5}
&\res_{x_{1}}x_{2}^{-1}\delta\left(
\frac{x_{1}-x_{0}x_{2}}{x_{2}}\right)
x_{1}^{h_{2}^{\mu}-h_{3}^{\nu}+l}
\Y^{0}(Y_{W}(x_{2}^{L_{V}(0)}v, -x_{0}x_{2})x_{1}^{L_{W}(0)_{S}}w, 
x_{1})w_{2}\nn
&\quad =-\res_{x_{1}}x_{0}^{-1}
\delta\left(\frac{x_{2}-x_{1}}{-x_{0}x_{2}}\right)
 x_{1}^{h_{2}^{\mu}-h_{3}^{\nu}+l}
Y_{W_{3}}(x_{2}^{L_{V}(0)}v, x_{2})\Y^{0}(x_{1}^{L_{W}(0)_{S}}w, x_{1})w_{2}\nn
&\quad\quad +\res_{x_{1}}x_{0}^{-1}\delta\left(\frac{x_{1}-x_{2}}{x_{0}x_{2}}\right)
 x_{1}^{h_{2}^{\mu}-h_{3}^{\nu}+l}
\Y^{0}(x_{1}^{L_{W}(0)_{S}}w, x_{1})Y_{W_{2}}(x_{2}^{L_{V}(0)}v, x_{2})w_{2}.
\end{align}
Since $w_{2}\in (W_{2})_{[h_{2}^{\mu}+l]}$, we have 
\begin{align}\label{l-r-action-comm-6.7}
\res_{x_{1}}x_{1}^{h_{2}^{\mu}-h_{3}^{\nu}+l+j}\Y^{0}(x_{1}^{L_{W}(0)_{S}}w, x_{1})w_{2}
\in (W_{2})_{[h_{3}^{\nu}-1-j]}=0
\end{align}
for $j\in \N$.
Using (\ref{l-r-action-comm-6.5})  and (\ref{l-r-action-comm-6.7}), we see that
(\ref{l-r-action-comm-6}) is equal to 
\begin{align}\label{l-r-action-comm-7}
&\sum_{\nu\in \Gamma(W_{3})}
\res_{x_{0}}\res_{x_{2}}\res_{x_{1}}
x_{0}^{-1}\delta\left(\frac{x_{1}-x_{2}}{x_{0}x_{2}}\right)T_{k+l+1}((x_{0}+1)^{-k+n-l-1})
\cdot\nn
&\quad\quad\quad\cdot
 x_{1}^{h_{2}^{\mu}-h_{3}^{\nu}+l}
x_{2}^{-k-2}
\Y^{0}(x_{1}^{L_{W}(0)_{S}}w, x_{1})Y_{W_{2}}(x_{2}^{L_{V}(0)}v, x_{2})w_{2}\nn
&\quad=\sum_{\nu\in \Gamma(W_{3})}\sum_{m=0}^{n}\binom{-k+n-l-1}{m}\res_{x_{2}}\res_{x_{1}}
\left(\frac{x_{1}-x_{2}}{x_{2}}\right)^{-k+n-l-m-1}
\cdot\nn
&\quad\quad\quad\cdot
 x_{1}^{h_{2}^{\mu}-h_{3}^{\nu}+l}
x_{2}^{-k-2}
\Y^{0}(x_{1}^{L_{W}(0)_{S}}w, x_{1})Y_{W_{2}}(x_{2}^{L_{V}(0)}v, x_{2})w_{2}\nn
&\quad=\sum_{\nu\in \Gamma(W_{3})}\sum_{m=0}^{n}\sum_{r=0}^{N}\sum_{j\in \N}
\binom{-k+n-l-1}{m}\binom{-k+n-l-m-1}{j}(-1)^{j}
\cdot\nn
&\quad\quad\quad\cdot \res_{x_{2}}\res_{x_{1}}x_{1}^{h_{2}^{\mu}-h_{3}^{\nu}-k+n-m-1-j}x_{2}^{-n+l+m-1+j}
\Y^{0}(x_{1}^{L_{W}(0)_{S}}w, x_{1})Y_{W_{2}}(x_{2}^{L_{V}(0)}v, x_{2})w_{2}.
\end{align}
Since $w_{2}\in (W_{2})_{[h_{2}^{\mu}+l]}$, 
$\res_{x}x^{q-1}Y_{W_{2}}(x^{L_{V}(0)}v, x)\in 
(W_{2})_{[h_{2}^{\mu}+l-q]}$ and thus is equal to $0$ if $q>l$.
In the case $j>n-m$, we have $-n+l+m+j>l$ and hence 
$$\res_{x_{2}}x_{2}^{-n+l+m-1+j}Y_{W_{2}}(x_{2}^{L_{V}(0)}v, x_{2})w_{2}=0.$$
In particular,  those terms in the right-hand side of (\ref{l-r-action-comm-7}) with $j>n-m+r$ is $0$.
So the right-hand side of (\ref{l-r-action-comm-7}) is equal to 
\begin{align}\label{l-r-action-comm-10}
&\sum_{\nu\in \Gamma(W_{3})}\sum_{m=0}^{n}\sum_{j=0}^{n-m}
\binom{-k+n-l-1}{m}\binom{-k+n-l-m-1}{j}(-1)^{j}
\cdot\nn
&\quad\quad\quad\cdot \res_{x_{2}}\res_{x_{1}}x_{1}^{h_{2}^{\mu}-h_{3}^{\nu}-k+n-m-1-j}x_{2}^{-n+l+m-1+j}
\Y^{0}(x_{1}^{L_{W}(0)_{S}}w, x_{1})Y_{W_{2}}(x_{2}^{L_{V}(0)}v, x_{2})w_{2}\nn
&\quad=\sum_{\nu\in \Gamma(W_{3})}\sum_{m=0}^{n}\sum_{q=m}^{n}
\binom{-k+n-l-1}{m}\binom{-k+n-l-m-1}{q-m}(-1)^{q-m}
\cdot\nn
&\quad\quad\quad\cdot \res_{x_{2}}\res_{x_{1}}x_{1}^{h_{2}^{\mu}-h_{3}^{\nu}-k+n-1-q}x_{2}^{-n+l-1+q}
\Y^{0}(x_{1}^{L_{W}(0)_{S}}w, x_{1})Y_{W_{2}}(x_{2}^{L_{V}(0)}v, x_{2})w_{2}\nn
&\quad=\sum_{\nu\in \Gamma(W_{3})}\sum_{q=0}^{n}\left(\sum_{m=0}^{q}
\binom{-k+n-l-1}{m}\binom{-k+n-l-m-1}{q-m}(-1)^{q-m}\right)
\cdot\nn
&\quad\quad\quad\cdot \res_{x_{2}}\res_{x_{1}}x_{1}^{h_{2}^{\mu}-h_{3}^{\nu}-k+n-1-q}x_{2}^{-n+l-1+q}
\Y^{0}(x_{1}^{L_{W}(0)_{S}}w, x_{1})Y_{W_{2}}(x_{2}^{L_{V}(0)}v, x_{2})w_{2}.
\end{align}
Using (\ref{binom-identity}), we see that the right-hand side of (\ref{l-r-action-comm-10}) is equal to 
\begin{align}\label{l-r-action-comm-11}
&\sum_{\nu\in \Gamma(W_{3})}\res_{x_{2}}\res_{x_{1}}x_{1}^{h_{2}^{\mu}-h_{3}^{\nu}-k+n-1}x_{2}^{-n+l-1}
\Y^{0}(x_{1}^{L_{W}(0)_{S}}w, x_{1})Y_{W_{2}}(x_{2}^{L_{V}(0)}v, x_{2})w_{2}\nn
&\quad=\vartheta_{\Y}([w]_{kn})
\res_{x_{2}}x_{2}^{-n+l-1}Y_{W_{2}}(x_{2}^{L_{V}(0)}v, x_{2})w_{2}\nn
&\quad=\vartheta_{\Y}([w]_{kn})\vartheta_{W_{2}}([v]_{nl})w_{2}.
\end{align}

The calculations from (\ref{l-r-action-comm-5}) to (\ref{l-r-action-comm-11}) gives
$$\vartheta_{\Y}([w]_{kn}\diamond [v]_{nl})w_{2}
=\vartheta_{\Y}([w]_{kn})\vartheta_{W_{2}}([v]_{nl})w_{2}.$$
Thus we have proved 
$$\vartheta_{\Y}([w]_{km}\diamond [v]_{nl})w_{2}=
\vartheta_{\Y}([w]_{km})\vartheta_{W_{2}}([v]_{nl})w_{2}$$
for $k, m, n, l\in \N$, $v\in V$, $w\in W$ and $w_{2}\in W_{2}$.  
This shows that $\vartheta_{\Y}$ commutes with the right actions of $U^{\infty}(V)$.

We know that $U^{\infty}(W)/\ker \vartheta_{\Y}$ is linearly isomorphic 
to  the image $\vartheta_{\Y}(U^{\infty}(W))$ of  $U^{\infty}(W)$ 
in $\hom(W_{2}, W_{3})$ under $\vartheta_{\Y}$. 
Since $\vartheta_{\Y}$ commutes with both the left and right actions of $U^{\infty}(V)$,
$U^{\infty}(W)/\ker \vartheta_{\Y}$ is in fact equivalent to $\vartheta_{\Y}(U^{\infty}(W))$ 
as $U^{\infty}(V)$-bimodules. But $\hom( W_{2}, W_{3})$ 
is an $A^{\infty}(V)$-bimodule and $\vartheta_{\Y}(U^{\infty}(W))$ is an $A^{\infty}(V)$-subbimodule
of $\hom(W_{2}, W_{3})$. 
Thus $U^{\infty}(W)/\ker \vartheta_{\Y}$ is an $A^{\infty}(V)$-bimodule equivalent to 
the $A^{\infty}(V)$-bimodule $\vartheta_{\Y}(U^{\infty}(W))$.
\epfv

Let $Q^{\infty}(W)$ be the intersection of $\ker \vartheta_{\Y}$ for all 
lower-bounded generalized $V$-modules $W_{2}$ and $W_{3}$ and all intertwining operators $\Y$ of type 
$\binom{W_{3}}{WW_{2}}$. 

\begin{thm}
The quotient $A^{\infty}(W)=U^{\infty}(W)/Q^{\infty}(W)$ is an $A^{\infty}(V)$-bimodule. 
\end{thm}
\pf
For every pair of  
lower-bounded generalized $V$-modules $W_{2}$ and $W_{3}$ and every intertwining operator $\Y$ of type 
$\binom{W_{3}}{WW_{2}}$, 
$U^{\infty}(W)/\ker \vartheta_{\Y}$ is an $A^{\infty}(V)$-bimodule. Hence 
$U^{\infty}(W)$ is an $A^{\infty}(V)$-bimodule modulo elements of $\ker \vartheta_{\Y}$ for all 
lower-bounded generalized $V$-modules $W_{2}$ and $W_{3}$  and all  intertwining operator $\Y$ of type 
$\binom{W_{3}}{WW_{2}}$.
Thus $U^{\infty}(W)$ is an $A^{\infty}(V)$-bimodule modulo elements of the intersection 
$Q^{\infty}(W)$ of all such $\ker \vartheta_{\Y}$. This shows that $A^{\infty}(W)
=U^{\infty}(W)/Q^{\infty}(W)$ is an 
$A^{\infty}(V)$-bimodule. 
\epfv

\begin{prop}\label{alg-bimod}
In the case $W=V$, the $A^{\infty}(V)$-bimodule $A^{\infty}(W)$ is 
equal to the associative algebra $A^{\infty}(V)$ introduced in \cite{H-aa-va}
(see also Section 2).
\end{prop}
\pf
We need only prove that $Q^{\infty}(W)$ for $W=V$ as a lower-bounded generalized $V$-module
 is equal to $\Q^{\infty}(V)$ 
in \cite{H-aa-va} (see also Section 2). 
By Remark \ref{alt-defn-Q-infty-V}, $Q^{\infty}(V)$ in \cite{H-aa-va} 
is equal to the intersection of 
$\ker \vartheta_{W}$ for all lower-bounded generalized $V$-modules $W$. 
Note that $\ker \vartheta_{W_{2}}=\ker \vartheta_{Y_{W_{2}}}$.
If for every pair of lower-bounded generalized $V$-modules 
$W_{2}$ and $W_{3}$ and every intertwining operator $\Y$ of type $\binom{W_{3}}{VW_{2}}$,
$\ker \vartheta_{W_{2}}\subset \ker \vartheta_{\Y}$, then the intersection of 
$\ker \vartheta_{W_{2}}$ for all lower-bounded generalized $V$-module $W_{2}$ 
is equal to the intersection of all $\ker \vartheta_{\Y}$ for all 
lower-bounded generalized $V$-modules 
$W_{2}$ and $W_{3}$ and all intertwining operator $\Y$ of type $\binom{W_{3}}{VW_{2}}$,
that is, the proposition is true.

We now prove this fact. From the $L(-1)$-derivative property and $L_{V}(-1)\one=0$, 
we see that $\Y(\one, x)$ must be independent of $x$. We denote it as $f$ and then $f$ is a linear map
from $W_{2}$ to $W_{3}$. It is easy to verify that $f$ is in fact a $V$-module map. 

Using the associativity between the intertwining operator $\Y$ and the vertex operator
maps $Y_{V}$ and $Y_{W_{2}}$, we have 
\begin{align}\label{alg-bimod-1}
\langle w_{3}', \Y(v, z_{2})w_{2}\rangle&=\langle w_{3}', \Y(Y_{V}(\one, z_{1}-z_{2})v, z_{2})w_{2}\rangle\nn
&=\langle w_{3}', \Y(\one, z_{1})Y_{W_{2}}(v, z_{2})w_{2}\rangle\nn
&=\langle w_{3}', f(Y_{W_{2}}(v, z_{2})w_{2})\rangle
\end{align}
for $v\in V$, $w_{2}\in W_{2}$ and $w_{3}'\in W_{3}'$. Since every term in (\ref{alg-bimod-1}) 
is defined for all $z_{1}$ and $z_{2}\ne 0$, (\ref{alg-bimod-1})  holds for all such $z_{1}$ and $z_{2}$. 
In particular, we obtain 
$$\Y(v, x)w_{2}=f(Y_{W_{2}}(v, x)w_{2})$$
for $v\in V$ and $w_{2}\in W_{2}$. By the definition of $\vartheta_{W_{2}}$ and $\vartheta_{\Y}$
we have 
\begin{align*}
\vartheta_{\Y}([v]_{kl})w_{2}&=\res_{x}x^{l-k-1}\Y(x^{L_{V}(0)}v, x)w_{2}\nn
&=\res_{x}x^{l-k-1}
f(Y_{W_{2}}(x^{L_{V}(0)}v, x)w_{2})\nn
&=f(\vartheta_{W_{2}}([v]_{kl})w_{2})
\end{align*}
for $k, l\in \N$, $v\in V$ and $w_{2}\in (W_{2})_{\l l\r}$. Then we have 
$$\vartheta_{\Y}(\mathfrak{v})=f\circ \vartheta_{W_{2}}(\mathfrak{v})$$
for $\mathfrak{v}\in U^{\infty}(V)$. Now $\ker \vartheta_{W_{2}}\subset \ker \vartheta_{\Y}$
follows immediately. 
\epfv

Using the Jacobi identity  of intertwining operators, we have the following 
result giving some particular elements of $Q^{\infty}(W)$:

\begin{prop}\label{jacobi-ker-theta}
For $k, l, n\in \N$, $p\in \Z$ such that $l+p\in \N$,
homogeneous $v\in V$ and $w\in W$, the element
\begin{align}\label{jacobi-ker-theta-0}
&\sum_{j\in \N}(-1)^{j}\binom{p}{j}[v]_{k,n+p-j}\diamond
[w]_{n+p-j,l+p}\nn
&\quad -\sum_{\stackrel{j\in \N}{l-n+k+p-j\ge 0}}(-1)^{p-j}\binom{p}{j}
[w]_{k,l-n+k+p-j}\diamond [v]_{l-n+k+p-j,l+p}\nn
&\quad-\sum_{j\in \N}\binom{\wt v+n-k-1}{j}
[(Y_{W})_{p+j}(v)w]_{k,l+p}
\end{align}
of $U^{\infty}(W)$ is in fact  in $Q^{\infty}(W)$.
\end{prop}
\pf
Let $W_{2}$ and $W_{3}$ be lower-bounded generalized $V$-modules and 
$\Y$ an intertwining operator of type $\binom{W_{3}}{WW_{2}}$. 
For $k, l, n\in \N$, $p\in \Z$ such that $l+p\in \N$, $\mu\in \Gamma(W_{2})$, 
 homogeneous $v\in V$, $w\in W$ and $w_{2}\in (W_{2})_{[h_{2}^{\mu}+l+p]}$,
\begin{align}\label{jacobi-ker-theta-1}
&\sum_{j\in \N}(-1)^{j}\binom{p}{j}\vartheta_{\Y}([v]_{k,n+p-j}\diamond
[w_{1}]_{n+p-j,l+p}) w_{2}\nn
&\quad =\sum_{j\in \N}(-1)^{j}\binom{p}{j}\vartheta_{W_{3}}([v]_{k,n+p-j})
\vartheta_{\Y}([w]_{n+p-j,l+p})w_{2}\nn
&\quad =\sum_{\nu\in \Gamma(W_{3})}
\sum_{j\in \N}(-1)^{j}\binom{p}{j}\vartheta_{W_{3}}([v]_{k,n+p-j})\cdot\nn
&\quad\quad\quad\quad\quad\quad \cdot 
\text{Coeff}_{\log x_{2}}^{0}\res_{x_{2}}x_{2}^{h_{2}^{\mu}-h_{3}^{\nu}+(l+p)-(n+p-j)-1}
\Y(x_{2}^{L_{W}(0)}w, x_{2})w_{2}\nn
&\quad =\sum_{\nu\in \Gamma(W_{3})}\sum_{j\in \N}\binom{p}{j}(-1)^{j}
\res_{x_{1}}\text{Coeff}_{\log x_{2}}^{0}
\res_{x_{2}}x_{1}^{(n+p-j)-k-1}\cdot\nn
&\quad\quad\quad\quad\quad\quad \cdot 
x_{2}^{h_{2}^{\mu}-h_{3}^{\nu}+(l+p)-(n+p-j)-1}
Y_{W_{3}}(x_{1}^{L_{V}(0)}v, x_{1})
\Y(x_{2}^{L_{W}(0)}w, x_{2})w_{2}\nn
&\quad =\sum_{\nu\in \Gamma(W_{3})}\res_{x_{1}}\text{Coeff}_{\log x_{2}}^{0}
\res_{x_{2}}(x_{1}-x_{2})^{p}\cdot\nn
&\quad\quad\quad\quad\quad\quad \cdot 
x_{1}^{n-k-1}x_{2}^{h_{2}^{\mu}-h_{3}^{\nu}+l-n-1}
Y_{W_{3}}(x_{1}^{L_{V}(0)}v, x_{1})
\Y(x_{2}^{L_{W}(0)}w, x_{2})w_{2}\nn
&\quad =\sum_{\nu\in \Gamma(W_{3})}\res_{x_{1}}\text{Coeff}_{\log x_{2}}^{0}
\res_{x_{2}}\res_{x_{0}}x_{0}^{p}x_{0}^{-1}\delta\left(\frac{x_{1}-x_{2}}{x_{0}}\right)\cdot\nn
&\quad\quad\quad\quad\quad\quad \cdot 
x_{1}^{n-k-1}x_{2}^{h_{2}^{\mu}-h_{3}^{\nu}+l-n-1}
Y_{W_{3}}(x_{1}^{L_{V}(0)}v, x_{1})
\Y(x_{2}^{L_{W}(0)}w, x_{2})w_{2}.
\end{align}
Using the Jacobi identity for the intertwining operator $\Y$, we see that the right-hand side of 
(\ref{jacobi-ker-theta-1}) is equal to
\begin{align}\label{jacobi-ker-theta-2}
&\sum_{\nu\in \Gamma(W_{3})}\res_{x_{1}}\text{Coeff}_{\log x_{2}}^{0}
\res_{x_{2}}\res_{x_{0}}x_{0}^{p}x_{0}^{-1}\delta\left(\frac{x_{2}-x_{1}}{-x_{0}}\right)\cdot\nn
&\quad\quad\quad\quad\quad\quad \cdot x_{1}^{n-k-1}x_{2}^{h_{2}^{\mu}-h_{3}^{\nu}+l-n-1}
\Y(x_{2}^{L_{W}(0)}w, x_{2})Y_{W_{2}}(x_{1}^{L_{V}(0)}v, x_{1})w_{2}\nn
&\quad \quad +\sum_{\nu\in \Gamma(W_{3})}\res_{x_{1}}\text{Coeff}_{\log x_{2}}^{0}
\res_{x_{2}}\res_{x_{0}}x_{0}^{p}x_{1}^{-1}\delta\left(\frac{x_{2}+x_{0}}{x_{1}}\right)\cdot\nn
&\quad\quad\quad\quad\quad\quad \cdot 
x_{1}^{n-k-1}x_{2}^{h_{2}^{\mu}-h_{3}^{\nu}+l-n-1}
\Y(Y_{W_{1}}(x_{1}^{L_{V}(0)}v, x_{0})x_{2}^{L_{W}(0)}w, x_{2})w_{2}\nn
&\quad =\sum_{\nu\in \Gamma(W_{3})}\res_{x_{1}}\text{Coeff}_{\log x_{2}}^{0}
\res_{x_{2}}(-x_{2}+x_{1})^{p}\cdot\nn
&\quad\quad\quad\quad\quad\quad \cdot x_{1}^{n-k-1}x_{2}^{h_{2}^{\mu}-h_{3}^{\nu}+l-n-1}
\Y(x_{2}^{L_{W}(0)}w, x_{2})Y_{W_{2}}(x_{1}^{L_{V}(0)}v, x_{1})w_{2}\nn
&\quad \quad +\sum_{\nu\in \Gamma(W_{3})}\res_{x_{1}}\text{Coeff}_{\log x_{2}}^{0}
\res_{x_{2}}\res_{x_{0}}x_{0}^{p}x_{1}^{-1}\delta\left(\frac{x_{2}+x_{0}}{x_{1}}\right)\cdot\nn
&\quad\quad\quad\quad\quad\quad \cdot 
x_{1}^{n-k-1}x_{2}^{h_{2}^{\mu}-h_{3}^{\nu}+l-n-1}
\Y(x_{2}^{L_{W}(0)}Y_{W_{1}}(x_{1}^{L_{V}(0)}x_{2}^{-L_{V}(0)}v, x_{0}x_{2}^{-1})w, x_{2})w_{2}\nn
&\quad =\sum_{\nu\in \Gamma(W_{3})}\sum_{j\in \N}(-1)^{p-j}\binom{p}{j}\text{Coeff}_{\log x_{2}}^{0}
\res_{x_{2}} x_{2}^{h_{2}^{\mu}-h_{3}^{\nu}+(l-n+k+p-j)-k-1}\cdot\nn
&\quad\quad\quad\quad\quad\quad \cdot 
\Y(x_{2}^{L_{W}(0)}w, x_{2})\res_{x_{1}}x_{1}^{(l+p)-(l-n+k+p-j)-1}
Y_{W_{2}}(x_{1}^{L_{V}(0)}v, x_{1})w_{2}\nn
&\quad \quad +\sum_{\nu\in \Gamma(W_{3})}\res_{x_{1}}\text{Coeff}_{\log x_{2}}^{0}
\res_{x_{2}}\res_{x_{0}}x_{0}^{p}x_{1}^{-1}\delta\left(\frac{x_{2}+x_{0}}{x_{1}}\right)\cdot\nn
&\quad\quad\quad\quad\quad\quad \cdot 
x_{1}^{n-k-1}x_{2}^{h_{2}^{\mu}-h_{3}^{\nu}+l-n-1}
\Y(x_{2}^{L_{W}(0)}Y_{W_{1}}(x_{1}^{L_{V}(0)}x_{2}^{-L_{V}(0)}v, x_{0}x_{2}^{-1})w, x_{2})w_{2}.
\end{align}
Using the definitions of $\vartheta_{\Y}$ and $\vartheta_{W_{2}}$, 
the property of the formal $\delta$-function, noting that $v$ is homogeneous and 
$\res_{x_{1}}x_{1}^{(l+p)-(l-n+k+p-j)-1}
Y_{W_{2}}(x_{1}^{L_{V}(0)}v, x_{1})w_{2}=0$ when  $l-n+k+p-j<0$,
and changing the variable $x_{0}$ to $y=x_{0}x_{2}^{-1}$ in the second term, 
we see that the right-hand side of (\ref{jacobi-ker-theta-2}) is equal to 
\begin{align}\label{jacobi-ker-theta-3}
&\sum_{\stackrel{j\in \N}{l-n+k+p-j\ge 0}}(-1)^{p-j}\binom{p}{j}
\vartheta_{\Y}( [w]_{k,l-n+k+p-j})\cdot\nn
&\quad\quad\quad\quad\quad\quad\quad \cdot
 \res_{x_{1}}x_{1}^{(l+p)-(l-n+k+p-j)-1}Y_{W_{2}}(x_{1}^{L_{V}(0)}v, x_{1})w_{2}\nn
&\quad\quad +\sum_{\nu\in \Gamma(W_{3})}\text{Coeff}_{\log x_{2}}^{0}
\res_{x_{2}}\res_{x_{0}}x_{0}^{p}
(x_{2}+x_{0})^{n-k-1}\cdot\nn
&\quad\quad\quad\quad\quad\quad\quad \cdot x_{2}^{h_{2}^{\mu}-h_{3}^{\nu}+l-n-1}
\Y(x_{2}^{L_{W}(0)}Y_{W_{1}}((1+x_{0}x_{2}^{-1})^{L_{V}(0)}v, x_{0}x_{2}^{-1})w, x_{2})w_{2}\nn
& \quad=\sum_{\stackrel{j\in \N}{l-n+k+p-j\ge 0}}(-1)^{p-j}\binom{p}{j}
\vartheta_{\Y}( [w]_{k,l-n+k+p-j})
 \vartheta_{W_{2}}([v]_{l-n+k+p-j,l+p})w_{2}\nn
&\quad \quad+\sum_{\nu\in \Gamma(W_{3})}\text{Coeff}_{\log x_{2}}^{0}
\res_{x_{2}}\res_{x_{0}}x_{0}^{p}
(1+x_{0}x_{2}^{-1})^{\swt v+n-k-1}\cdot\nn
&\quad\quad\quad\quad\quad\quad\quad \cdot x_{2}^{h_{2}^{\mu}-h_{3}^{\nu}+l-k-2}
\Y(x_{2}^{L_{W}(0)}Y_{W_{1}}(v, x_{0}x_{2}^{-1})w, x_{2})w_{2}\nn
&\quad =\sum_{\stackrel{j\in \N}{l-n+k+p-j\ge 0}}(-1)^{p-j}\binom{p}{j}
\vartheta_{\Y}( [w]_{k,l-n+k+p-j}\diamond [v]_{l-n+k+p-j,l+p})w_{2}\nn
&\quad \quad+\sum_{\nu\in \Gamma(W_{3})}\text{Coeff}_{\log x_{2}}^{0}
\res_{x_{2}}\res_{y}y^{p}
(1+y)^{\swt v+n-k-1}\cdot\nn
&\quad\quad\quad\quad\quad\quad\quad \cdot x_{2}^{h_{2}^{\mu}-h_{3}^{\nu}+(l+p)-k-1}
\Y(x_{2}^{L_{W}(0)}Y_{W_{1}}(v, y)w, x_{2})w_{2}.
\end{align}
Expanding $(1+y)^{\swt v+n-k-1}$ and $Y_{W_{1}}(v, y)$ and evaluating $\res_{y}$, 
we see that the second term in the right-hand side
of (\ref{jacobi-ker-theta-3}) is equal to
\begin{align}\label{jacobi-ker-theta-4}
&\sum_{\nu\in \Gamma(W_{3})}\sum_{j\in \N}\binom{\wt v+n-k-1}{j}\text{Coeff}_{\log x_{2}}^{0}
\res_{x_{2}}x_{2}^{h_{2}^{\mu}-h_{3}^{\nu}+(l+p)-k-1}
\Y(x_{2}^{L_{W}(0)}(Y_{W_{1}})_{p+j}(v)w, x_{2})w_{2}\nn
&\quad =\sum_{j\in \N}\binom{\wt v+n-k-1}{j}
\vartheta_{\Y}([(Y_{W_{1}})_{p+j}(v)w]_{k,l+p})
w_{2}.
\end{align}
From  (\ref{jacobi-ker-theta-1})--(\ref{jacobi-ker-theta-4}), we obtain
\begin{align}\label{jacobi-ker-theta-0.5}
&\vartheta_{\Y}\Biggl(\sum_{j\in \N}(-1)^{j}\binom{p}{j}[v]_{k,n+p-j}\diamond
[w_{1}]_{n+p-j,l+p}\nn
&\quad \quad-\sum_{\stackrel{j\in \N}{l-n+k+j\ge 0}}(-1)^{p-j}\binom{p}{j}
[w_{1}]_{k,l-n+k+j}\diamond [v]_{l-n+k+j,l+p}\nn
&\quad\quad -\sum_{j\in \N}\binom{\wt v+n-k-1}{j}
[(Y_{W_{1}})_{j}(v)w_{1}]_{k,l+p}\Biggr)w_{2}\nn
&\quad =0.
\end{align}
Since $w_{2}$, $l$ and $p$ are arbitrary, we see from 
 (\ref{jacobi-ker-theta-0.5}) that 
 (\ref{jacobi-ker-theta-0}) is in $\ker \vartheta_{\Y}$. Since $W_{2}, W_{3}$
and $\Y$ are arbitrary, (\ref{jacobi-ker-theta-0}) is in the intersection $Q^{\infty}(W)$ of all $\ker \vartheta_{\Y}$.
\epfv

\begin{rema}\label{L(0)-L(-1)-comm}
{\rm In the case that $V$ is a vertex operator algebra, that is, we have a conformal vector $\omega\in V$. 
Let $\omega^{\infty}(0)$ and $\omega^{\infty}(-1)$ be the elements of $U^{\infty}(V)$ with diagonal entries being $\omega \in V$ and
all the other entries being $0$ and with the $(k+1, k)$-entries being $\omega$ for $k\in \N$ and 
all the other entries being $0$, respectively. Then 
$$\omega^{\infty}(0)=\sum_{k\in \N}[\omega]_{kk}, \;\omega^{\infty}(-1)=\sum_{k\in \N}[\omega]_{k+1, k}.$$
From the definitions of the left and right actions of $U^{\infty}(V)$ on $U^{\infty}(W)$,
we have $[\omega]_{kk}\diamond [w]_{kl}=\omega^{\infty}(0)\diamond [w]_{kl}$, 
$[w]_{kl}\diamond [\omega]_{ll}=[w]_{kl}\diamond \omega^{\infty}(0)$,  
$[\omega]_{k+1, k}\diamond [w]_{kl}=\omega^{\infty}(-1)\diamond [w]_{kl}$ and
$[w]_{kl}\diamond [\omega]_{l+1, l}=[w]_{kl}\diamond \omega^{\infty}(-1)$ for $k, l\in \N$ and $w\in W$. 
Taking $n=k$, $p=0$ and $v=\omega$ in (\ref{jacobi-ker-theta-0}), 
we see that elements of the form
\begin{align*}
&\omega^{\infty}(0)\diamond [w]_{nl}-[w]_{nl}\diamond \omega^{\infty}(0)
-[(L_{W}(-1)+L_{W}(0))w]_{nl}\nn
&\quad =[\omega]_{nn}\diamond [w]_{nl}-[w]_{nl}\diamond [\omega]_{ll}
-[(L_{W}(-1)+L_{W}(0))w]_{nl}
\end{align*}
for $n, l\in \N$ and $w\in W$ are in $Q^{\infty}(W)$. Taking $k=n+1$, $p=0$ and $v=\omega$ in (\ref{jacobi-ker-theta-0}),
we see that elements of the form 
\begin{align*}
&\omega^{\infty}(-1)\diamond [w]_{nl}-[w]_{n+1, l+1}\diamond \omega^{\infty}(-1)
-[L_{W}(-1)w]_{n+1, l}\nn
&\quad =[\omega]_{n+1, n}\diamond [w]_{nl}-[w]_{n+1, l+1}\diamond [\omega]_{l+1, l}
-[L_{W}(-1)w]_{n+1, l}
\end{align*}
for $n, l\in \N$ and $w\in W$ are in $Q^{\infty}(W)$. In particular, we have 
\begin{align*}
&(\omega^{\infty}(0)+Q^{\infty}(W))\diamond ([w]_{nl}+Q^{\infty}(W))-([w]_{nl}+Q^{\infty}(W))
\diamond (\omega^{\infty}(0)+Q^{\infty}(W))\nn
&\quad =[(L_{W}(-1)+L_{W}(0))w]_{nl}+Q^{\infty}(W),\\%\label{lt-rt-L(0)-comm}\\
&(\omega^{\infty}(-1)+Q^{\infty}(W))\diamond ([w]_{nl}+Q^{\infty}(W))-([w]_{n+1, l+1}+Q^{\infty}(W))
\diamond (\omega^{\infty}(-1)+Q^{\infty}(W))\nn
&\quad =[L_{W}(-1)w]_{n+1, l}+Q^{\infty}(W)%\label{lt-rt-L(-1)-comm}
\end{align*}
in $A^{\infty}(W)$. } 
\end{rema}

In the special case that $W=V$, we have:

\begin{cor}\label{jacobi-ker-theta-V}
For $k, l, n\in \N$, $p\in \Z$ such that $l+p\in \N$, homogeneous $u\in V$ and $v\in V$, 
the elements 
\begin{align}\label{jacobi-ker-theta-0-V}
&\sum_{j\in \N}(-1)^{j}\binom{p}{j}[u]_{k,n+p-j}\diamond
[v]_{n+p-j,l+p}\nn
&\quad -\sum_{\stackrel{j\in \N}{l-n+k+p-j\ge 0}}(-1)^{p-j}\binom{p}{j}
[v]_{k,l-n+k+p-j}\diamond [u]_{l-n+k+p-j,l+p}\nn
&\quad-\sum_{j\in \N}\binom{\wt u+n-k-1}{j}
[(Y_{V})_{p+j}(u)v]_{k,l+p}
\end{align}
is in $Q^{\infty}(V)$.
\end{cor}
\pf
This corollary follows immediately from Propositions \ref{alg-bimod} and 
\ref{jacobi-ker-theta}. 
\epfv

\begin{rema}
{\rm We can also give a subspace $O^{\infty}(W)$ of $Q^{\infty}(W)$
analogous to $O^{\infty}(V)$ in \cite{H-aa-va}. 
We conjecture that $Q^{\infty}(W)$ is spanned 
by $O^{\infty}(W)$
and the elements of the form (\ref{jacobi-ker-theta-0})
for $k, l, n, p\in \N$, homogeneous $v\in V$ and $w\in W$. In particular, 
we also conjecture that $Q^{\infty}(V)$ is spanned 
by $O^{\infty}(V)$ in \cite{H-aa-va} 
and elements of the form (\ref{jacobi-ker-theta-0-V}) of $Q^{\infty}(V)$. 
We shall
discuss these in a future paper. }
\end{rema}

\renewcommand{\theequation}{\thesection.\arabic{equation}}
\setcounter{equation}{0}
\setcounter{thm}{0}
\section{Isomorphisms between spaces of intertwining operators and $A^{\infty}(V)$-module maps}

We formulate and prove the first main theorem of the present paper in this section. For lower-bounded generalized  $V$-modules
$W_{1}$, $W_{2}$ and $W_{3}$, we
define a linear map $\rho: \mathcal{V}_{W_{1}W_{2}}^{W_{3}}\to 
\hom_{A^{\infty}(V)}(A^{\infty}(W_{1})\otimes_{A^{\infty}(V)} W_{2}, W_{3})$,
where $W_{1}$, $W_{2}$ and $W_{3}$ are lower-bounded generalized $V$-module 
and $\mathcal{V}_{W_{1}W_{2}}^{W_{3}}$ is the space of 
intertwining operators of type $\binom{W_{3}}{W_{1}W_{2}}$.
Our first main theorem states that $\rho$ is an isomorphism. 

Before we formulate and prove this theorem, we prove first 
that the category lower-bounded generalized $V$-modules
 and the category of graded $A^{\infty}(V)$-modules are isomorphic 
(not just equivalent since the underlying vector space are the same). 
In Section 2, we have obtained a functor from the category of lower-bounded generalized $V$-modules
to the category of graded $A^{\infty}(V)$-modules. We now have:

\begin{thm}\label{cat-isom}
The functor from the category lower-bounded generalized $V$-modules
to the category of graded $A^{\infty}(V)$-modules is in fact an isomorphism
of categories. 
\end{thm}
\pf
Given a graded $A^{\infty}(V)$-module $W$, we need to construct a
lower-bounded generalized $V$-module structure on $W$. We define 
$$(Y_{W})_{\swt v+l-k-1}(v)w=\vartheta_{W}([v]_{kl})w$$
for $k, l\in \N$, homogeneous $v\in V$ and $w\in W_{\l l\r}$ and 
$$Y_{W}(v, x)w=\sum_{k\in \N}(Y_{W})_{\swt v+l-k-1}(v)w x^{-\swt v-l+k}$$
for $l\in \N$, homogeneous $v\in V$ and $w\in W_{\l l\r}$.
It is clear that $Y_{W}(v, x)w$ is lower truncated. For any lower-bounded generalized 
$V$-module $W_{0}$, by definition, 
$$\vartheta_{W_{0}}([\one]_{kl})w_{0}=\res_{x}x^{l-k-1}w_{0}=\delta_{kl}[\one]_{ll}w_{0}$$
for $k, l\in \N$ and $w_{0}\in W_{0}$. So we have $[\one]_{kl}-\delta_{kl}[\one]_{ll}
\in \ker \vartheta_{W_{0}}$. Since $W_{0}$ is arbitrary, 
by Remark \ref{alt-defn-Q-infty-V}, we see that 
$[\one]_{kl}-\delta_{kl}[\one]_{ll}
\in Q^{\infty}(V)$. In particular, $\vartheta_{W}([\one]_{kl})w=\delta_{kl}w$
for $k, l\in \N$ and $w\in W_{\l l\r}$. Thus we have 
$Y_{W}(\one, x)w=w$. The $L(0)$- and $L(-1)$-commutator formulas and the $L(0)$-grading 
condition follow 
hold since $W$ is a graded $A^{\infty}(V)$-algebra. The only remaining axiom to be proved is the 
Jacobi identity. 

We now prove the Jacobi identity using Corollary \ref{jacobi-ker-theta-V}.
From Corollary \ref{jacobi-ker-theta-V}, we know that (\ref{jacobi-ker-theta-0-V}) with  $w$ replaced by $w_{1}$ 
is an element of $Q^{\infty}(V)$. 
In particular, 
\begin{align}\label{cat-isom-1}
&\sum_{\in \N}(-1)^{j}\binom{p}{j}\vartheta_{W}([u]_{k,n+p-j}\diamond
[v]_{n+p-j,l+p})w_{2}\nn
&\quad \quad -\sum_{\stackrel{j\in \N}{l-n+k+p-j\ge 0}}(-1)^{p-j}\binom{p}{j}
\vartheta_{W}([v]_{k,l-n+k+p-j}\diamond [u]_{l-n+k+p-j,l+p})w\nn
&\quad \quad-\sum_{j\in \N}\binom{\wt u+n-k-1}{j}
\vartheta_{W}([(Y_{V})_{p+j}(u)v]_{k,l+p})w\nn
&\quad =0
\end{align}
for $k, l, n\in \N$, $p\in \Z$ such that $l+p\in \N$, homogeneous $u, v\in V$ and $w\in W_{\l l+p\r}$. 
The left-hand side of (\ref{cat-isom-1}) is equal to 
\begin{align}\label{cat-isom-2}
&\sum_{j\in \N}(-1)^{j}\binom{p}{j}\vartheta_{W}([u]_{k,n+p-j})
\vartheta_{W}([v]_{n+p-j,l+p})w\nn
&\quad -\sum_{\stackrel{j\in \N}{l-n+k+p-j\ge 0}}(-1)^{p-j}\binom{p}{j}
\vartheta_{W}([v]_{k,l-n+k+p-j})\vartheta_{W}([v]_{l-n+k+p-j,l+p})w_{2}\nn
&\quad-\sum_{j\in \N}\binom{\wt u+n-k-1}{j}
\vartheta_{W}([(Y_{V})_{p+j}(u)v]_{k,l+p})w.
\end{align}
From the definitions of $\vartheta_{W}$, 
we see that (\ref{cat-isom-2}) is equal to 
\begin{align}\label{cat-isom-3}
&\sum_{j\in \N}(-1)^{j}\binom{p}{j}\res_{x_{1}}\res_{x_{2}}x_{1}^{n+p-j-k-1}
x_{2}^{l-n+j-1}Y_{W}(x_{1}^{L_{V}(0)}u, x_{1})
Y_{W}(x_{2}^{L_{V}(0)}v, x_{2})w\nn
&\quad -\sum_{\stackrel{j\in \N}{l-n+k+p-j\ge 0}}(-1)^{p-j}\binom{p}{j}
\res_{x_{1}}\res_{x_{2}}x_{1}^{n-k+j-1} \cdot\nn
&\quad\quad\quad\quad\quad\quad\quad\quad\quad\quad\quad\cdot  x_{2}^{l-n+p-j-1}
Y_{W}(x_{2}^{L_{V}(0)}v, x_{2})Y_{W_{2}}(x_{1}^{L_{V}(0)}u, x_{1})w\nn
&\quad-\sum_{j\in \N}\binom{\wt u+n-k-1}{j}
\res_{x_{2}}x_{2}^{l+p-k-1}Y_{W}(x_{2}^{L_{V}(0)}
(Y_{V})_{p+j}(u)v, x_{2})w\nn
&=\res_{x_{1}}\res_{x_{2}}x_{1}^{n-k-1}
x_{2}^{l-n-1}(x_{1}-x_{2})^{p}Y_{W}(x_{1}^{L_{V}(0)}u, x_{1})
Y_{W}(x_{2}^{L_{V}(0)}v, x_{2})w\nn
&\quad -
\res_{x_{1}}\res_{x_{2}}x_{1}^{n-k-1}x_{2}^{l-n-1}
(-x_{2}+x_{1})^{p}
Y_{W}(x_{2}^{L_{V}(0)}v, x_{2})Y_{W}(x_{1}^{L_{V}(0)}u, x_{1})w\nn
&\quad-
\res_{x_{2}}\res_{y}y^{p}(1+y)^{\swt u+n-k-1}
x_{2}^{l+p-k-1}Y_{W}(x_{2}^{L_{V}(0)}
Y_{V}(u,  y)v, x_{2})w.
\end{align}
Since the left-hand side of (\ref{cat-isom-3}) is equal to the left-hand side of  (\ref{cat-isom-1}), we see that 
the right-hand side of (\ref{cat-isom-3}) is $0$ for $k, l, n\in \N$, $p\in \Z$ such that $l+p\in \N$, homogeneous $u, v\in V$
and $w\in W_{\l l+p\r}$.  In fact, the right-hand side of (\ref{cat-isom-3}) is also $0$ for all $p\in \Z$ 
even if $l+p<0$ 
since $w=0$. Multiplying the right-hand side 
of  (\ref{cat-isom-3}) by $x_{0}^{-p-1}$ and then take sum over $p\in \Z$, we obtain
\begin{align}\label{cat-isom-4}
&\res_{x_{1}}\res_{x_{2}}x_{1}^{n-k-1}
x_{2}^{l-n-1}x_{0}^{-1}\delta\left(\frac{x_{1}-x_{2}}{x_{0}}\right)
 Y_{W}(x_{1}^{L_{V}(0)}u, x_{1})
Y_{W}(x_{2}^{L_{V}(0)}v, x_{2})w\nn
&\quad -
\res_{x_{1}}\res_{x_{2}}x_{1}^{n-k-1}x_{2}^{l-n-1}
x_{0}^{-1}\delta\left(\frac{x_{2}-x_{1}}{-x_{0}}\right)
Y_{W}(x_{2}^{L_{V}(0)}v, x_{2})Y_{W}(x_{1}^{L_{V}(0)}u, x_{1})w\nn
&\quad-
\res_{x_{2}}\res_{y}x_{0}^{-1}\delta\left(\frac{yx_{2}}{x_{0}}\right)(1+y)^{n-k-1}
x_{2}^{l-k-1}Y_{W}(x_{2}^{L_{V}(0)}
Y_{V}((1+y)^{L_{V}(0)}u, y)v, x_{2})w\nn
&=0.
\end{align}
The third term in the left-hand side of (\ref{cat-isom-4}) is equal to 
\begin{align}\label{cat-isom-5}
&-\res_{x_{2}}\res_{y}x_{0}^{-1}\delta\left(\frac{yx_{2}}{x_{0}}\right)(1+y)^{n-k-1}
x_{2}^{l-k-1}  Y_{W}(x_{2}^{L_{V}(0)}
Y_{V}((1+y)^{L_{V}(0)}u, y)v, x_{2})w\nn
&\quad =-\res_{x_{2}}\res_{y}x_{0}^{-1}\delta\left(\frac{yx_{2}}{x_{0}}\right)(1+x_{0}x_{2}^{-1})^{n-k-1}
x_{2}^{l-k-1}\cdot\nn
&\quad\quad\quad\quad\quad\quad\quad\cdot  
Y_{W}(x_{2}^{L_{V}(0)}
Y_{V}((1+x_{0}x_{2}^{-1})^{L_{V}(0)}u, x_{0}x_{2}^{-1})v, x_{2})w\nn
&\quad =-\res_{x_{2}}(x_{2}+x_{0})^{n-k-1}
x_{2}^{l-n-1}
Y_{W}(
Y_{V}((x_{2}+x_{0})^{L_{V}(0)}u, x_{0})x_{2}^{L_{V}(0)}v, x_{2})w\nn
&\quad =-\res_{x_{1}}\res_{x_{2}}x_{1}^{n-k-1}
x_{2}^{l-n-1}x_{1}^{-1}\delta\left(\frac{x_{2}+x_{0}}{x_{1}}\right)
Y_{W}(
Y_{V}(x_{1}^{L_{V}(0)}u, x_{0})x_{2}^{L_{V}(0)}v, x_{2})w.
\end{align}
Using (\ref{cat-isom-5}) and substituting $x_{1}^{-L_{V}(0)}v$ and $x_{2}^{-L_{W_{1}}(0)_{S}}w_{1}$
for $v$ and $w_{1}$, respectively, we see that (\ref{cat-isom-4}) becomes
\begin{align}\label{cat-isom-6}
&\res_{x_{1}}\res_{x_{2}}x_{1}^{n-k-1}
x_{2}^{l-n-1}x_{0}^{-1}\delta\left(\frac{x_{1}-x_{2}}{x_{0}}\right)
Y_{W}(u, x_{1})
Y_{W}(v, x_{2})w\nn
&\quad -
\res_{x_{1}}\res_{x_{2}}x_{1}^{n-k-1}x_{2}^{l-n-1}
x_{0}^{-1}\delta\left(\frac{x_{2}-x_{1}}{-x_{0}}\right)
Y_{W}(v, x_{2})Y_{W_{2}}(u, x_{1})w\nn
&\quad-\res_{x_{1}}\res_{x_{2}}x_{1}^{n-k-1}
x_{2}^{l-n-1}x_{1}^{-1}\delta\left(\frac{x_{2}+x_{0}}{x_{1}}\right)
Y_{W}(
Y_{V}(u, x_{0})v, x_{2})w\nn
&=0.
\end{align}
Since $k, l,n\in \N$ and $w\in W$ are arbitrary, (\ref{main-9}) gives the Jacobi identity for $Y_{W}$. 

This construction of a lower-bounded generalized $V$-module structure from 
a graded $A^{\infty}(V)$-module gives us a functor from 
the category of lower-bounded generalized $V$-modules to the category of 
graded $A^{\infty}(V)$-modules. 
It is clear that the graded $A^{\infty}(V)$-module structure on $W$ obtained from 
$Y_{W}$ is the same as the graded $A^{\infty}(V)$-module structure that we start with.
So this functor is the inverse of the functor from the category of 
graded $A^{\infty}(V)$-modules to the category of lower-bounded generalized $V$-modules.
The proposition is proved. 
\epfv

Let $V$ be a grading-restricted vertex algebra and 
$W_{1}$ , $W_{2}$ and $W_{3}$ lower-bounded generalized $V$-modules.
Let  $A^{\infty}(W_{1})\otimes_{A^{\infty}(V)} W_{2}$ be the tensor product over 
$A^{\infty}(V)$ 
of the $A^{\infty}(V)$-bimodule $A^{\infty}(W_{1})$ and the left $A^{\infty}(V)$-module
$W_{2}$. 

We first prove some lemmas which we shall need later.

\begin{lemma}\label{lemma-4.1}
The $A^{\infty}(V)$-module $A^{\infty}(W_{1})\otimes_{A^{\infty}(V)} W_{2}$ 
is spanned by elements of the form 
$$([w_{1}]_{kl}+Q^{\infty}(W_{1}))
\otimes_{A^{\infty}(V)} w_{2}$$
for $k, l\in \N$, 
$w_{1}\in W_{1}$, $w_{2}\in (W_{2})_{\l l\r}$.
\end{lemma}
\pf
We know that elements of the form 
$$([w_{1}]_{kn}+Q^{\infty}(W_{1}))
\otimes_{A^{\infty}(V)} w_{2}$$ 
for $k, l,  n\in \N$, 
$w_{1}\in W_{1}$ span $A^{\infty}(W_{1})\otimes_{A^{\infty}(V)} W_{2}$.
If $l\ne n$, then by definition
$[w_{1}]_{kn}\diamond [\one]_{ll}=0$.
Also for $ l\in \N$ and $w_{2}\in (W_{2})_{\l l\r}$,
by definition, $\vartheta_{W_{2}}([\one]_{ll})w_{2}=w_{2}$. 
Hence for $k, l,  n\in \N$, 
$w_{1}\in W_{1}$, $w_{2}\in (W_{2})_{\l l\r}$, if $n\ne l$, then
\begin{align*}
&([w_{1}]_{kn}+Q^{\infty}(W_{1}))
\otimes_{A^{\infty}(V)} w_{2}\nn
&\quad =([w_{1}]_{kn}+Q^{\infty}(W_{1}))
\otimes_{A^{\infty}(V)} \vartheta_{W_{2}}([\one]_{ll})w_{2}\nn
&\quad =([w_{1}]_{kn}\diamond [\one]_{ll}+Q^{\infty}(W_{1}))
\otimes_{A^{\infty}(V)} w_{2}\nn
&\quad =0.
\end{align*}
Hence the lemma is true.
\epfv

\begin{lemma}\label{lemma-4.2}
Let 
$$f\in \hom_{A^{\infty}(V)}(A^{\infty}(W_{1})\otimes_{A^{\infty}(V)} 
W_{2}, W_{3}).$$
For $k, l\in \N$, 
$w_{1}\in W_{1}$, $w_{2}\in (W_{2})_{\l l\r}$, 
$$f(([w_{1}]_{kl}+Q^{\infty}(W_{1}))\otimes_{A^{\infty}(V)} w_{2})\in (W_{3})_{\l k\r}.$$
\end{lemma}
\pf
For $k, l\in \N$, 
$w_{1}\in W_{1}$, $w_{2}\in (W_{2})_{\l l\r}$, since $f$ is an 
$A^{\infty}(V)$-module map, 
\begin{align*}
&f(([w_{1}]_{kl}+Q^{\infty}(W_{1}))\otimes_{A^{\infty}(V)} w_{2})\nn
&\quad=f(([\one]_{kk}\diamond [w_{1}]_{kl}
+Q^{\infty}(W_{1}))\otimes_{A^{\infty}(V)} w_{2})\nn
&\quad=\vartheta_{W_{3}}([\one]_{kk})f(([w_{1}]_{kl}
+Q^{\infty}(W_{1}))\otimes_{A^{\infty}(V)} w_{2}).
\end{align*}
By the definition of $\vartheta_{W_{3}}$, $\vartheta_{W_{3}}([\one]_{kk})w_{3}=0$
for $w_{3}\in (W_{3})_{\l n\r}$, $n\ne k$. Hence  the lemma is true.
\epfv

Let $\Y$ be an 
intertwining operator of type $\binom{W_{3}}{W_{1}W_{2}}$. 
We define a linear map 
$$\rho(\Y): A^{\infty}(W_{1})\otimes W_{2} \to W_{3}$$
by 
$$(\rho(\Y))((\mathfrak{w}_{1}+Q^{\infty}(W_{1}))\otimes w_{2})=
\vartheta_{\Y}(\mathfrak{w}_{1})w_{2}$$
for $\mathfrak{w}_{1}\in U^{\infty}(W_{1})$ 
and $w_{2}\in W_{2}$.  Since $Q^{\infty}(W_{1})\subset \ker \vartheta_{\Y}$, 
$\rho(\Y)$ is well defined.

Using Proposition \ref{l-r-action-comm}, we have:

\begin{prop}\label{rho-y-mod-m}
The linear map $\rho(\Y)$ is in fact an $A^{\infty}(V)$-module map from 
$A^{\infty}(W_{1})\otimes_{A^{\infty}(V)} W_{2}$ to $W_{3}$, that is, 
$$\rho(\Y)\in \hom_{A^{\infty}(V)}(A^{\infty}(W_{1})\otimes_{A^{\infty}(V)} 
W_{2}, W_{3}).$$
\end{prop}
\pf
By Proposition \ref{l-r-action-comm}, for $\mathfrak{v}\in U^{\infty}(V)$, $\mathfrak{w}_{1}\in U^{\infty}(W_{1})$
and $w_{2}\in W_{2}$, we have 
\begin{align*}
&(\rho(\Y))(((\mathfrak{v}+Q^{\infty}(V))\diamond (\mathfrak{w}_{1}+Q^{\infty}(W_{1})))\otimes 
w_{2})\nn
&\quad =(\rho(\Y))((\mathfrak{v}\diamond \mathfrak{w}_{1}+Q^{\infty}(W_{1}))\otimes 
w_{2})\nn
&\quad=\vartheta_{\Y}(\mathfrak{v}\diamond \mathfrak{w}_{1})w_{2}\nn
&\quad =\vartheta_{W_{3}}(\mathfrak{v})\vartheta_{\Y}(\mathfrak{w}_{1})w_{2}\nn
&\quad =\vartheta_{W_{3}}(\mathfrak{v})(\rho(\Y))
(([w_{1}]_{nl}+Q^{\infty}(W_{1}))\otimes 
w_{2})
\end{align*}
and 
\begin{align*}
&(\rho(\Y))(((\mathfrak{w}_{1}+Q^{\infty}(W_{1}))\diamond (\mathfrak{v}+Q^{\infty}(V)))\otimes 
w_{2})\nn
&\quad =(\rho(\Y))((\mathfrak{w}_{1}\diamond \mathfrak{v}+Q^{\infty}(V))\otimes 
w_{2})\nn
&\quad=\vartheta_{\Y}(\mathfrak{w}_{1}\diamond \mathfrak{v})w_{2}\nn
&\quad =\vartheta_{\Y}(\mathfrak{w}_{1})\vartheta_{W_{2}}(\mathfrak{v})w_{2}\nn
&\quad=(\rho(\Y))((\mathfrak{w}_{1}+Q^{\infty}(W_{1}))\otimes 
\vartheta_{W_{2}}(\mathfrak{v})w_{2}),
\end{align*}
proving that $\rho(\Y)$ is indeed an $A^{\infty}(V)$-module map from 
$A^{\infty}(W_{1})\otimes_{A^{\infty}(V)} W_{2}$ to $W_{3}$.
\epfv

We now have a linear map
\begin{align*}
\rho: \mathcal{V}_{W_{1}W_{2}}^{W_{3}}&\to 
\hom_{A^{\infty}(V)}(A^{\infty}(W_{1})\otimes_{A^{\infty}(V)} W_{2}, W_{3})\nn
\Y&\mapsto \rho(\Y).
\end{align*}

In the proof of Theorem \ref{main} below (in fact only in the surjectivity part), we need a conformal vector of $V$
(that is, $V$ is a vertex operator algebra) so that we do not have to verify separately 
the $L(0)$- and $L(-1)$-commutator formulas 
when we construct an intertwining operator. 

\begin{thm}\label{main}
Let $V$ be a vertex operator algebra. Then the linear map $\rho$ is an isomorphism.
\end{thm}
\pf
Let $\Y\in  \mathcal{V}_{W_{1}W_{2}}^{W_{3}}$. 
Assume that $\rho(\Y)=0$. Then for $k, l\in \N$, $\mu\in \Gamma(W_{2})$, 
$w_{1}\in W_{1}$, $w_{2}\in (W_{2})_{[h^{\mu}_{2}+l]}$, we have 
\begin{align*}
&\sum_{\nu\in \Gamma(W_{3})}
\res_{x}x^{h_{2}^{\mu}-h_{3}^{\nu}+l-k-1} \Y^{0}(x^{L_{W_{1}}(0)_{S}}w_{1}, x)w_{2}\nn
&\quad=\vartheta_{\Y}([w]_{kl})w_{2}\nn
&\quad=(\rho(\Y))(([w]_{kl}+Q^{\infty})\otimes
w_{2}) \nn
&\quad=0.
\end{align*}
So for $w_{3}'\in (W_{3}')_{[h_{3}^{\nu}+k]}$, 
we have 
\begin{equation}\label{main-1}
\langle w_{3}', \res_{x}x^{h_{2}^{\mu}-h_{3}^{\nu}+l-k-1} \Y^{0}(x^{L_{W_{1}}(0)_{S}}w_{1}, x)
w_{2}\rangle=0.
\end{equation}
On the other hand, 
\begin{align}\label{rho-o-n-3.5}
\Y^{0}(w, x)&=\text{Coeff}_{\log x}^{0}\Y(w, x)\nn
& =\text{Coeff}_{\log x}^{0}x^{L_{W_{3}}(0)}
\Y(x^{-L_{W}(0)}w, 1)x^{-L_{W_{2}}(0)}\nn
& =x^{L_{W_{3}}(0)_{S}}
\Y^{0}(x^{-L_{W}(0)_{S}}w, 1)x^{-L_{W_{2}}(0)_{S}}.
\end{align}
By (\ref{rho-o-n-3.5}),  
\begin{align}\label{main-2}
\langle w_{3}', & \res_{x}x^{h_{2}^{\mu}-h_{3}^{\nu}+l-k-1} \Y^{0}(x^{L_{W_{1}}(0)_{S}}w_{1}, x)
w_{2}\rangle\nn
&=\res_{x}x^{-1}\langle x^{-L_{W_{3}'}(0)_{S}}w_{3}', \Y^{0}(x^{L_{W_{1}}(0)_{S}}w_{1}, x)
x^{L_{W_{2}}(0)_{S}}w_{2}\rangle\nn
&=\res_{x}x^{-1}\langle w_{3}', \Y^{0}(w_{1}, 1)w_{2}\rangle\nn
&=\langle w_{3}', \Y^{0}(w_{1}, 1)w_{2}\rangle.
\end{align}
From (\ref{main-1}) and (\ref{main-2}), we obtain 
\begin{equation}\label{main-3}
\langle w_{3}', \Y^{0}(w_{1}, 1)w_{2}\rangle=0
\end{equation}
for $k, l\in \N$, $\mu\in \Gamma(W_{2})$, 
$\nu\in \Gamma(W_{3})$, $w_{1}\in W_{1}$, $w_{2}\in (W_{2})_{[h_{2}^{\mu}+l]}$
and $w_{3}'\in (W_{3}')_{[h_{3}^{\nu}+k]}$. Since $\mu, \nu, k, l$ are also arbitrary, 
we see that  (\ref{main-3}) holds for $w_{1}\in W_{1}$,
$w_{2}\in W_{2}$ and $w_{3}'\in W_{3}'$. Hence
$\Y^{0}(w_{1}, 1)=0$ for $w_{1}\in W_{1}$. Then by the $L(0)$-conjugation property of intertwining operators,
we have 
\begin{align*}
\Y(w_{1}, x)&=x^{L_{W_{3}}(0)}\Y(x^{-L_{W_{1}}(0)}w_{1}, 1)x^{-L_{W_{2}}(0)}\nn
&=x^{L_{W_{3}}(0)}\Y^{0}(x^{-L_{W_{1}}(0)}w_{1}, 1)x^{-L_{W_{2}}(0)}\nn
&=0
\end{align*}
for $w_{1}\in W_{1}$. Since $w_{1}$ is arbitrary, we obtain 
$\Y=0$, proving the injectivity of $\rho$.

We now prove the surjectivity of $\rho$. Let 
$$f\in \hom_{A^{\infty}(V)}(A^{\infty}(W_{1})\otimes_{A^{\infty}(V)} W_{2}, W_{3}).$$
We want to construct an intertwining operator $\Y^{f}$ of type $\binom{W_{3}}{W_{1}W_{2}}$ such that 
$\rho(Y^{f})=f$. For simplicity, we construct $\Y^{f}$  and prove the surjectivity only 
in the case $W_{2}=W_{2}^{\mu}\ne 0$ and $W_{3}=W_{3}^{\nu}\ne 0$ for $\mu, \nu\in \C/\Z$.
The surjectivity in the general case follows immediately. 

We define 
$$\Y_{h_{2}^{\mu}-h_{3}^{\nu}+l-k+\swt w_{1}-1, 0}^{f}(w_{1})w_{2}
=f(([w_{1}]_{kl}+Q^{\infty}(W_{1}))
\otimes_{A^{\infty}(V)} w_{2})$$
for $k, l\in \N$,  homogeneous $w_{1}\in W_{1}$ and $w_{2}\in 
(W_{2})_{[h_{2}^{\mu}+l]}$. Then by Lemma \ref{lemma-4.2},
$$\Y_{h_{2}^{\mu}-h_{3}^{\nu}+l-k+\swt w_{1}-1, 0}^{f}(w_{1})w_{2}\in 
(W_{3})_{\l k\r}=(W_{3})_{[h^{\nu}_{3}+k]}.$$
We define a map 
\begin{align*}
(\Y^{f})^{0}: W_{1}\otimes W_{2}&\to W_{3}\{x\}\nn
w_{1}\otimes w_{2}&\mapsto (\Y^{f})^{0}(w_{1}, x)w_{2}
\end{align*}
by
$$(\Y^{f})^{0}(w_{1}, x)w_{2}=\sum_{k\in \N}
\Y_{h_{2}^{\mu}-h_{3}^{\nu}+l-k+\swt w_{1}-1, 0}^{f}(w_{1})x^{-h_{2}^{\mu}+h_{3}^{\nu}-l+k-\swt w_{1}}$$
for $w_{1}\in W_{1}$, $w_{2}\in (W_{2})_{[h_{2}^{\mu}+l]}$ and $l\in \N$. 
We then define 
\begin{align}\label{defn-Y-f}
\Y^{f}(w_{1}, x)w_{2}&=x^{L_{W_{3}}(0)}(\Y^{f})^{0}(x^{-L_{W_{1}}(0)}w_{1}, 1)x^{-L_{W_{2}}(0)}w_{2}\nn
&=x^{L_{W_{3}}(0)_{N}}(\Y^{f})^{0}(x^{-L_{W_{1}}(0)_{N}}w_{1}, x)x^{-L_{W_{2}}(0)_{N}}w_{2}
\end{align}
for $w_{1}\in W_{1}$, $w_{2}\in W_{2}$ and $w_{3}\in W_{3}$. 
Since $L_{W_{1}}(0)_{N}$ and $L_{W_{2}}(0)_{N}$ are locally nilpotent operators, 
$x^{-L_{W_{1}}(0)_{N}}w_{1}\in W_{1}[\log x]$ and $x^{-L_{W_{2}}(0)_{N}}w_{2}\in W_{2}[\log x]$. 
Then the coefficient of $x$ to a power in 
$$(\Y^{f})^{0}(x^{-L_{W_{1}}(0)_{N}}w_{1}, x)x^{-L_{W_{2}}(0)_{N}}w_{2}$$ 
is an element of $W_{3}[\log x]$. Since $L_{W_{3}}(0)_{N}$ is also locally nilpotent,
the coefficient of $x$ to a power 
in the right-hand side of (\ref{defn-Y-f}) is also in $W_{3}[\log x]$. This shows that 
$\Y^{f}(w_{1}, x)w_{2}$ is in fact in $W_{3}[\log x]\{x\}$. In particular, we obtain a linear map
\begin{align*}
\Y^{f}: W_{1}\otimes W_{2}&\to W_{3}[\log x]\{x\}\nn
w_{1}\otimes w_{2}&\mapsto \Y^{f}(w_{1}, x)w_{2}.
\end{align*}

We now prove that $\Y^{f}$ is indeed an intertwining operator of type $\binom{W_{3}}{W_{1}W_{2}}$. 
From the definition of $\Y^{f}$, we see that it is lower truncated. 
Also by definition, for $w_{1}\in W_{1}$, 
$w_{2}\in (W_{2})_{[h_{2}^{\mu}+l]}$, $l\in \N$,
\begin{align*}
(\Y^{f})^{0}(w_{1}, 1)w_{2}&=\sum_{k\in \N}
\Y_{h_{2}^{\mu}-h_{3}^{\nu}+l-k+\swt w_{1}-1, 0}^{f}(w_{1})w_{2}\nn
&=\sum_{k\in \N}
f(([w_{1}]_{kl}+Q^{\infty}(W_{1}))
\otimes_{A^{\infty}(V)} w_{2}).
\end{align*}
Then we have 
$$\Y^{f}(w_{1}, x)w_{2}=
\sum_{k\in \N}x^{L_{W_{3}}(0)}
f(([x^{-L_{W_{1}}(0)}w_{1}]_{kl}+Q^{\infty}(W_{1}))
\otimes_{A^{\infty}(V)} x^{-L_{W_{1}}(0)}w_{2})$$
and 
\begin{align}\label{main-3.5}
&f(([w_{1}]_{kl}+Q^{\infty}(W_{1}))
\otimes_{A^{\infty}(V)} w_{2})\nn
&\quad=\res_{x}x^{h_{2}^{\mu}-h_{3}^{\nu}+l-k-1}(\Y^{f})^{0}(x^{L_{W_{1}}(0)_{S}}w_{1}, x)w_{2}\nn
&\quad=\text{Coeff}_{\log x}^{0}\res_{x}x^{h_{2}^{\mu}-h_{3}^{\nu}+l-k-1}\Y^{f}(x^{L_{W_{1}}(0)}w_{1}, x)w_{2}.
\end{align}

From Proposition \ref{jacobi-ker-theta}, we know that (\ref{jacobi-ker-theta-0}) with  $w$ replaced by $w_{1}$ 
is an element of $Q^{\infty}(W_{1})$. 
In particular, 
\begin{align}\label{main-4}
&\sum_{j\in \N}(-1)^{j}\binom{p}{j}f(([v]_{k,n+p-j}\diamond
[w_{1}]_{n+p-j,l+p}+Q^{\infty}(W_{1}))
\otimes_{A^{\infty}(W_{1})} w_{2})\nn
&\quad \quad -\sum_{\stackrel{i\in \N}{l-n+k+p-j\ge 0}}(-1)^{p-j}\binom{p}{j}\cdot \nn
&\quad \quad\quad \quad\cdot
f(([w_{1}]_{k,l-n+k+p-j}\diamond [v]_{l-n+k+p-j,l+p}
+Q^{\infty}(W_{1}))
\otimes_{A^{\infty}(W_{1})}w_{2})\nn
&\quad \quad-\sum_{j\in \N}\binom{\wt v+n-k-1}{j}
f(([(Y_{W_{1}})_{p+j}(v)w_{1}]_{k,l+p}+Q^{\infty}(W_{1}))
\otimes_{A^{\infty}(W_{1})}w_{2})\nn
&\quad =0
\end{align}
for $k, l, n\in \N$, $p\in \Z$ such that $l+p\in \N$,
$v\in V$, $w_{1}\in W_{1}$ and $w_{2}\in (W_{2})_{[h_{2}^{\mu}+l+p]}$. 
Since 
$$f\in \hom_{A^{\infty}(V)}(A^{\infty}(W_{1})\otimes_{A^{\infty}(V)} W_{2}, W_{3}),$$
the left-hand side of (\ref{main-4}) is equal to 
\begin{align}\label{main-5}
&\sum_{j\in \N}(-1)^{j}\binom{p}{j}\vartheta_{W_{3}}([v]_{k,n+p-j})
f(([w_{1}]_{n+p-j,l+p}+Q^{\infty}(W_{1}))
\otimes_{A^{\infty}(W_{1})}w_{2})\nn
&\quad -\sum_{\stackrel{j\in \N}{l-n+k+p-j\ge 0}}(-1)^{p-j}\binom{p}{j}\cdot \nn
&\quad \quad\quad \quad\cdot
f(([w_{1}]_{k,l-n+k+p-j}+Q^{\infty}(W_{1}))
\otimes_{A^{\infty}(W_{1})}\vartheta_{W_{2}}([v]_{l-n+k+p-j,l+p})w_{2})\nn
&\quad-\sum_{j\in \N}\binom{\wt v+n-k-1}{j}
f(([(Y_{W_{1}})_{p+j}(v)w_{1}]_{k,l+p}+Q^{\infty}(W_{1}))
\otimes_{A^{\infty}(W_{1})}w_{2}).
\end{align}
From the definitions of $\vartheta_{W_{3}}$, $\vartheta_{W_{2}}$ and (\ref{main-3.5}), 
we see that (\ref{main-5}) is equal to 
\begin{align}\label{main-6}
&\sum_{j\in \N}(-1)^{j}\binom{p}{j}\res_{x_{1}}\res_{x_{2}}x_{1}^{n+p-j-k-1}\cdot\nn
&\quad\quad\quad\quad\quad\quad\quad\quad\quad\cdot
x_{2}^{h_{2}^{\mu}-h_{3}^{\nu}+l-n+j-1}Y_{W_{3}}(x_{1}^{L_{V}(0)}v, x_{1})
(\Y^{f})^{0}(x_{2}^{L_{W_{1}}(0)_{S}}w_{1}, x_{2})w_{2}\nn
&\quad -\sum_{\stackrel{j\in \N}{l-n+k+p-j\ge 0}}(-1)^{p-j}\binom{p}{j}
\res_{x_{1}}\res_{x_{2}}x_{1}^{n-k+j-1}\cdot\nn
&\quad\quad\quad\quad\quad\quad\quad\quad\quad\cdot x_{2}^{h_{2}^{\mu}-h_{3}^{\nu}+l-n+p-j-1}
(\Y^{f})^{0}(x_{2}^{L_{W_{1}}(0)_{S}}w_{1}, x_{2})Y_{W_{2}}(x_{1}^{L_{V}(0)}v, x_{1})w_{2}\nn
&\quad-\sum_{j\in \N}\binom{\wt v+n-k-1}{j}
\res_{x_{2}}x_{2}^{h_{2}^{\mu}-h_{3}^{\nu}+l+p-k-1}(\Y^{f})^{0}(x_{2}^{L_{W_{1}}(0)_{S}}
(Y_{W_{1}})_{p+j}(v)w_{1}, x_{2})w_{2}\nn
&=\res_{x_{1}}\res_{x_{2}}x_{1}^{n-k-1}
x_{2}^{h_{2}^{\mu}-h_{3}^{\nu}+l-n-1}(x_{1}-x_{2})^{p}Y_{W_{3}}(x_{1}^{L_{V}(0)}v, x_{1})
(\Y^{f})^{0}(x_{2}^{L_{W_{1}}(0)_{S}}w_{1}, x_{2})w_{2}\nn
&\quad -
\res_{x_{1}}\res_{x_{2}}x_{1}^{n-k-1}x_{2}^{h_{2}^{\mu}-h_{3}^{\nu}+l-n-1}
(-x_{2}+x_{1})^{p}
(\Y^{f})^{0}(x_{2}^{L_{W_{1}}(0)_{S}}w_{1}, x_{2})Y_{W_{2}}(x_{1}^{L_{V}(0)}v, x_{1})w_{2}\nn
&\quad-
\res_{x_{2}}\res_{y}y^{p}(1+y)^{\swt v+n-k-1}
x_{2}^{h_{2}^{\mu}-h_{3}^{\nu}+l+p-k-1}(\Y^{f})^{0}(x_{2}^{L_{W_{1}}(0)_{S}}
Y_{W_{1}}(v,  y)w_{1}, x_{2})w_{2}.
\end{align}
Since the left-hand side of (\ref{main-6}) is equal to the left-hand side of  (\ref{main-4}), we see that 
the right-hand side of (\ref{main-6}) is $0$ for $k, l, n\in \N$, $p\in \Z$ such that $l+p\in \N$, 
$v\in V$, $w_{1}\in W_{1}$ 
and $w_{2}\in (W_{2})_{[h_{2}^{\mu}+l+p]}$.  The right-hand side of (\ref{main-6}) is also 
$0$ when $l+p<0$ since in this case $w_{2}=0$. 
Multiplying the right-hand side 
of  (\ref{main-6}) by $x_{0}^{-p-1}$ and then take sum over $p\in \Z$, we obtain
\begin{align}\label{main-7}
&\res_{x_{1}}\res_{x_{2}}x_{1}^{n-k-1}
x_{2}^{h_{2}^{\mu}-h_{3}^{\nu}+l-n-1}x_{0}^{-1}\delta\left(\frac{x_{1}-x_{2}}{x_{0}}\right)
Y_{W_{3}}(x_{1}^{L_{V}(0)}v, x_{1})
(\Y^{f})^{0}(x_{2}^{L_{W_{1}}(0)_{S}}w_{1}, x_{2})w_{2}\nn
&\quad -
\res_{x_{1}}\res_{x_{2}}x_{1}^{n-k-1}x_{2}^{h_{2}^{\mu}-h_{3}^{\nu}+l-n-1}
\cdot\nn
&\quad\quad\quad\quad\quad\quad\quad\quad\quad\cdot x_{0}^{-1}\delta\left(\frac{x_{2}-x_{1}}{-x_{0}}\right)
(\Y^{f})^{0}(x_{2}^{L_{W_{1}}(0)_{S}}w_{1}, x_{2})Y_{W_{2}}(x_{1}^{L_{V}(0)}v, x_{1})w_{2}\nn
&\quad-
\res_{x_{2}}\res_{y}x_{0}^{-1}\delta\left(\frac{yx_{2}}{x_{0}}\right)(1+y)^{n-k-1}
\cdot\nn
&\quad\quad\quad\quad\quad\quad\quad\quad\quad\cdot x_{2}^{h_{2}^{\mu}-h_{3}^{\nu}+l-k-1}(\Y^{f})^{0}(x_{2}^{L_{W_{1}}(0)_{S}}
Y_{W}((1+y)^{L_{V}(0)}v, y)w_{1}, x_{2})w_{2}\nn
&=0.
\end{align}
The third term in the left-hand side of (\ref{main-7}) is equal to 
\begin{align}\label{main-8}
&-\res_{x_{2}}\res_{y}x_{0}^{-1}\delta\left(\frac{yx_{2}}{x_{0}}\right)(1+y)^{n-k-1}
x_{2}^{h_{2}^{\mu}-h_{3}^{\nu}+l-k-1}\cdot\nn
&\quad\quad\quad\quad\quad\quad\quad\cdot  (\Y^{f})^{0}(x_{2}^{L_{W_{1}}(0)_{S}}
Y_{W}((1+y)^{L_{V}(0)}v, y)w_{1}, x_{2})w_{2}\nn
&\quad =-\res_{x_{2}}\res_{y}x_{0}^{-1}\delta\left(\frac{yx_{2}}{x_{0}}\right)(1+x_{0}x_{2}^{-1})^{n-k-1}
x_{2}^{h_{2}^{\mu}-h_{3}^{\nu}+l-k-1}\cdot\nn
&\quad\quad\quad\quad\quad\quad\quad\cdot  
(\Y^{f})^{0}(x_{2}^{L_{W_{1}}(0)_{S}}
Y_{W}((1+x_{0}x_{2}^{-1})^{L_{V}(0)}v, x_{0}x_{2}^{-1})w_{1}, x_{2})w_{2}\nn
&\quad =-\res_{x_{2}}(x_{2}+x_{0})^{n-k-1}
x_{2}^{h_{2}^{\mu}-h_{3}^{\nu}+l-n-1}\cdot\nn
&\quad\quad\quad\quad\quad\quad\quad\cdot  
(\Y^{f})^{0}(
Y_{W}((x_{2}+x_{0})^{L_{V}(0)}v, x_{0})x_{2}^{L_{W_{1}}(0)_{S}}w_{1}, x_{2})w_{2}\nn
&\quad =-\res_{x_{1}}\res_{x_{2}}x_{1}^{n-k-1}
x_{2}^{h_{2}^{\mu}-h_{3}^{\nu}+l-n-1}x_{1}^{-1}\delta\left(\frac{x_{2}+x_{0}}{x_{1}}\right)\cdot\nn
&\quad\quad\quad\quad\quad\quad\quad\cdot  
(\Y^{f})^{0}(
Y_{W}(x_{1}^{L_{V}(0)}v, x_{0})x_{2}^{L_{W_{1}}(0)_{S}}w_{1}, x_{2})w_{2}.
\end{align}
Using (\ref{main-8}) and substituting $x_{1}^{-L_{V}(0)}v$ and $x_{2}^{-L_{W_{1}}(0)_{S}}w_{1}$
for $v$ and $w_{1}$, respectively, we see that (\ref{main-7}) becomes
\begin{align}\label{main-9}
&\res_{x_{1}}\res_{x_{2}}x_{1}^{n-k-1}
x_{2}^{h_{2}^{\mu}-h_{3}^{\nu}+l-n-1}x_{0}^{-1}\delta\left(\frac{x_{1}-x_{2}}{x_{0}}\right)
Y_{W_{3}}(v, x_{1})
(\Y^{f})^{0}(w_{1}, x_{2})w_{2}\nn
&\quad -
\res_{x_{1}}\res_{x_{2}}x_{1}^{n-k-1}x_{2}^{h_{2}^{\mu}-h_{3}^{\nu}+l-n-1}
x_{0}^{-1}\delta\left(\frac{x_{2}-x_{1}}{-x_{0}}\right)
(\Y^{f})^{0}(w_{1}, x_{2})Y_{W_{2}}(v, x_{1})w_{2}\nn
&\quad-\res_{x_{1}}\res_{x_{2}}x_{1}^{n-k-1}
x_{2}^{h_{2}^{\mu}-h_{3}^{\nu}+l-n-1}x_{1}^{-1}\delta\left(\frac{x_{2}+x_{0}}{x_{1}}\right)
(\Y^{f})^{0}(
Y_{W}(v, x_{0})w_{1}, x_{2})w_{2}\nn
&=0.
\end{align}
Since $k, l,n\in \N$ are arbitrary, (\ref{main-9}) gives 
\begin{align}\label{main-10}
&x_{0}^{-1}\delta\left(\frac{x_{1}-x_{2}}{x_{0}}\right)
Y_{W_{3}}(v, x_{1})
(\Y^{f})^{0}(w_{1}, x_{2})w_{2}-
x_{0}^{-1}\delta\left(\frac{x_{2}-x_{1}}{-x_{0}}\right)
(\Y^{f})^{0}(w_{1}, x_{2})Y_{W_{2}}(v, x_{1})w_{2}\nn
&\quad=x_{1}^{-1}\delta\left(\frac{x_{2}+x_{0}}{x_{1}}\right)
(\Y^{f})^{0}(
Y_{W}(v, x_{0})w_{1}, x_{2})w_{2}.
\end{align}
This is the Jacobi identity for $(\Y^{f})^{0}$. 
To obtain the Jacobi identity for $\Y^{f}$, we replace $x_{0}$, $x_{1}$, $x_{2}$, $v$, $w_{1}$ and $w_{2}$
by $x_{0}x_{2}^{-1}$, 
$x_{1}x_{2}^{-1}$, $1$, $x_{2}^{-L_{V}(0)}v$,
$x_{2}^{-L_{W_{1}}(0)}w_{1}$ and $x^{-L_{W_{2}}(0)}w_{2}$, respectively
 in (\ref{main-10}) and then multiply $x^{L_{W_{3}}(0)}$ from the left of both sides of (\ref{main-10}).
Then we obtain 
\begin{align}\label{main-11}
&x_{0}^{-1}x_{2}\delta\left(\frac{x_{1}x_{2}^{-1}-1}{x_{0}x_{2}^{-1}}\right)
x_{2}^{L_{W_{1}}(0)}Y_{W_{3}}(x_{2}^{-L_{V}(0)}v, x_{1}x_{2}^{-1})
(\Y^{f})^{0}(x_{2}^{-L_{W_{1}}(0)}w_{1}, 1)x_{2}^{-L_{W_{1}}(0)}w_{2}\nn
&\quad \quad -
x_{0}^{-1}x_{2}\delta\left(\frac{1-x_{1}x_{2}^{-1}}{-x_{0}x_{2}^{-1}}\right)x_{2}^{L_{W_{1}}(0)}
(\Y^{f})^{0}(x_{2}^{-L_{W_{1}}(0)}w_{1}, 1)Y_{W_{2}}(x_{2}^{-L_{V}(0)}v, x_{1}
x_{2}^{-1})x_{2}^{-L_{W_{1}}(0)}w_{2}\nn
&\quad=x_{1}^{-1}x_{2}\delta\left(\frac{1+x_{0}x_{2}^{-1}}{x_{1}x_{2}^{-1}}\right)
x_{2}^{L_{W_{1}}(0)}(\Y^{f})^{0}(
Y_{W}(x_{2}^{-L_{V}(0)}v, x_{0}x_{2}^{-1})x_{2}^{-L_{W_{1}}(0)}w_{1}, 1)w_{2}.
\end{align}
Using the $L(0)$-conjugation formulas for vertex operators and the definition 
$$\Y^{f}(w_{1}, x_{2})=x_{2}^{L_{W_{1}}(0)}(\Y^{f})^{0}(x_{2}^{-L_{W_{1}}(0)}w_{1}, 1)
x_{2}^{-L_{W_{1}}(0)}$$
of 
$\Y^{f}$, we see that (\ref{main-11}) is equivalent to
\begin{align*}
&x_{0}^{-1}\delta\left(\frac{x_{1}-x_{2}}{x_{0}}\right)
Y_{W_{3}}(v, x_{1})
\Y^{f}(w_{1}, x_{2})w_{2}
-x_{0}^{-1}\delta\left(\frac{x_{2}-x_{1}}{-x_{0}}\right)
\Y^{f}(w_{1}, x_{2})Y_{W_{2}}(v, x_{1})w_{2}\nn
&\quad=x_{1}^{-1}\delta\left(\frac{x_{2}+x_{0}}{x_{1}}\right)
\Y^{f}(Y_{W}(v, x_{0})w_{1}, x_{2})w_{2},
\end{align*}
the Jacobi identity for $\Y$. 

We also need to show that $\Y^{f}$ satisfies the $L(-1)$-derivative property. 
This follows from (\ref{defn-Y-f}) and the $L(0)$-commutator formula for
$(\Y^{f})^{0}$ which is a special case of the Jacobi identity for $(\Y^{f})^{0}$. 
In fact, applying $\frac{d}{dx}$ to both sides of (\ref{defn-Y-f}), we obtain
\begin{align}\label{der-Y-f}
\frac{d}{dx}\Y^{f}(w_{1}, x)w_{2}
&=x^{-1}x^{L_{W_{3}}(0)}L_{W_{3}}(0)
(\Y^{f})^{0}(x^{-L_{W_{1}}(0)}w_{1}, 1)x^{-L_{W_{2}}(0)}w_{2}\nn
&\quad
-x^{-1}x^{L_{W_{3}}(0)}(\Y^{f})^{0}(L_{W_{1}}(0)x^{-L_{W_{1}}(0)}w_{1}, 1)x^{-L_{W_{2}}(0)}w_{2}\nn
&\quad
-x^{-1}x^{L_{W_{3}}(0)}(\Y^{f})^{0}(x^{-L_{W_{1}}(0)}w_{1}, 1)L_{W_{2}}(0)x^{-L_{W_{2}}(0)}w_{2}
\end{align}
for $w_{1}\in W_{1}$ and $w_{2}\in W_{2}$. 
Using the $L(0)$-commutator formula 
\begin{align*}
L_{W_{3}}(0)&(\Y^{f})^{0}(x^{-L_{W_{1}}(0)}w_{1}, 1)
-(\Y^{f})^{0}(x^{-L_{W_{1}}(0)}w_{1}, 1)
L_{W_{2}}(0)\nn
&=(\Y^{f})^{0}((L_{W_{1}}(-1)+L_{W_{1}}(0))x^{-L_{W_{1}}(0)}w_{1}, 1)
\end{align*}
and (\ref{defn-Y-f}) with $w_{1}$ replaced by $L_{W_{1}}(-1)w_{1}$, 
the right-hand side of (\ref{der-Y-f}) is equal to 
\begin{align}\label{der-Y-f-1}
&x^{-1}x^{L_{W_{3}}(0)}
(\Y^{f})^{0}(L_{W_{1}}(-1)x^{-L_{W_{1}}(0)}w_{1}, 1)x^{-L_{W_{2}}(0)}w_{2}\nn
&\quad =x^{L_{W_{3}}(0)}
(\Y^{f})^{0}(x^{-L_{W_{1}}(0)}L_{W_{1}}(-1)w_{1}, 1)x^{-L_{W_{2}}(0)}w_{2}\nn
&\quad =\Y^{f}(L_{W_{1}}(-1)w_{1}, x)w_{2}.
\end{align}
From (\ref{der-Y-f}) and (\ref{der-Y-f-1}) , we obtain the $L(-1)$-derivative property for $\Y^{f}$. 

Since $V$ is a vertex operator algebra and $\Y^{f}$ satisfies the lower-truncation property, the Jacobi identity
and the $L(-1)$-derivative property, $\Y^{f}$ is an intertwining operator of type $\binom{W_{3}}{W_{1}W_{2}}$
(see Definition 3.10 in \cite{HLZ} for the definition of intertwining operator  in the case that $V$ is a vertex operator algebra).

Finally we prove $\rho(Y^{f})=f$. By definition, 
\begin{align}\label{main-20}
(\rho(\Y^{f}))&(([w_{1}]_{kl}+Q^{\infty}(W_{1}))
\otimes_{A^{\infty}(V)}w_{2})\nn
&=\vartheta_{\Y^{f}}([w_{1}]_{kl})w_{2}\nn
&=\res_{x}x^{h_{2}^{\mu}-h_{3}^{\nu}+l-k-1}(\Y^{f})^{0}(x^{L_{W_{1}}(0)}w_{1}, x)w_{2}\nn
&=f(([w_{1}]_{kl}+Q^{\infty}(W_{1}))
\otimes_{A^{\infty}(V)}w_{2})
\end{align}
for $k, l\in \N$,
$w_{1}\in W_{1}$ and $w_{2}\in (W_{2})_{[h_{2}^{\mu}+l]}$. 
 From Lemma \ref{lemma-4.1} and (\ref{main-20}), we obtain
$\rho(Y^{f})=f$. This finishes the proof that $\rho$ is surjective.
\epfv

\renewcommand{\theequation}{\thesection.\arabic{equation}}
\setcounter{equation}{0}
\setcounter{thm}{0}
\section{$A^{N}(V)$-bimodules and intertwining operators}

In this section, we formulate and prove the second main theorem of the present paper. 
For $N\in \N$, we use the results obtained in the preceding two sections to give an
$A^{N}(V)$-bimodule $A^{N}(W)$ for a lower-bounded generalized $V$-module $W$. For 
lower bounded generalized $V$-modules $W_{1}$, $W_{2}$ and $W_{3}$,  
we obtain a linear map $\rho^{N}: \mathcal{V}_{W_{1}W_{2}}^{W_{3}}\to 
\hom_{A^{N}(V)}(A^{N}(W_{1})\otimes_{A^{N}(V)}\Omega_{N}^{0}(W_{2}), 
\Omega_{N}^{0}(W_{3})$ induced from the linear isomorphism $\rho$ in the preceding section (see below for the definition of $\Omega_{n}^{0}(W)$
for a lower-bounded generalized $V$-module).
We prove that $\rho^{N}$ is surjective. 
Our second main theorem states that $\rho^{N}$ is an isomorphism when $W_{2}$ and $W_{3}'$ are
certain universal lower-bounded generalized $V$-modules generated by 
$\Omega_{N}^{0}(W_{2})$ and $\Omega_{N}^{0}(W_{3}')$, respectively.

Let $V$ be a grading-restricted vertex algebra and 
$W$ a lower-bounded generalized $V$-module. For $n\in \N$, let
$$\Omega_{n}^{0}(W)=\coprod_{m=0}^{n}W_{\l m\r}=\coprod_{m=0}^{n}\coprod_{\mu\in \Gamma(W)}
W_{[h^{\mu}+m]}.$$ 
For homogeneous $v\in V$, $ k, l, n\in \N$, $\mu\in \Gamma(W)$ 
and $w\in W_{[h^{\mu}+l]}$, we have 
$$\vartheta_{W}([v]_{kn})w=\delta_{nl}
\res_{x}x^{l-k-1}Y_{W}(x^{L_{V}(0)}v, x)w\in W_{[h^{\mu}+k]} \subset W_{\l k\r}.$$

We fix $N\in \N$ in the rest of this section. 
Recall the associative algebra $A^{N}(V)$ (see Section 2
and Subsection 4.2 in \cite{H-aa-va}). 
Then the discussion above in particular shows that when restricted to $A^{N}(V)$, $\vartheta_{W}$
gives 
$\Omega_{N}^{0}(W)$ a graded 
$A^{N}(V)$-module structure. But in general, $\Omega_{N}^{0}(W)$ is not nondegenerate 
as a graded 
$A^{N}(V)$-module (see Definition \ref{gr-A-N-mod} for the definition of 
 nondegenerate graded 
$A^{N}(V)$-module).

Let $U^{N}(W)$ be the space of all $(N+1)\times (N+1)$ matrices with entries in $W$.
Then $U^{N}(W)$ can be viewed as subspaces of $U^{\infty}(W)$. 
In this paper, we shall always view $U^{N}(W)$
as subspaces of  $U^{\infty}(W)$ and we shall use the 
notations for elements of $U^{\infty}(W)$ to denote elements
of $U^{N}(W)$. Using these notations, we see that 
$U^{N}(W)$ are spanned by elements of the forms 
$[w]_{kl}$ for $w\in W$ and $k, l=0, \dots, N$. 

Let 
$A^{N}(W)$ be the subspace of $A^{\infty}(W)$ consisting of elements of the form 
$\mathfrak{w}+Q^{\infty}(W)$ for $\mathfrak{w}\in U^{N}(W)$. 
Since $A^{N}(V)$ is a subalgebra of $A^{\infty}(V)$, 
$A^{\infty}(W)$ as an $A^{\infty}(V)$-bimodule is also an $A^{N}(V)$-bimodule. 
For $v\in V$, $w\in W$ and $k, m, n, l=0, \dots, N$, by
definition, $[v]_{km}\diamond [w]_{nl}$ and 
$[w]_{km}\diamond [v]_{nl}$
are still in $U^{N}(W)$. 
Thus $A^{N}(W)$ is in fact 
an $A^{N}(V)$-subbimodule of $A^{\infty}(W)$.

Let $W_{1}$, $W_{2}$ and $W_{3}$ be lower-bounded generalized $V$-modules and 
$\Y$ an intertwining operator of type $\binom{W_{3}}{W_{1}W_{2}}$. 
Let  
$$\eta^{N}: A^{N}(W_{1})\otimes_{A^{N}(V)} \Omega_{N}^{0}(W_{2})\to  
A^{\infty}(W_{1})
\otimes_{A^{\infty}(V)}W_{2}$$
be the $A^{N}(V)$-module map defined by 
$$\eta^{N}((\mathfrak{w}_{1}+Q^{\infty}(V))\otimes_{A^{N}(V)} w_{2})= 
(\mathfrak{w}_{1}+Q^{\infty}(V))\otimes_{A^{\infty}(V)} w_{2}$$
for $\mathfrak{w}_{1}\in U^{N}(W_{1})$  and $w_{2}\in \Omega_{N}^{0}(W_{2})$. 
Given
$$f\in \hom_{A^{\infty}(V)}(A^{\infty}(W_{1})\otimes_{A^{\infty}(V)}W_{2}, 
W_{3}),$$ 
we have a map
$f^{N}=f\circ \eta^{N}$. 
By Lemma \ref{lemma-4.2}, 
the image of 
$A^{N}(W_{1})\otimes_{A^{N}(V)}\Omega_{N}^{0}(W_{2})$ under $f^{N}$
is in fact in $\Omega_{N}^{0}(W_{3})$. 
So  $f^{N}$ is an element of 
$$\hom_{A^{N}(V)}(A^{N}(W_{1})\otimes_{A^{N}(V)} \Omega_{N}^{0}(W_{2}), 
\Omega_{N}^{0}(W_{3})).$$ 
Hence $f\mapsto f^{N}$ gives a linear map from 
$$\hom_{A^{\infty}(V)}(A^{\infty}(W_{1})\otimes_{A^{\infty}(V)}W_{2}, 
W_{3})$$ 
to 
$$\hom_{A^{N}(V)}(A^{N}(W_{1})\otimes_{A^{N}(V)} \Omega_{N}^{0}(W_{2}), 
\Omega_{N}^{0}(W_{3})).$$
In particular, the image of 
$$\rho(\Y)\in  \hom_{A^{\infty}(V)}
(A^{\infty}(W_{1})\otimes_{A^{\infty}(V)} W_{2}, W_{3})$$
under this map is an element 
$$\rho^{N}(\Y)=\rho(\Y)\circ \eta^{N}\in \hom_{A^{N}(V)}
(A^{N}(W_{1})\otimes_{A^{N}(V)} \Omega_{N}^{0}(W_{2}), \Omega_{N}^{0}(W_{3}))$$
More explicitly,  $(\rho^{N}(\Y))$ is given by
\begin{align*}
&(\rho^{N}(\Y))([w_{1}]_{kl}+Q^{\infty}(W_{1}))\otimes_{A^{N}(V)} w_{2})\nn
&\quad =(\rho(\Y))([w_{1}]_{kl}+Q^{\infty}(W_{1}))\otimes_{A^{\infty}(V)} w_{2})\nn
&\quad =\vartheta_{\Y}([w_{1}]_{kl})w_{2}\nn
&\quad =\sum_{\nu\in \Gamma(W_{3})}
\res_{x}x^{h_{2}^{\mu}-h_{3}^{\nu}+l-k-1}
\Y^{0}(x^{L_{W_{1}}(0)}w_{1}, x)w_{2}
\end{align*}
for $k, l=0, \dots, N$, $w_{1}\in W_{1}$, $w_{2}\in (W_{2})_{[h_{2}^{\mu}+l]}$ and $\mu\in \Gamma(W_{2})$.

We now have a linear map 
\begin{align*}
\rho^{N}: \mathcal{V}_{W_{1}W_{2}}^{W_{3}}&\to 
\hom_{A^{N}(V)}(A^{N}(W_{1})\otimes_{A^{N}(V)} \Omega_{N}^{0}(W_{2}), 
\Omega_{N}^{0}(W_{3}))\\
\Y&\mapsto \rho^{N}(\Y).
\end{align*}

\begin{prop}\label{rho-N-inj}
Assume that $W_{2}$ and $W_{3}'$ are generated by $\Omega_{N}^{0}(W_{2})$ and $\Omega_{N}^{0}(W_{3}')$.
Then the linear map $\rho^{N}$ is injective.
\end{prop}
\pf
Assume that $\rho^{N}(\Y)=0$. Then for $\mathfrak{w}_{1}\in U^{N}(W_{1})$, $w_{2}\in \Omega_{N}^{0}(W_{2})$, 
\begin{equation}\label{rho-N-inj-1}
(\rho(\Y))(\mathfrak{w}_{1}+Q^{\infty}(W_{1}))\otimes_{A^{\infty}(V)} w_{2})=0.
\end{equation}
Since $W_{2}$ is generated by $\Omega_{N}^{0}(W_{2})$, $W_{2}$ is spanned by elements of the form
$(Y_{W_{2}})_{\swt v+l-k-1}(v)w_{2}$ for $k\in \N$, $0\le l\le N$,  
homogeneous $v\in V$ and  $w_{2}\in (W_{2})_{\l l\r}$. For $k\in \N$, $0\le l\le N$,  $\mathfrak{w}_{1}\in U^{N}(W_{1})$
homogeneous $v\in V$ and  $w_{2}\in (W_{2})_{\l l\r}\subset \Omega_{N}^{0}(W_{2})$, we have 
\begin{align*}
&(\rho(\Y))(\mathfrak{w}_{1}+Q^{\infty}(W_{1}))\otimes_{A^{\infty}(V)} (Y_{W_{2}})_{\swt v+l-k-1}(v)w_{2})\nn
&\quad =(\rho(\Y))((\mathfrak{w}_{1}+Q^{\infty}(W_{1}: X))\otimes_{A^{\infty}(V)} \vartheta_{W_{2}}([v]_{kl})w_{2})\nn
&\quad =(\rho(\Y))((\mathfrak{w}_{1}\diamond [v]_{kl}+Q^{\infty}(W_{1}))\otimes_{A^{\infty}(V)} w_{2})\nn
&\quad =0.
\end{align*}
So (\ref{rho-N-inj-1}) holds for $\mathfrak{w}_{1}\in U^{N}(W_{1})$ and $w_{2}\in W_{2}$.

Since $W_{3}'$ is generated by $\Omega_{N}^{0}(W_{3}')$, 
every element of $W_{3}'$ is a linear combination of elements of the form
$(Y_{W_{3}})_{\swt v+k-n-1}(v)'w_{3}'$
for $n\in \N$, $0\le k\le N$,  homogeneous $v\in V$, $w_{3}'\in (W_{3}')_{\l n\r}$, 
where $(Y_{W_{3}})_{\swt v+k-n-1}(v)'$ is the adjoint of 
$(Y_{W_{3}})_{\swt v+k-n-1}(v)$. Then 
for $n\in \N$, $0\le k, l\le N$, $\mathfrak{w}_{1}\in U^{\infty}(W_{1})$
homogeneous $v\in V$, $w_{2}\in W_{2}$ and $w_{3}'\in (W_{3}')_{\l k\r}$, 
we have 
\begin{align*}
&\langle (Y_{W_{3}})_{\swt v+k-n-1}(v)'w_{3}', (\rho(\Y))(\mathfrak{w}_{1}+Q^{\infty}(W_{1}))
\otimes_{A^{\infty}(V)} w_{2})\rangle\nn
&\quad =\langle w_{3}', (Y_{W_{3}})_{\swt v+k-n-1}(v)(\rho(\Y))(\mathfrak{w}_{1}+Q^{\infty}(W_{1}))
\otimes_{A^{\infty}(V)} w_{2})\rangle\nn
&\quad =\langle w_{3}', \vartheta_{W_{3}}([v]_{nk})(\rho(\Y))(\mathfrak{w}_{1}+Q^{\infty}(W_{1}))
\otimes_{A^{\infty}(V)} w_{2})\rangle\nn
&\quad =\langle w_{3}', (\rho(\Y))([v]_{nk}\diamond \mathfrak{w}_{1}+Q^{\infty}(W_{1}))
\otimes_{A^{\infty}(V)} w_{2})\rangle\nn
&\quad =0.
\end{align*}
Then we see that (\ref{rho-N-inj-1}) holds for $\mathfrak{w}_{1}\in U^{\infty}(W_{1})$ and $w_{2}\in W_{2}$.
Thus we obtain $\rho(\Y)=0$. 
By Theorem \ref{main}, $\rho$ is injective. So $\Y=0$, proving the injectivity of $\rho^{N}$. 
\epfv

In general, $\rho^{N}$ is not surjective. But we shall prove that in the case that $W_{2}$ and $W_{3}'$ 
are equivalent to certain universal lower-bounded generalized $V$-modules, it is also surjective. 
Such a universal lower-bounded generalized $V$-module $S^{N}(M)$ has been constructed 
from an $A^{N}(V)$-module $M$ in 
Section 5 of \cite{H-aa-va} using the construction in Section 5 of \cite{H-const-twisted-mod}
for a grading-restricted vertex algebra. 

Since in our result below on the surjectivity, we assume for simplicity that 
$V$ is a vertex operator algebra, we shall instead use the modified construction 
for vertex operator algebra in subsection 4.2 of \cite{H-affine-twisted-mod}. 

In the remaining part of this section, $V$ is a vertex operator algebra.
Let $M=\coprod_{n=0}^{N}G_{n}(M)$ be a 
graded $A^{N}(V)$-module given by a linear map $\vartheta_{M}: A^{N}(V)\to {\rm End}\;M$
and operators 
$L_{M}(0)$ and $L_{M}(-1)$ (see Definition \ref{gr-A-N-mod}). Take $g=1_{V}$ and $B\in \R$ a lower 
bound of the real parts of the eigenvalues of $L_{M}(0)$.
From Subsection 4.2 of \cite{H-affine-twisted-mod}, we have a 
lower-bounded generalized $V$-module $\arc{M}_{B}^{1_{V}}$ satisfying 
the universal property given by Theorem 4.7 in \cite{H-affine-twisted-mod}. 
For simplicity, we shall denote it simply by 
$\arc{M}$. 

Using the construction in Subsection 4.2 of \cite{H-affine-twisted-mod}
and the results in \cite{H-const-twisted-mod} and \cite{H-exist-twisted-mod},
we see that $\arc{M}$ is generated by $M$. In particular, we identify
$M$ as a subspace of $\arc{M}$. Let $J_{M}$  be  the generalized $V$-submodule 
of $\arc{M}$ generated by elements of the forms
\begin{equation}\label{J-M-0}
\res_{x}x^{l-k-1}Y_{\arc{M}}(x^{L_{V}(0)}v, x)w
\end{equation}
for $ l=0, \dots, N$, $k\in -\Z_{+}$ and $w\in G_{l}(M)$,
\begin{equation}\label{J-M}
\res_{x}x^{l-k-1}Y_{\arc{M}}(x^{L_{V}(0)}v, x)w-\vartheta_{M}([v]_{kl})w
\end{equation}
for $v\in V$, $k, l=0, \dots, N$ and $w\in G_{l}(M)$. 

Let $S^{N}_{\rm voa}(M)=\arc{M}/J_{M}$. 
Then $S^{N}_{\rm voa}(M)$ is a lower-bounded generalized $V$-module. 
From (\ref{J-M-0}) and (\ref{J-M}), we see that $J_{M}\cap M=0$. Then
we can also identify $M$ as a subspace of $S^{N}_{\rm voa}(M)$.
Since $\arc{M}$ is generated by $M$, 
$S^{N}_{\rm voa}(M)$ is also generated by $M$ and hence is spanned by elements of the form 
\begin{equation}\label{S-N-M}
\res_{x}x^{l-n-1}Y_{S^{N}_{\rm voa}(M)}(x^{L_{V}(0)}v, x)w
\end{equation}
for $v\in V$, $l=0, \dots, N$, $n\in \N$ and $w\in G_{l}(M)$. 
In Section 5 of \cite{H-aa-va}, the nondegenerate graded $A^{\infty}(V)$-module $Gr(S^{N}(M))$ is studied.
What we are interested here is the graded $A^{\infty}(V)$-module structure on $S^{N}_{\rm voa}(M)$
and the graded $A^{N}(V)$-module structure on $\Omega_{N}^{0}(S^{N}_{\rm voa}(M))$.

\begin{prop}
For a graded $A^{N}(V)$-module $M$, the graded $A^{N}(V)$-module
$\Omega_{N}^{0}(S^{N}_{\rm voa}(M))$ is equal to $M$.
\end{prop}
\pf 
We prove only the case that the eigenvalues of $L_{M}(0)$ are all congruent to each other modulo $\Z$. 
The general case can be obtained by taking direct sums. In this case, there exists $\mu\in \C/\Z$ and 
$h^{\mu}\in \C$ such that $M=\coprod_{n=0}^{N}M_{[h^{\mu}+l]}$, where for 
$l=0, \dots, N$, $M_{[h^{\mu}+l]}$
is the generalized eigenspace of $L_{M}(0)$ with eigenvalue $h^{\mu}+l$. 

Since $S^{N}_{\rm voa}(M)$ is spanned by elements of the form (\ref{S-N-M}),
$(S^{N}_{\rm voa}(M))_{\l n\r}$ for $n\in \N$ is spanned by elements of the form
(\ref{S-N-M}) for $v\in V$, $l=0, \dots, N$ and $w\in G_{l}(M)$. But from the definition of 
$J_{M}$, when $0\le n\le N$, 
an element of this form is equal to $\vartheta_{M}([v]_{nl})w\in M_{[h^{\mu}+n]}$.
So $(S^{N}_{\rm voa}(M))_{\l n\r}\subset M$ for $n=0, \dots, N$. 
In particular, we have $\Omega_{N}^{0}(S^{N}_{\rm voa}(M))\subset M$. But by the construction, 
$M\subset \Omega_{N}^{0}(S^{N}_{\rm voa}(M))$. Thus we obtain 
$\Omega_{N}^{0}(S^{N}_{\rm voa}(M))=M$.
\epfv

From the construction 
of $S^{N}_{\rm voa}(M)$, we have the following universal property:

\begin{prop}\label{univ-prop}
Let $W$ be a lower-bounded generalized $V$-module and $f: M\to \Omega_{N}^{0}(W)$
be a graded $A^{N}(V)$-module map. Then there exists a unique $V$-module map 
$f^{\vee}: S^{N}_{\rm voa}(M)\to W$ such that $f^{\vee}|_{M}=f$. 
\end{prop}
\pf
By Theorem 4.7 in  \cite{H-affine-twisted-mod}, there is a unique $V$-module map 
$\arc{f}: \arc{M}\to W$ such that $\arc{f}|_{M}=f$. 
Since $\arc{f}$ is a $V$-module map, the image of $f$ is in $\Omega_{N}^{0}(W)$
and $f$ is a graded $A^{N}(V)$-module map, 
$$\arc{f}(\res_{x}x^{l-n-1}Y_{\arc{M}}(x^{L_{V}(0)}v, x)w)
=\res_{x}x^{l-n-1}Y_{W}(x^{L_{V}(0)}v, x)f(w)=0$$
for $ l=0, \dots, N$, $k\in -\Z_{+}$ and $w\in G_{l}(M)$ and 
\begin{align*}
&\arc{f}(\res_{x}x^{l-k-1}Y_{\arc{M}}(x^{L_{V}(0)}v, x)w-\vartheta_{M}([v]_{kl})w)\nn
&\quad =\res_{x}x^{l-k-1}Y_{W}(x^{L_{V}(0)}v, x)f(w)-\vartheta_{W}([v]_{kl})f(w)\nn
&\quad=0
\end{align*}
for $v\in V$, $k, l=0, \dots, N$ and $w\in G_{l}(M)$. So we obtain $J_{M}\subset\ker \arc{f}$.
In particular, $\arc{f}$ induces a $V$-module map 
$f^{\vee}: S^{N}_{\rm voa}(M)\to W$ such that $f^{\vee}|_{M}=f$. The uniqueness of $f^{\vee}$ follows from 
the uniqueness of $\arc{f}$. 
\epfv

Using Propositions \ref{cat-isom} and \ref{univ-prop}, we have the following consequence:

\begin{cor}\label{univ-prop-A-infty}
Let $W$ be an $A^{\infty}(V)$-module. Let $f: M\to W$
be an $A^{N}(V)$-module map with $W$ viewed as an $A^{N}(V)$-module. 
Then there exists an $A^{\infty}(V)$-module map 
$f^{\vee}: S^{N}_{\rm voa}(M)\to W$  such that $f^{\vee}|_{M}=f$. 
\end{cor}
\pf
Since $V$ is a vertex operator algebra, we have a conformal vector $\omega$ such that 
$L_{S^{N}_{\rm voa}(M)}(0)$ and $L_{S^{N}_{\rm voa}(M)}(-1)$ are given by the actions of 
$\omega^{\infty}(0)$ and $\omega^{\infty}(-1)$ on $S^{N}_{\rm voa}(M)$ (see Remark \ref{L(0)-L(-1)-comm}). 
The actions of 
$\omega^{\infty}(0)$ and $\omega^{\infty}(-1)$ on $W$ also give operators $L_{W}(0)$ and $L_{W}(-1)$
on $W$. In particular, $L_{W}(0)$ and $L_{W}(-1)$ acts on $M$ and on $\coprod_{n=0}^{N-1}G_{n}(M)$,
respectively (recall from \cite{H-aa-va} that $G_{n}(M)$ is the homogeneous subspace of $M$ of level $n$).
Since $f$ is an $A^{N}(V)$-module map, it in particular commutes with the actions of 
$\omega^{\infty}(0)$ and $\omega^{\infty}(-1)$. Then $f(M)$ is also a graded 
$A^{\infty}(V)$-module. Let $W_{0}$ be the $A^{\infty}(V)$-submodule of $W$ generated by 
$f(M)$. Then $W_{0}$ is a graded $A^{\infty}(V)$-module containing $f(M)$.

We now have a graded $A^{\infty}(V)$-module $W_{0}$ and an $A^{N}(V)$-module map 
from $M$ to $f(M)\subset W_{0}$. By Propositions \ref{cat-isom} and \ref{univ-prop}, 
there exists a unique $A^{\infty}(V)$-module map  $f^{\vee}: S^{N}_{\rm voa}(M)\to W_{0}$
such that $f^{\vee}|_{M}=f$. Since $W_{0}\subset W$, we can also view $f^{\vee}$ 
as an $A^{\infty}(V)$-module map  from $S^{N}_{\rm voa}(M)$ to $W$. 
\epfv

We also need  a right $A^{\infty}(V)$-module structure on $S^{N}_{\rm voa}(M)$ 
and a universal property of this right $A^{\infty}(V)$-module structure. These can be obtained 
using the following result and Corollary \ref{univ-prop-A-infty}:

\begin{prop}\label{opp-alg}
The linear map $O: U^{\infty}(V)\to U^{\infty}(V)$ given by 
$$O([v]_{kl})= [-e^{L_{V}(1)}(-1)^{-L_{V}(0)}v]_{kl}$$
for $k, l\in \N$ and $v\in V$ induces an isomorphism from $A^{\infty}(V)$ to its opposite algebra.
\end{prop}
\pf
We need only prove 
$$O(\mathfrak{u}\diamond \mathfrak{v})-
O(\mathfrak{v})\diamond O(\mathfrak{u})\in Q^{\infty}(V)$$
for $\mathfrak{u}, \mathfrak{v}\in U^{\infty}(V)$.

Let $W$ be a lower-bounded generalized $V$-module. Then 
\begin{align}\label{opp-alg-1}
&\langle \vartheta_{W'}([v]_{kl})w', w\rangle\nn
&\quad =\res_{x}x^{l-k-1}\langle Y_{W'}(x^{L_{V}(0)}v, x)w', w\rangle\nn
&\quad=\res_{x}x^{l-k-1}\langle w', Y_{W}(e^{xL_{V}(1)}(-x^{-2})^{L_{V}(0)}x^{L_{V}(0)}v, x^{-1})w\rangle\nn
&\quad=\res_{x}x^{l-k-1}\langle w', Y_{W}(e^{xL_{V}(1)}x^{-L_{V}(0)}(-1)^{-L_{V}(0)}v, x^{-1})w\rangle\nn
&\quad=\res_{x}x^{l-k-1}\langle w', Y_{W}(x^{-L_{V}(0)}e^{L_{V}(1)}(-1)^{-L_{V}(0)}v, x^{-1})w\rangle\nn
&\quad=-\res_{y}y^{k-l-1}\langle w, Y_{W}(y^{L_{V}(0)}e^{L_{V}(1)}(-1)^{-L_{V}(0)}v, y)w\rangle\nn
&\quad=\langle  w', \vartheta_{W}([-e^{L_{V}(1)}(-1)^{-L_{V}(0)}v]_{lk})w\rangle\nn
&\quad=\langle  w', \vartheta_{W}(O([v]_{kl}))w\rangle
\end{align}
for $k, l\in \N$, $v\in V$, $w\in W$ and $w'\in W'$. From 
(\ref{opp-alg-1}), we obtain 
\begin{equation}\label{opp-alg-2}
\langle  w', \vartheta_{W}(O(\mathfrak{v}))w\rangle=\langle \vartheta_{W'}(\mathfrak{v})w', w\rangle
\end{equation}
for $\mathfrak{v}\in U^{\infty}(V)$, $w\in W$ and $w'\in W'$.

Using (\ref{opp-alg-2}), we have 
\begin{align}\label{opp-alg-3}
&\langle  w', \vartheta_{W}(O(\mathfrak{u}\diamond \mathfrak{v}))w\rangle\nn
&\quad =\langle \vartheta_{W'}(\mathfrak{u}\diamond \mathfrak{v})w', w\rangle\nn
&\quad =\langle \vartheta_{W'}(\mathfrak{u})\vartheta_{W'}(\mathfrak{v})w', w\rangle\nn
&\quad =\langle \vartheta_{W'}(\mathfrak{v})w', \vartheta_{W}(O(\mathfrak{u}))w\rangle\nn
&\quad =\langle w', \vartheta_{W}(O(\mathfrak{v}))\vartheta_{W}(O(\mathfrak{u}))w\rangle\nn
&\quad =\langle w', \vartheta_{W}(O(\mathfrak{v})\diamond O(\mathfrak{u}))w\rangle
\end{align}
for $\mathfrak{u}, \mathfrak{v}\in U^{\infty}(V)$, $w\in W$ and $w'\in W'$. Since $w'$ and $w$ 
are arbitrary, 
we obtain from (\ref{opp-alg-3})
$$O(\mathfrak{u}\diamond \mathfrak{v})-O(\mathfrak{v})\diamond O(\mathfrak{u})\in \ker \vartheta_{W}$$
for $\mathfrak{u}, \mathfrak{v}\in U^{\infty}(V)$. Since $W$ is arbitrary, 
we obtain
$$O(\mathfrak{u}\diamond \mathfrak{v})-O(\mathfrak{v})\diamond O(\mathfrak{u})\in Q^{\infty}(V),$$
proving that $O$ indeed induces an isomorphism from $A^{\infty}(V)$ to its opposite algebra. 
\epfv

Since $A^{\infty}(V)$ is isomorphic to its opposite algebra, a left $A^{\infty}(V)$-module also give
a right $A^{\infty}(V)$-module and vice versa. In particular, $S^{N}_{\rm voa}(M)$ is also a right 
$A^{\infty}(V)$-module. Then from Corollary \ref{univ-prop-A-infty}, we obtain the following
consequence immediately:

\begin{cor}\label{univ-prop-A-infty-right}
Let $W$ be a right $A^{\infty}(V)$-module. Let $f: M\to W$
be a right $A^{N}(V)$-module map when $W$ is viewed as a right $A^{N}(V)$-module. 
Then there exists a right $A^{\infty}(V)$-module map 
$f^{\vee}: S^{N}_{\rm voa}(M)\to W$  such that $f^{\vee}|_{M}=f$. 
\end{cor}

We now prove our second main theorem. 

\begin{thm}\label{main-N}
Let $V$ be a vertex operator algebra. 
Assume that $W_{2}$ and $W_{3}'$ are equivalent to 
 $S^{N}_{\rm voa}(\Omega_{N}^{0}(W_{2}))$ and $S^{N}_{\rm voa}(\Omega_{N}^{0}(W_{3}'))$,
respectively.
Then $\rho^{N}$ is a linear isomorphism.
\end{thm}
\pf
By Theorem \ref{main}, $\rho$ is a linear isomorphism.
Then the composition $\rho^{N} \circ \rho^{-1}$ is a linear map 
from 
$$\hom_{A^{\infty}(V)}(A^{\infty}(W_{1})\otimes_{A^{\infty}(V)}W_{2}, 
W_{3})$$ 
to 
$$\hom_{A^{N}(V)}(A^{N}(W_{1})\otimes_{A^{N}(V)} \Omega_{N}^{0}(W_{2}), 
\Omega_{N}^{0}(W_{3})).$$
We need only prove that  $\rho^{N} \circ \rho^{-1}$ is an isomorphism. 
Then from the definition of $\rho^{N}$ and Theorem \ref{main},
$\rho^{N}$ is also a linear isomorphism.

By Proposition \ref{rho-N-inj}, $\rho^{N}$ is injective. 
Since both $\rho^{N}$ and $\rho^{-1}$ are injective, $\rho^{N} \circ \rho^{-1}$
is also injective. We still need to prove that $\rho^{N} \circ \rho^{-1}$ is surjective. 
Let 
$$f^{N}\in \hom_{A^{N}(V)}(A^{N}(W_{1})\otimes_{A^{N}(V)} \Omega_{N}^{0}(W_{2}), 
\Omega_{N}^{0}(W_{3})).$$
We want to find 
$$f\in \hom_{A^{\infty}(V)}(A^{\infty}(W_{1})\otimes_{A^{\infty}(V)} W_{2}, 
W_{3})$$
such that the $(\rho^{N} \circ \rho^{-1})(f)$ 
is $f^{N}$.

Let $A^{N,\infty}(W_{1})$ be the subspace of $A^{\infty}(W_{1})$
obtained by taking sums of elements of the form 
$[w_{1}]_{np}+Q^{\infty}(W_{1})$
for $w_{1}\in W_{1}$, $0\le n\le N$ and $p\in \N$,  including certain infinite sums  
as in the case of $U^{\infty}(W_{1})$. Then $A^{N,\infty}(W_{1})$ is an
$A^{N}(V)$-$A^{\infty}(V)$-bimodule. 
Since 
$\Omega_{N}^{0}(W_{3}')\otimes_{A^{N}(V)} A^{N, \infty}(W_{1})$
is a right $A^{\infty}(V)$-module, its dual space
$(\Omega_{N}^{0}(W_{3}')\otimes_{A^{N}(V)} A^{N, \infty}(W_{1}))^{*}$
a left $A^{\infty}(V)$-module. The map $f^{N}$ gives an $A^{N}(V)$-module map from 
$\Omega_{N}^{0}(W_{2})$ to $(\Omega_{N}^{0}(W_{3}')\otimes_{A^{N}(V)} A^{N, \infty}(W_{1}))^{*}$
as follows: For $w_{2}\in \Omega_{N}^{0}(W_{2})$,
we define the image of $w_{2}$ in $(\Omega_{N}^{0}(W_{3}')\otimes_{A^{N}(V)} A^{N, \infty}(W_{1}))^{*}$
to be the linear functional given by 
$$w_{3}'\otimes_{A^{N}(V)}([w_{1}]_{np}+Q^{\infty}(W_{1}))\mapsto 
\langle w_{3}', f^{N}(([w_{1}]_{np}+Q^{\infty}(W_{1}))\otimes_{A^{N}(V)}w_{2})\rangle$$ 
for $0\le n, p\le N$, $w_{1}\in W_{1}$ and $w_{3}'\in \Omega_{N}^{0}(W_{3}')$
and 
$$w_{3}'\otimes_{A^{N}(V)}([w_{1}]_{np}+Q^{\infty}(W_{1}))\mapsto 0$$
for $0\le n\le N$, $p\in N+1+\N$, $w_{1}\in W_{1}$ and $w_{3}'\in \Omega_{N}^{0}(W_{3}')$.
For $0\le k, l\le N$, $v\in V$, we have 
\begin{align*}
&\langle w_{3}', f^{N}(([w_{1}]_{np}+Q^{\infty}(W_{1}))\otimes_{A^{N}(V)}\vartheta_{W_{2}}([v]_{kl})w_{2})\rangle\nn
&\quad =\langle w_{3}', f^{N}(([w_{1}]_{np}+Q^{\infty}(W_{1}))\diamond ( [v]_{kl}+Q^{\infty}(V))
\otimes_{A^{N}(V)}w_{2})\rangle.
\end{align*}
This means that the linear map from 
$\Omega_{N}^{0}(W_{2})$ to $(\Omega_{N}^{0}(W_{3}')\otimes_{A^{N}(V)} A^{N, \infty}(W_{1}))^{*}$
is an $A^{N}(V)$-module map. 

Since $W_{2}$ is equivalent to $S^{N}_{\rm voa}(\Omega_{N}^{0}(W_{2}))$, 
by Corollary \ref{univ-prop-A-infty}, there exists a unique $A^{\infty}(V)$-module map from 
$W_{2}$ to $(\Omega_{N}^{0}(W_{3}')\otimes_{A^{N}(V)} A^{N, \infty}(W_{1}))^{*}$
such that when restricted to $\Omega_{N}^{0}(W_{2})$, it is equal to 
the $A^{N}(V)$-module map given above. But such an 
$A^{\infty}(V)$-module map is equivalent to an element
$$f^{N, \infty}\in (\Omega_{N}^{0}(W_{3}')\otimes_{A^{N}(V)}A^{N, \infty}(W_{1})\otimes_{A^{\infty}(V)}W_{2})^{*}.$$

Since
$A^{\infty}(W_{1})\otimes_{A^{\infty}(V)}W_{2}$
is a left $A^{\infty}(V)$-module, its dual space
$(A^{\infty}(W_{1})\otimes_{A^{\infty}(V)}W_{2})^{*}$
is a right $A^{\infty}(V)$-module. 
The map $f^{N, \infty}$ gives a right $A^{N}(V)$-module map from 
$\Omega_{N}^{0}(W_{3}')$ to $(A^{N, \infty}(W_{1})\otimes_{A^{\infty}(V)}W_{2})^{*}$
as follows: For $w_{3}'\in \Omega_{N}^{0}(W_{3}')$,
we define the image of $w_{3}'$ in $(A^{N, \infty}(W_{1})\otimes_{A^{\infty}(V)}W_{2})^{*}$
to be the linear functional given by 
$$([w_{1}]_{np}+Q^{\infty}(W_{1}))\otimes_{A^{\infty}(V)}w_{2}\mapsto 
f^{N, \infty}(w_{3}' \otimes_{A^{N}(V)}([w_{1}]_{np}+Q^{\infty}(W_{1}))\otimes_{A^{N}(V)}w_{2})\rangle$$ 
for $0\le n\le N$, $p\in \N$, $w_{1}\in W_{1}$ and $w_{2}\in W_{2}$
and 
$$([w_{1}]_{np}+Q^{\infty}(W_{1}))\otimes_{A^{\infty}(V)}w_{2}\mapsto 0$$
for $n\in N+1+\N$, $p\in \N$, $w_{1}\in W_{1}$ and $w_{2}\in W_{2}$.
We use $w_{3}'\vartheta_{W_{3}'}^{r}(\mathfrak{v})$ to denote the right action of 
$\mathfrak{v}$ on $W_{3}'$. Then for $0\le k, l\le N$, $v\in V$, we have 
\begin{align*}
&f^{N, \infty}(w_{3}' \vartheta_{W_{3}'}^{r}([v]_{kl})
\otimes_{A^{N}(V)}([w_{1}]_{np}+Q^{\infty}(W_{1}))\otimes_{A^{N}(V)}w_{2})\rangle\nn
&\quad =f^{N, \infty}(w_{3}'
\otimes_{A^{N}(V)}([v]_{kl}\diamond [w_{1}]_{np}+Q^{\infty}(W_{1}))\otimes_{A^{N}(V)}w_{2})\rangle.
\end{align*}
This means that the linear map from 
$\Omega_{N}^{0}(W_{3}')$ to $(A^{N, \infty}(W_{1})\otimes_{A^{\infty}(V)}W_{2})^{*}$
is a right $A^{N}(V)$-module map. 

Since $W_{3}'$ is equivalent to $S^{N}_{\rm voa}(\Omega_{N}^{0}(W_{3}'))$, 
by Corollary \ref{univ-prop-A-infty-right}, there exists a unique right $A^{\infty}(V)$-module map from 
$W_{3}'$ to $(A^{\infty}(W_{1})\otimes_{A^{\infty}(V)}W_{2})^{*}$
such that when restricted to $\Omega_{N}^{0}(W_{3}')$, it is equal to 
the right $A^{N}(V)$-module map given above. 
But such a right
$A^{\infty}(V)$-module map is equivalent to a left $A^{\infty}(V)$-module map 
from $(A^{\infty}(W_{1})\otimes_{A^{\infty}(V)}W_{2})^{**}$ to $W_{3}$.
Since $A^{\infty}(W_{1})\otimes_{A^{\infty}(V)}W_{2}$ is a $A^{\infty}(V)$-submodule of
$(A^{\infty}(W_{1})\otimes_{A^{\infty}(V)}W_{2})^{**}$, we obtain a unique 
left $A^{\infty}(V)$-module map  from $A^{\infty}(W_{1})\otimes_{A^{\infty}(V)}W_{2}$ to 
$W_{3}$, or equivalently, an element
$$f\in \hom_{A^{\infty}(V)}(A^{\infty}(W_{1})\otimes_{A^{\infty}(V)}W_{2}, W_{3}).$$
From our construction of $f$, it is clear that $(\rho^{N} \circ \rho^{-1})(f)=f^{N}$. 
This proves the surjectivity of $\rho^{N} \circ \rho^{-1}$.
\epfv

\noindent {\small \sc Department of Mathematics, Rutgers University,
110 Frelinghuysen Rd., Piscataway, NJ 08854-8019}

\noindent {\em E-mail address}: yzhuang@math.rutgers.edu

\end{document}